\documentclass[11pt]{article}
\usepackage{amsmath,amssymb,amsthm,graphicx,amsfonts,bbm}
\usepackage{subfig}

\theoremstyle{plain}
\begingroup
 
\newtheorem{theorem}{Theorem}

\newtheorem{remark}{Remark}
\endgroup

\newcommand{\R}{\mathbb{R}}

\newcommand{\tr}{\text{tr}}

\newcommand{\pr}{\text{pr}}

\newcommand{\one}{\mathbbm{1}}

\newcommand{\Jac}{\mathrm{Jac}}

\newcommand{\sign}{\text{sign}}

\newcommand{\eps}{\epsilon}
\newcommand{\be}{\begin{equation}}
\newcommand{\ee}{\end{equation}}
\newcommand{\ba}{\begin{eqnarray}}
\newcommand{\ea}{\end{eqnarray}}
\newcommand{\qref}[1]{(\ref{#1})}
\newcommand{\deftime}{\tau}
\newcommand{\contdeftime}{\tau}
\newcommand{\resctime}{T}
\newcommand{\defindex}{\iota}
\newcommand{\Symm}{\mathrm{Symm}}
\newcommand{\specmap}{\mathcal{S}}
\newcommand{\diag}{\mathrm{diag}}

\title{How long does it take to compute the eigenvalues of a random symmetric matrix?\thanks{The work of the first and third author was supported in part by  NSF grant DMS 07-48482. The work of the second author was supported in part by NSF grant DMS 10-01886. The second author
also acknowledges support by a grant from The Simonyi Fund at the Institute for Advanced Study in Princeton}}

\author{Christian W. Pfrang \thanks{J.P. Morgan Securities LLC, 383 Madison Ave., New York, NY 10033, USA,  ({\tt christian.w.pfrang@jpmorgan.com})},
\and Percy Deift \thanks{Courant Institute for Mathematical Sciences, New York University, 251 Mercer St., New York, NY, 10012, USA,   ({\tt deift@cims.nyu.edu})},
\and Govind Menon \thanks{Division of Applied Mathematics, Brown University, 182 George St., Providence, RI 02912, USA, ({\tt menon@dam.brown.edu}).}}

\newtheorem{thm:dlnt}{Theorem}
\newtheorem{thm:discrete-signum-algo}[thm:dlnt]{Theorem}

\begin{document}
\maketitle


\begin{abstract}
We present the results of an empirical study of the performance of the QR 
algorithm (with and without shifts) and the Toda algorithm
on random symmetric matrices. The random matrices are chosen from six ensembles, four of which lie in the Wigner class.  For all three algorithms, we observe a form of universality for the deflation time  statistics for random matrices within the Wigner class. For these ensembles, the empirical distribution of a normalized deflation time is found to collapse onto a curve that depends only on the algorithm, but not on the matrix size or deflation tolerance provided the matrix size is large enough (see Figure~\ref{fig:QRnormalized-runtimes-wigner}, Figure~\ref{fig:S-QRnormalized-runtimes-wigner} and Figure~\ref{fig:Todanormalized-runtimes-wigner}). For the QR algorithm with the Wilkinson shift, the observed universality is even stronger and includes certain non-Wigner ensembles. Our experiments also provide a quantitative statistical picture of the accelerated convergence with shifts.
\end{abstract}




\section{Introduction}
\label{sec:intro}
We present the results of a statistical study of the performance of the QR and Toda eigenvalue algorithms on random symmetric matrices. Our work is mainly inspired by progress in quantifying the ``probability of difficulty'' and ``typical behavior'' for several numerical algorithms~\cite{DemmelDifficult,vonNeumann}. This approach has led to a deeper understanding of the efficacy of fundamental numerical algorithms such as Gaussian elimination and the simplex method~\cite{Rudelson-Vershynin, SmoothedAnalysisCondition,Smale,Tao-Vu-2010}. It has also stimulated new ideas in random matrix theory~\cite{MatrixModelsBetaEnsembles,RandMatCondition,Edelman-Sutton}. Testing eigenvalue algorithms with random input continues this effort. In related work~\cite{ChristianDiss}, we have also studied the performance of a version of the matrix sign algorithm. However, these results are of a different character, and apart from some theoretical observations, we do not present any experimental results for this algorithm (see~\cite{ChristianDiss} for more information). Our study is empirical -- a study of the eigenvalue problem from the viewpoint of complexity theory is presented in~\cite{Armentano}.

\subsection{Algorithms and Ensembles}
It is natural to study the QR algorithm because of its elegance and fundamental practical importance. But in fact all the algorithms we study are linked by a common framework. In each case, an initial matrix $L_0$ is diagonalized via a sequence of isospectral iterates $L_m$. The gist of the framework is that the $L_m$ correspond exactly to the flow of a completely integrable Hamiltonian system evaluated at integer times. The Hamiltonian for these flows is of the form $\tr \,G(L)$ where $G$ is a real-valued function defined on an interval. Different choices of $G$ generate different algorithms: $G(x)= x(\log x -1)$ yields unshifted QR, $G(x)= x^2/2$ yields Toda, and $G(x) = |x|$ yields the matrix sign algorithm. As noted above, we will not present any numerical experiments on the matrix sign algorithm (but see Section~\ref{sec:background}). We note that the practical implementation of the QR algorithm requires an efficient shifting strategy. Our work includes a study of the QR algorithm with the Wilkinson shift as discussed below. 

Initial matrices are drawn from six ensembles that  arise  in random matrix theory. These are listed below in Section~\ref{subsec:RMT}. For many random matrix ensembles,  as the size of the matrix grows, the density of eigenvalues and suitably rescaled fluctuations have limiting distributions that may be computed explicitly. Four of the ensembles we study consist of random matrices with independent entries subject to the constraint of symmetry. The law of these entries is chosen so that these ensembles have the Wigner semicircle law as limiting spectral density. We say that these ensembles are in the {\em Wigner class\/}. Numerical experiments with these ensembles are contrasted with  two ensembles that do not belong to the Wigner class.  

\subsection{Deflation and QR with the Wilkinson shift}
In evaluating these algorithms we focus on the statistics of deflation. Given a real, symmetric, $n \times n$ matrix $L$ and an integer $k$ between $1$ and $n$, we write 
\be
\label{eq:pert1}
L = \left( \begin{array}{ll} L_{11} & L_{12} \\ L_{12}^T & L_{22} \end{array} \right),\quad \tilde{L}= \left( \begin{array}{ll} L_{11} & 0 \\ 0 & L_{22} \end{array} \right),
\ee
where $L_{11}$ is a  $k \times k$ block. Let $\lambda_j$ and $\tilde{\lambda}_j$, $j=1,\ldots, n$, denote the eigenvalues of $L$ and $\tilde{L}$. For a fixed tolerance $\eps>0$ we say that $L$ is deflated to $\tilde{L}$ when the off-diagonal block $L_{12}$ is so small that $\max_j |\lambda_j -\tilde{\lambda}_j| < \eps$.  The {\em deflation time\/} is the number of iterations $m$ before $L_m$ can be deflated by a tolerance $\eps>0$ at some index $k$. The {\em deflation index\/} is this value of $k$.  Since the iterative eigenvalue algorithms correspond to Hamiltonian flows, there is also a natural notion of deflation time for the Hamiltonian flows (see equations~\qref{eq:runtime} and \qref{eq:cont-runtime} below).

Let us now explain why deflation serves as a useful measure of the time required to compute the eigenvalues of a matrix. The cost of practical computation requires an analysis of algorithms, hardware and software. In our study, we only focus on the algorithm, and ``time'' is taken to mean the number of iterations required for convergence.   In our experiments we have observed that the QR and Toda algorithms deflate a matrix at the upper-left or lower-right corner with high probability. The deflation index for the shifted QR algorithm is $n-1$ with overwhelming probability. The deflation index for unshifted QR is also typically $n-1$ (see  Figure~\ref{fig:QR-deflation-indices} and Figure~\ref{fig:Toda-deflation-indices}). As a consequence, the deflation time is typically the same as the time taken to compute an eigenvalue. We then expect that the time taken to compute all eigenvalues with these algorithms is determined by $n$ deflations. By contrast, we find that the matrix sign algorithm typically deflates a matrix in the middle and  does not immediately yield any eigenvalues. Instead, these are obtained after a divide-and-conquer procedure that consists of approximately  $\log_2 n$ deflations. Thus, for all these algorithms a finite sequence of deflation times determines the number of iterations necessary to compute eigenvalues. We must note however, that we do not track all deflations in our experiments, only the first. This restriction is necessary to keep the datasets manageable as $n$ increases. A more extensive study that tracks all deflation times for these algorithms will certainly yield further interesting information. Finally, as we show in Section~\ref{subsec:example} below, the notion of deflation time is also of theoretical value since it is the starting point for an analysis of the expected number of iterations for eigenvalue algorithms that is similar in spirit to~\cite{Smale}.

The convergence of the QR algorithm is greatly accelerated by shifts. We will only consider the Wilkinson shift, i.e. the shift  is the eigenvalue of the $2\times 2$ lower diagonal corner of the matrix that is closer to $L_{nn}$.  The QR algorithm on tridiagonal matrices is cubically convergent with this choice of shift (this is generically true~\cite{Wilkinson}, see also~\cite{Tomei} for a more careful analysis). As noted above, the unshifted QR algorithm deflates at index $n-1$ with very high probability. Since the Wilkinson shift utilizes the lower $2\times 2$ block of the matrix, the number of the iterations required for shifted QR, as opposed to unshifted QR, to deflate is far smaller. While such acceleration of convergence is well-known, some features of our experiments still come as a surprise.  For example, a striking feature of Figure~\ref{fig:goe-qr} and Figure~\ref{fig:goe-qrw} is that the number of iterations required to deflate a random matrix with the QR algorithm (shifted and unshifted)  is almost independent of $n$ for matrices as large as $190\times 190$.

\subsection{Universality}
Our main empirical findings concern universal fluctuations in the deflation time distribution for the QR algorithm (shifted and unshifted) and the Toda algorithm for ensembles in the  Wigner class.  We sample the deflation time for a range of matrix size and deflation tolerance combinations and normalize these empirical distributions to  mean zero and variance one. The resulting histograms have the same general shape and in particular, the same tails on the right side (see in particular  Figure~\ref{fig:QRnormalized-runtimes-wigner}, Figure~\ref{fig:S-QRnormalized-runtimes-wigner} and Figure~\ref{fig:Todanormalized-runtimes-wigner}). In other words, {\em the fluctuations in deflation time are universal\/}. For the Toda and unshifted QR algorithm, the observed limiting fluctuations for Wigner and non-Wigner ensembles are distinct (see Figure~\ref{fig:combined-histograms} and Figure~\ref{fig:S-combined-histograms}). In addition, we find that the universal distributions for Wigner ensembles have exponential tails for unshifted QR and Gaussian tails for Toda (Figure~\ref{fig:combined-histograms} and Figure~\ref{fig:Toda-combined-histograms}). Universality of the tails is quantified with a statistical methodology developed in \cite{powerlawEmpiricalData}.  Quite remarkably, for the (Wilkinson) shifted QR algorithm, the observed universality is stronger: to a good approximation {\em all\/} tested ensembles show the same limiting distribution (see Figure \ref{fig:S-combined-histograms}).

The origin of such universality is not clear. We do not understand fully if our results are connected with the now familiar universality theorems of random matrix theory such as those that describe fluctuations in the bulk and at the edge of the spectrum for the Wigner ensembles~\cite{Yau,MehtaRMT,Tracy-Widom}. Unlike these universality theorems, where the mean and variance are known theoretically,  in our work the mean and variance of the deflation time are computed empirically and we have not yet been able to determine analytically how these depend on $n$. It does appear however that the mean deflation time is linearly proportional to $\log \eps$ (see Figure~\ref{fig:QR-means} and Figure~\ref{fig:Toda-means}). 

More broadly, our experiments are suggestive of a wider class of questions concerning universality of fluctuations for computations in numerical linear algebra. For example, in similar experiments to be reported elsewhere, one of the authors (P.D.) and Sheehan Olver have studied the solution $x$ to the linear equation $Ax=b$ empirically,  when $A$ is a random positive symmetric matrix and $b$ is a random vector. They compute the solution using the conjugate gradient method and observe universal fluctuations in the number of iterations required for convergence, independent of the choice of ensemble for $A$ and $b$.

We now discuss the algorithms and ensembles in greater detail.  This is followed by a description of the results in Section~\ref{sec:numerical-results}. The implementation of the algorithms is discussed briefly in Section~\ref{sec:methods}.

\section{Algorithms, ensembles and deflation statistics}
\label{sec:background}
\subsection{Notation}
We denote the space of real, symmetric $n \times n$ matrices by $\Symm(n)$ and the space of real $n \times m$ matrices by $\R^{n \times m}$. Matrices in $\Symm(n)$ are denoted  $L$ or $M$ and  the iterates of an eigenvalue algorithm  are denoted $L_m$ or $M_m$, $m=0,1,2,\ldots$. We use $Q$ to denote an orthogonal matrix and $R$ an upper triangular matrix with positive diagonal entries, typically with reference to a QR factorization. We use $\sigma(L)$ to denote the spectrum of $L$. The spectral decomposition of $L \in \Symm(n)$ is written $L = U \Lambda U^T$. Here $U$ denotes the orthogonal matrix of eigenvectors of $L$ and $\Lambda =\diag (\lambda_1,\lambda_2, \ldots, \lambda_n)$ the matrix of eigenvalues. 
We also use the following standard notation. The $n\times n$ identity matrix  is $I$; the  standard basis in $\R^n$ is $(e_1, \ldots, e_n)$; the unit sphere in $\mathbb{R}^{n+1}$ and  its positive orthant are $S^n$ and $S^{n}_+$ respectively; the symmetric group of order $n$ is  $S_n$. 

When $L \in \Symm(n)$ is tridiagonal, we denote its diagonal entries by $(a_1, \ldots, a_n)$, and its off-diagonal entries by $(b_1,\ldots, b_{n-1})$. A Jacobi matrix is a tridiagonal matrix with $b_i>0$. The space of Jacobi matrices is denoted $\Jac(n)$. We use $u=U^Te_1$ to denote the first row of the matrix of eigenvectors. When $L \in \Jac(n)$, $\sigma(L)$ is simple and we may assume that $\lambda_1> \lambda_2 > \ldots > \lambda_n$. Moreover, all the components $u_i$ are non-zero and we may assume that $u_i>0$ (this is also generically true for $L \in \Symm(n)$).  
Let $\mathcal{M}$ denote the manifold $\{ (\Lambda,u) \in \R^n \times S_+^n| \lambda_1 > \lambda_2 > \ldots > \lambda_n\}$.  It is a basic result in the spectral/ inverse spectral theory for Jacobi matrices that the matrix $L$ can be reconstructed if $\Lambda$ and $u$ are  given. More precisely, the spectral map $\specmap: \Jac(n)\to \mathcal{M}$ defined by $L \mapsto (\Lambda,u)$ is a diffeomorphism~\cite[Thm 2]{ODESymmetricEig}. 

We use the following standard notation for probabilistic notions. The phrase  independent and identically distributed is abbreviated to  iid. A normal random variable with mean $\mu$ and variance $\sigma^2$ is denoted $\mathcal{N}(\mu,\sigma^2)$;  a Bernoulli random variable that is $\pm 1$ with probability $1/2$ is denoted $\mathcal{B}$; a random variable with the $\chi$-distribution with parameter $k$ is denoted $\chi_k$.  The notation $X \sim Y$ means that $X$ has the same law as $Y$.

\subsection{The QR algorithm and Hamiltonian eigenvalue algorithms}
\label{sec:dlnt}
We assume the reader is familiar with the QR algorithm (excellent textbook presentations are \cite{AppliedNumLA,GVL,Trefethen}). In the unshifted QR algorithm the iterates $M_m$ are generated through QR factorizations and matrix multiplication in the reverse order:
\begin{equation}
\label{eq:qr-basic}
Q_m R_m = M_{m}, \quad M_{m+1}=R_m Q_m, \quad m=0,1,2, \ldots
\end{equation} 
The shifted QR algorithm relies on a shift $\mu_m$ at each step, and the modified steps
\begin{equation}
\label{eq:qr-shift}	
Q_m R_m = M_m -\mu_m I, \quad M_{m+1}= R_m Q_m +\mu_m I, \quad m=0,1,2, \ldots
\end{equation}
Typical shifts, such as the Wilkinson shift, are constructed from the lower $2 \times 2$ block of $M_m$~\cite[p.418]{GVL}.

In the early 80's it was discovered that the QR algorithm is intimately connected with integrable Hamiltonian systems~\cite{HamiltonianGroup,QRAlgoScattering,TodaFlowGeneric,ODESymmetricEig,DiffEqQRAlgo}.  We summarize these results below. An expanded presentation of these connections may be found in~\cite{NotesPercy,Deift-Li-Symplectic,ChristianDiss}. A different exposition that explains these ideas in a fashion ``intrinsic'' to  numerical linear algebra is~\cite{IsospectralFlows}.

Assume $G$ is a piecewise smooth real-valued function defined on an interval, and set $g=G'$. If $g$ is defined on $\sigma(L)$, we define $g(L):=Ug(\Lambda)U^T$. Let $M_-$ denote the strictly lower
triangular part of the square matrix $M$, and $\pr_{\mathfrak{k}}M:=M_-^T-M_-$, the projection of $M$ onto skew-symmetric matrices. We then consider  the ordinary differential equation
\begin{align} \label{eq:dlnt-equation}
\dot{L} = [\pr_\mathfrak{k}g(L),L].
\end{align}
Equation \qref{eq:dlnt-equation} defines a completely integrable Hamiltonian flow on the space of (generic) symmetric matrices with Hamiltonian $H(L)=\tr \, G(L)$ and symplectic structure detailed in~\cite{TodaFlowGeneric}. This flow is connected to the unshifted QR algorithm as follows. 
\medskip

\begin{theorem}
\label{mark:thm-dlnt-big}
Let $g$ be a real-valued function defined on $\sigma(L_0)$. Then 
\begin{enumerate}
\item[(a)] The solution to equation (\ref{eq:dlnt-equation}) with initial condition $L_0$ is an isospectral deformation 
\be
\label{eq:qr-thm1}
L(t) = Q(t)^T L_0 Q(t), 
\ee
where the orthogonal matrix $Q(t)$ is given by the unique QR factorization 
\begin{equation}
\label{eq:qr-thm}
e^{tg(L_0)} =Q(t) R(t), \quad t \geq 0,
\end{equation}
that has $Q(0)=I$ and depends smoothly on $t$. 
\item[(b)] At integer times $m=0, 1,2 \ldots$ the solution $L(m)$ satisfies
\begin{equation}
	\label{eq7}
e^{g(L(m))}  = M_m,
\end{equation}
where $M_m$ is the $m$--th step of the QR algorithm applied to the initial matrix
$M_0=e^{g(L_0)}$.
\noindent
\item[(c)] Assume that the spectrum $\sigma(L_0)$ is simple and that $g$ is injective on $\sigma(L_0)$. Then  
$L_\infty =\lim_{t \to \infty} L(t)$ is a diagonal matrix consisting of the eigenvalues of $L_0$. 
\end{enumerate}
\end{theorem}
\medskip

The case of tridiagonal matrices is of practical and theoretical importance. When $L_0$ is tridiagonal, so is $L(t)$,  and the flow can be linearized using the  spectral map $\specmap$ for Jacobi matrices. 
\medskip

\begin{theorem}
\label{mark:thm-dlnt}
Assume $L_0\in \Jac(n)$. Then the solution $L(t)$ to \qref{eq:dlnt-equation}
is an isospectral deformation $L(t)=U(t)^T \Lambda U(t)$ and the evolution of $u(t)=U(t)^Te_1$ and $L(t)$ is given explicitly by  
\begin{equation} 
	\label{eq:dlnt-spectral-explicit}
u(t)  = \frac{e^{tg(\Lambda)}u_0}{\left\|e^{tg(\Lambda)}u_0\right\|}, \quad L(t) = \specmap^{-1}(\lambda,u(t)). 
\end{equation}
Assume  $g$ is injective on $\sigma(L_0)$. Then 
\begin{equation}
	\label{eq:permut-increase}
\lim_{t \to \infty} L(t) = \diag\left(\lambda_{\sigma_1}, \ldots, \lambda_{\sigma_n}\right),
\end{equation}
where $\sigma \in S_n$ is the permutation such that  $g(\lambda_{\sigma_1}) >\cdots>g(\lambda_{\sigma_n})$.
\end{theorem}

Theorem~\ref{mark:thm-dlnt-big} and Theorem~\ref{mark:thm-dlnt}  may be used
to develop numerical schemes. The main observation is that 
each choice of a Hamiltonian $H(L)= \tr \, G(L)$ corresponds to a choice of an algorithm. In particular, we have
\begin{enumerate}
\item The \textbf{unshifted QR algorithm}: $g(x)=\log(x)$, $G(x) = x(\log x -1)$ and
$H_\text{QR}(L)=\tr\,[L\log L-L]$ \cite{DiffEqQRAlgo}. 
\item The \textbf{Toda algorithm}: $g(x)=x$, $G(x)=x^2/2$ and $H_\text{Toda}(L)=\frac{1}{2}\tr \, L^2$. In this case, equation (\ref{eq:dlnt-equation}) describes the evolution of the Toda lattice\cite{MoserNotesDS}.
\item The \textbf{matrix sign algorithm}: $g(x)=\sign(x)$, $G(x)=|x|$ and $H_\text{sign}(L)=\tr\,|L|$. 
\end{enumerate}
\medskip
Of course, each step of the shifted QR algorithm, $L \mapsto L-\mu I$, is Hamiltonian, with Hamiltonian $H_{\mathrm{QR,shift}}(L) =H_{\mathrm{QR}}(L-\mu I)$.
While every function $G$ defines a Hamiltonian  not all choices are equally relevant. Since our goal is to find the spectral decomposition of $L_0$, we must assume that $U$ and $\Lambda$ are unknown. But then how are we to compute the matrix-valued functions $g(L)$ or $e^{g(L)}$ efficiently? The choices $g(x)=\log x$ and $g(x)=x$ are special since these give $e^{g(L)}=L$ and  $g(L)=L$ respectively. The first choice gives the QR algorithm (strictly speaking a branch of the logarithm must be chosen so that \qref{eq:dlnt-equation} is well-defined, but this does not affect the QR algorithm because of equation~(\ref{eq7})). 
For the second choice $g(x)=x$, the vector field \qref{eq:dlnt-equation} is faster to compute than the matrix exponential $e^{L(m)}$ and it is natural to use an ordinary differential equation solver for \qref{eq:dlnt-equation} to diagonalize $L$. This is the essence of the Toda algorithm. 

Our final choice $g(x)=\sign(x)$ requires further comment since the observation that the matrix sign algorithm is Hamiltonian seems to us to be new. Assume zero is not an eigenvalue of $L_0$ and  let $\Sigma_{\pm}$ denote the eigenspaces of $L_0$ corresponding to positive and negative eigenvalues respectively.
Consider matrices $Q_\pm$ whose columns form an orthonormal basis for $\Sigma_{\pm}$ respectively. Then the matrices $P_+=Q_+Q_+^T$ and $P_-=Q_-Q_-^T$ are orthogonal projections onto $\Sigma_\pm$ respectively and we find $\sign(L_0)= P_+- P_-$ and $\left(I \pm \sign(L_0)\right)/2 =P_\pm$. It is immediate that 
\be
\label{eq:sign-proj2} 
e^{t\,\sign(L_0)} = e^tP_+ - e^{-t} P_-, \quad\mathrm{and}\quad \lim_{t \to \infty} e^{-t} e^{t\,\sign(L_0)} =P_+. 
\ee
The projection $P_+$ has a rank-revealing QR factorization $P_+ = U_\infty R_\infty \Pi$~\cite[Ch. 2.5]{Higham}. The matrix sign algorithm rests on the fact that with $U_\infty$ as above, $U_\infty^T U_\infty=I$, and the matrix
\be
\label{eq:matrix-sign-new}
\tilde{L}=U_\infty^T L_0 U_\infty
\ee
is block-diagonal as in equation~\qref{eq:pert1}, where $L_{11}$ is $k\times k$ with $k=\mathrm{dim}(\Sigma_+)$. Clearly, $\sigma(\tilde{L})=\sigma(L_0)$.  

Thus, the procedure to deflate a matrix using the matrix sign algorithm is:
\begin{enumerate}
\item Given $L_0$, compute $\sign(L_0)$ and hence $P_+=\left(I + \sign(L_0)\right)/2$. 
\item Compute $U_\infty$ using a rank-revealing QR decomposition of $P_+$.
\item Compute $\tilde{L} = U_\infty^T L_0 U_\infty$.
\end{enumerate}
We note that $\sign(L_0)$ can be computed efficiently using a scaled Newton iteration and inverse-free modifications of this procedure~\cite{Bai-Demmel-Gu,Higham, Malyshev}. The complete spectral decomposition of $L_0$ may be determined in a sequence of deflation steps. At each stage, the number of iterations required to deflate the matrix depends on the number of iterations required to compute  $\sign(L_0)$. 

From the dynamical point of view, let $L(t)$ denote the solution to \qref{eq:dlnt-equation} with $g(L)=\sign(L)$. Then it may be shown that for generic initial data $\Pi=I$ and $\lim_{t \to \infty} L(t) =\tilde{L}$ where $\tilde{L}=U_\infty^TL_0 U_\infty$ is the block-diagonal matrix obtained above by the matrix sign algorithm. While this dynamical interpretation of the matrix sign algorithm is of theoretical interest, it is not clear how to implement the algorithm numerically in an effective manner.  

We have not tested the performance of the matrix sign algorithm with random input in full generality. Instead, we have tested the deflation behavior of this algorithm in a more restricted setting by first precomputing $\sign(L_0)$ and then using Theorem~\ref{mark:thm-dlnt}. These results are not presented in this paper: the interested reader is referred to~\cite{ChristianDiss}.

\subsection{Deflation criterion}
\label{sec:deflation}
Consider a symmetric matrix $A$ with eigenvalues $\lambda_1 \geq \ldots \geq \lambda_n$, a symmetric matrix $B$, $\eps>0$ and the perturbed matrix $A+\eps B$ with eigenvalues $\lambda_1(\eps) \geq \ldots \lambda_n(\eps)$. Standard perturbation theory~\cite[Thm 5.1]{AppliedNumLA} implies
\begin{align}
|\lambda_i - \lambda_i(\epsilon)| \leq \epsilon \|B\|_2.
\end{align}
When deflating Jacobi matrices the perturbation 
matrix is of the form
\begin{align}
B = \left(
\begin{array}{cc}
0 & E_{1k}^T \\ 
E_{1k} & 0 
\end{array}
\right),
\end{align}
where the only non-zero entry in $E_{1k} \in \mathbb{R}^{(n-k)\times k}$ is a one in the upper right corner. Clearly, $\|B\|_2=1$
in this case. For the deflation of full symmetric matrices, the
perturbation matrix has the structure
\begin{align}
B = \left(
\begin{array}{cc}
0 & B_{21}^T \\ 
B_{21} & 0 
\end{array}
\right), 
\end{align}
where again $B_{21}\in\mathbb{R}^{(n-k)\times k}$, but now
all entries of $B$ satisfy $|b_{ij}|\leq 1$. In this case, 
one may show that $\|B\|_{2}\leq \sqrt{k(n-k)}$.

We now define the deflation criterion. If $L$ is a Jacobi matrix define
\begin{align}
\hat{\epsilon}_k = b_k.
\end{align}
If $L=(l_{ij})$ is a full symmetric matrix set
\begin{align}
\hat{\epsilon}_k = \sqrt{k(n-k)} \max_{\substack{k<i\leq n \\ 1\leq j\leq k }} |l_{ij}|
\end{align}
Assume $L_m$ is a sequence of iterates (Jacobi or full symmetric) obtained through an iterative eigenvalue algorithm.  For a given tolerance $\eps>0$ and initial matrix $L_0$ we define the {\em deflation time\/}
 \begin{equation}
	\label{eq:runtime}
\deftime_{n,\eps}(L_0)  = \min \left\{ m \left| \,\hat{\epsilon}_k(L_m)<\eps \,\, \mathrm{for\,\, some\/}\, 1 \leq k \leq n-1  \right. \right\}.
\end{equation}
For calculations based on the Hamiltonian flow \qref{eq:dlnt-equation} it is more natural to consider the real valued deflation time 
\begin{equation}
	\label{eq:cont-runtime}
\contdeftime_{n,\eps}(L_0)  = \inf \left\{ t >0 \left| \,\hat{\eps}_k(L(t))<\eps \,\, \mathrm{for\,\, some\/}\, 1 \leq k \leq n-1  \right. \right\}.
\end{equation}
The location where the matrix deflates is called the {\em deflation index\/}
\begin{align}
\defindex_{n,\eps}(L_0) = \arg\min_{1\leq k \leq n-1} \hat{\epsilon}_k( L_{\deftime_{n,\eps}(L_0)}).
\end{align}

There is an important difference between deflation and the asymptotic convergence guaranteed by Theorem~\ref{mark:thm-dlnt-big}. While Theorem~\ref{mark:thm-dlnt-big} may be used to compute asymptotic rates of convergence as $t \to \infty$~\cite[Thm 3]{ODESymmetricEig}, in practice the rate of convergence is determined by deflation and transients play an important role. We illustrate this with a simple example. 

Fix $\lambda_1 > \lambda_2>0$, let $\Lambda=\diag(\lambda_1, \lambda_2)$ and consider the QR flow on $\Symm(2)$ with the initial matrix 
\be
\label{eq:deflation-example}
 L_0 = Q_0 \Lambda Q_0^T, \quad Q_0=\left( \begin{array}{rr} \cos \theta_0 & \sin \theta_0 \\ \sin \theta_0 & -\cos \theta_0 \end{array} \right). 
\ee
According to Theorem~\ref{mark:thm-dlnt}, $\lim_{t \to \infty} L(t) =\Lambda$ for every $\theta_0$. However, if $\theta_0 \approx \pi/2$, $L_0$ is a small perturbation of $\diag(\lambda_2, \lambda_1)$,  and in practice, the algorithm would immediately deflate and return $L_0$. But according to Theorem~\ref{mark:thm-dlnt}, $L(t)$ must evolve so that the inital diagonal terms ``turn  around'' and are presented in the correct order $\diag(\lambda_1,\lambda_2)$ as $t \to \infty$ (see \qref{eq:permut-increase}). More generally, consider $\Lambda =(\lambda_1, \ldots, \lambda_n)$ with $\lambda_1 > \lambda_2 > \ldots > \lambda_n>0$. Each permutation $\sigma \in S_n$ yields a distinct fixed point $\Lambda_\sigma= (\lambda_{\sigma_1}, \ldots, \lambda_{\sigma_n})$ for the QR and Toda algorithms. In a numerical calculation, an initial condition close to $\Lambda_\sigma$ is immediately deflated.
Alternatively, iterates may pass close to one of the permutations $\Lambda_\sigma$ and again deflation occurs at finite times.  However, only the equilibrium $(\lambda_1, \ldots, \lambda_n)$ attracts  generic initial conditions~\cite{ODESymmetricEig}. Thus the notion of convergence as $t \to \infty$ and deflation are completely distinct.

\subsection{Ensembles}
\label{subsec:RMT}
We now introduce the six ensembles of random matrices that we will analyze. For general introductions on random matrices see~\cite{OrthPolRH,Edelman-Rao,MehtaRMT}.
The simplest way to construct an ensemble of random matrices is to choose entries independently subject only to the constraint of symmetry. Such ensembles are called {\em Wigner ensembles\/.} We also say that an ensemble lies in the {\em Wigner class\/} if the limiting spectral distribution for this ensemble is the Wigner semicircle law (described below).  We consider four Wigner ensembles in the Wigner class:
\begin{enumerate}
	\item the Gaussian Orthogonal Ensemble (GOE) (independent entries with $M_{ii} \sim \sqrt{2}\mathcal{N}(0,1)$, $M_{ij} \sim \mathcal{N}(0,1)$, $i>j$);
	\item the Gaussian Wigner ensemble (iid $M_{ij} \sim \mathcal{N}(0,1)$, $i \geq j$);
	\item the Bernoulli ensemble (iid $M_{ij} \sim \mathcal{B}$, $i \geq j$ );
	\item the Hermite--1 ensemble on Jacobi matrices (iid $a_k \sim \mathcal{N}(0,1)$, $k=1, \ldots, n$ and independent $b_k \sim \chi_k$, $k=0, \ldots, n-1$. ).
\end{enumerate}
(1)--(3) are ensembles of full symmetric matrices. The distinction between (1) and (2) is that the variance of the diagonal and off-diagonal entries of matrices in GOE is different to ensure orthogonal invariance (see~\cite{MehtaRMT}).   Hermite--1 is an ensemble of Jacobi matrices obtained by applying the Householder tridiagonalization procedure to the GOE ensemble. It is a remarkable fact that the entries remain independent under tridiagonalization (this is not true when matrices from ensembles (2) and (3) are tridiagonalized).

A choice of an ensemble of random, symmetric matrices is a choice of a probability measure on the space of symmetric matrices. When the matrix entries are independent this measure is a product measure. For example, the measure corresponding to GOE has density
\be
\label{eq:goe-hermite}
P_{\mathrm{GOE}}(M) = 2^{2n/2} (2\pi)^{-n(n+1)/4}e^{-\frac{1}{4}\tr(M^2)}.
\ee
For all these ensembles, while the matrix entries are independent, the eigenvalues are not. The joint density of eigenvalues for GOE and Hermite--1 may be computed explicitly and is given by the determinantal formula~\cite[Ch.3]{MehtaRMT}
\be
\label{eq:goe-spec}
f_1(\Lambda)= \frac{1}{Z_n}  |\triangle_n(\lambda)| e^{-\frac{|\lambda|^2}{2}}, \quad \triangle_n(\lambda) = \prod_{i<j} (\lambda_i-\lambda_j).
\ee
The normalization constant $Z_n$ may be computed explicitly. 
By contrast, while the analogues of \qref{eq:goe-hermite} for ensembles (2) and (3) are clear, there is no explicit analogue for \qref{eq:goe-spec}. 

The ensembles (1)--(4) are in the Wigner class, i.e. for each of these ensembles
\be
\label{eq:wigner1}
\lim_{n \to \infty} \frac{1}{n}\#\{\lambda_i \in \sqrt{n}(a,b)\} = \int_{a}^b \nu(x) dx, 
\ee
where $\nu(x)$ denotes the density of the {\em Wigner semicircle law\/}
\be
\label{eq:wigner2}
\nu(x) = \frac{1}{2\pi}\sqrt{4-x^2} \,\,\one_{|x|\leq 2}.
\ee

We will contrast our results on these ensembles  with two  ensembles of Jacobi matrices that are not in the Wigner class. These are:
\begin{enumerate}
	\item[5.] The uniform doubly stochastic Jacobi ensemble (UDSJ).
	\item[6.] The Jacobi uniform ensemble (JUE).
\end{enumerate}

Doubly stochastic Jacobi matrices of dimension $n\times n$ form
a compact polytope in $\mathbb{R}^{n-1}$ which can be equipped with its uniform
measure~\cite{RandomDoublyStochasticJacobi}. This is the UDSJ ensemble. We can approximately sample from this ensemble using a Gibbs sampler. 

JUE is defined using the spectral map $\specmap$ for $\Jac(n)$. Since we may describe Jacobi matrices by their spectral data $(\Lambda,u)$, a probability measure on the spectral data pulls back under $\specmap^{-1}$ to a probability measure on $\Jac(n)$. For JUE, we replace \qref{eq:goe-spec} with eigenvalues chosen independently and uniformly on an interval and $u$ distributed uniformly on the orthant $S_+^{n-1}$.  In our numerical simulations we assume the eigenvalues are uniformly distributed on $[-2\sqrt{n},2\sqrt{n}]$ because this interval corresponds to the support of the semicircle law and allows a comparison between JUE and ensembles in the Wigner class.
A particularly important aspect of JUE is that the eigenvalues do not repel one another. This strongly affects the statistics of $\deftime_{n,\eps}$ as  shown below (for unshifted QR and Toda, but not for shifted QR!).

\subsection{The normalized deflation time}
We have now defined the algorithms, ensembles and deflation criterion. For a given algorithm and ensemble, $\deftime_{n,\eps}(L)$ and $\defindex_{n,\eps}(L)$ are random variables that depends on the random initial matrix $L$ and $\eps>0$. We explore the empirical distributions of $\deftime_{n,\eps}$ and $\defindex_{n,\eps}$ in simulations.  Our main empirical finding is that for each algorithm these empirical distributions collapse into a universal distribution for the Wigner ensembles (1)--(4). Let $\mu_{n,\eps}$ and $\sigma^2_{n,\eps}$ denote the empirically determined mean and variance of $\deftime_{n,\eps}(L)$ for a particular algorithm and ensemble. 

Our simulations suggest that the normalized deflation time
\be
\label{eq:uni-scaling}
\resctime_{n,\eps}= \frac{\deftime_{n,\eps} - \mu_{n,\eps}}{\sigma_{n,\eps}} 
\ee
converges in distribution as $n \to \infty$ and $\eps \to 0$ and that the limit is the same for ensembles in the Wigner class (see Figure~\ref{fig:QRnormalized-runtimes-wigner} and Figure~\ref{fig:Todanormalized-runtimes-wigner})). Both $\mu_{n,\eps}$ and $\sigma_{n,\eps}$ are computed empirically.  Our numerical calculations also suggest that $\mu_{n,\eps} \sim C |\log\eps|$ for all ensembles in the  Wigner class~(see Figure~\ref{fig:QR-means} and Figure~\ref{fig:Toda-means}). As already noted above, a suprising outcome of our simulations is that universality for shifted QR is more encompassing, and actually holds for all six ensembles (1)--(6).

In order to prove convergence in distribution of $\resctime_{n,\eps}$ it is first necessary to estimate the mean and variance of $\deftime$. We present below a calculation of $\mu_{2,\eps}$ that illustrates the subtle role of eigenvalue repulsion.

\subsection{The scaling of the expected deflation time}
\label{subsec:example}
In this section we  estimate the expected deflation time of the Toda flow on $\Symm(2)$. We show that 
\be
\mu_{2,\eps,\mathrm{GOE}} \sim C |\log \eps|, \quad\mathrm{but}\quad  \mu_{2,\eps,\mathrm{JUE}} \sim C |\log \eps|^2, \quad \eps \to 0.
\ee
The interval of support for the JUE density is chosen here to be $[-1,1]$. This choice only affects the prefactor $C$, not the term $|\log \eps|^2$.

In order to establish these asymptotics, we first determine the deflation time $\contdeftime_\eps$ as a function of the initial condition (for brevity we write $\contdeftime_\eps$ for $\contdeftime_{2,\eps}$ since $n=2$ is fixed).  Since $M(t) \in \Symm(2)$ we may write $M=U(t)\Lambda U(t)^T$, where $\Lambda=\diag(\lambda_1,\lambda_2)$, $\lambda_1>\lambda_2$, and
\be
\label{eq:inspec1}
 U(t) = \left( \begin{array}{rr} \cos \theta(t) & \sin \theta(t) \\  \sin \theta(t)  & -\cos \theta(t) \end{array}\right).
\ee
Note that $m_{12}>0$ corresponds to $\theta \in (0,\pi/2)$. We use Theorem~\ref{mark:thm-dlnt} to obtain
\be
m_{12}(t)= (\lambda_1-\lambda_2) \cos \theta(t) \sin \theta(t)
= (\lambda_1-\lambda_2)\cdot\frac{e^{t(\lambda_2-\lambda_1)}\cdot \tan \theta_0}{1 + e^{2t(\lambda_2-\lambda_1)}\tan^2\theta_0}.
\ee
Here $\theta_0=\theta(0)$. Now we set $m_{12}(\contdeftime_\eps)=\eps$ and solve to find 
\be
\label{eq:deftimeformula}
(\lambda_1-\lambda_2)\contdeftime_\eps = \left\{\begin{array}{cc}
0 & m_{12}(0)\leq\epsilon \\
\log\tan\theta_0 - \log\left[
\frac{\lambda_1-\lambda_2}{2\epsilon} - \sqrt{\frac{(\lambda_1-\lambda_2)^2}{4\epsilon^2}-1}
\right] & m_{12}(0) > \epsilon.
\end{array}
\right.
\ee
The asymptotics of $\contdeftime_\eps$ are easily determined. We have
\be
\label{eq:deftimasympt}
(\lambda_1 -\lambda_2)\contdeftime_\eps \sim -\log \eps + \log \tan \theta_0 + \log (\lambda_1-\lambda_2), \quad \eps \to 0.
\ee

To compute the mean deflation time for GOE and JUE we first change to spectral variables. As noted above,  the spectral map $\mathcal{S}$ is a diffeomorphism between the set of $2\times 2$ symmetric matrices with $m_{12}>0$ and the set $\{\lambda_1> \lambda_2\} \times (0, \pi/2)$. The Jacobian of this transformation is $\lambda_1 -\lambda_2$ so that 
\be
\label{eq:inspec2}
dm_{11}dm_{22} dm_{12} = (\lambda_1 -\lambda_2) d\lambda_1 d\lambda_2 d\theta.
\ee
The mean deflation time for GOE is then given by 
\be
\label{eq:meangoe}
\mu_{2,\eps,\mathrm{GOE}} = \frac{1}{Z_1}
\int_{-\infty}^\infty \int_{-\infty}^{\lambda_1} \int_0^{\pi/2} \contdeftime_\eps(\lambda_1, \lambda_2, \theta) e^{-(\lambda_1^2+\lambda_2^2)/4} (\lambda_1-\lambda_2) d\lambda_1 d\lambda_2 d\theta.
\ee
For JUE, the eigenvalues are chosen uniformly from $[-1,1]$ and we find
\be
\label{eq:meanjue}
\mu_{2,\eps,\mathrm{JUE}} = \frac{1}{Z_2}
\int_{-1}^1 \int_{-1}^{\lambda_1} \int_0^{\pi/2} \contdeftime_\eps(\lambda_1, \lambda_2, \theta)   d\lambda_2 d\lambda_2 d\theta.
\ee
Here $Z_1$ and $Z_2$ are normalizing constants for these probability densities.

The asymptotic behavior of \qref{eq:deftimasympt}, combined  with \qref{eq:meangoe} and \qref{eq:meanjue}, suggests the following leading order behavior as $\eps \to 0$:
\ba
\nonumber
\mu_{2,\eps,\mathrm{GOE}} && \sim   \frac{|\log \eps|}{Z_1} \int_{-\infty}^\infty \int_{-\infty}^{\lambda_1} \int_0^{\pi/2} e^{-(\lambda_1^2+\lambda_2^2)/4} d\lambda_1 d\lambda_2 d\theta \sim C_1 |\log \eps|\,,\\
\nonumber
\mu_{2,\eps,\mathrm{JUE}} && \sim   \frac{|\log \eps|}{Z_2} \int_{-1}^1 \int_{-1}^{\lambda_1} \int_0^{\pi/2} \frac{1}{\lambda_1 -\lambda_2}\one_{m_{12}>\eps} d\lambda_1 d\lambda_2 d\theta \sim C_2 |\log \eps|^2.
\ea
Here $C_i$ denote constants that may be computed explicitly.  The second integral is divergent without the cut-off $\one_{m_{12}>\eps}$: the cut-off gives rise to an additional factor of $|\log \eps|$. With more effort, these formal estimates may be made rigorous.

The analogous calculations for $M(t) \in \Jac(n)$, $n>2$ are quite subtle. For Jacobi matrices deflation occurs when $M(t)$ approaches the boundary  $\partial \Jac(n)$ of $\Jac(n)$ (see for example~\cite[Figs. 6 and 7]{ODESymmetricEig}).  A theoretical analysis of such deflations, which we have not carried out yet,  is a significant challenge as it requires a detailed understanding of the geometry of both the flow and the initial probability distribution in the vicinity of $\partial \Jac(n)$ in high dimensions. For this reason, we are reduced to using the empirical mean $\mu_{n,\eps}$ and variance $\sigma^2_{n,\eps}$ to define the normalized deflation time in~\qref{eq:uni-scaling}.

\section{Results}
\label{sec:numerical-results}

We generated a large number (typically 5000--10,000) of samples of the deflation time and the deflation index for each choice of the following parameters:
\begin{enumerate}
\item an eigenvalue algorithm  (QR without shift, QR with shift, Toda).
\item a random matrix ensemble.
\item matrix size $n$ (typically ranging from $10,30,\ldots 190$.)
\item tolerance $\eps$ (typically $10^{-k}$, $k=2,4,6,8$).
\end{enumerate}
We present a representative sample of our main results. Further statistical tests, figures and tables that amplify our conclusions may be found in~\cite{ChristianDiss}.

\subsection{Unscaled deflation time statistics for GOE}
\begin{figure}[h]
\subfloat[][]{
\includegraphics[width=7cm]{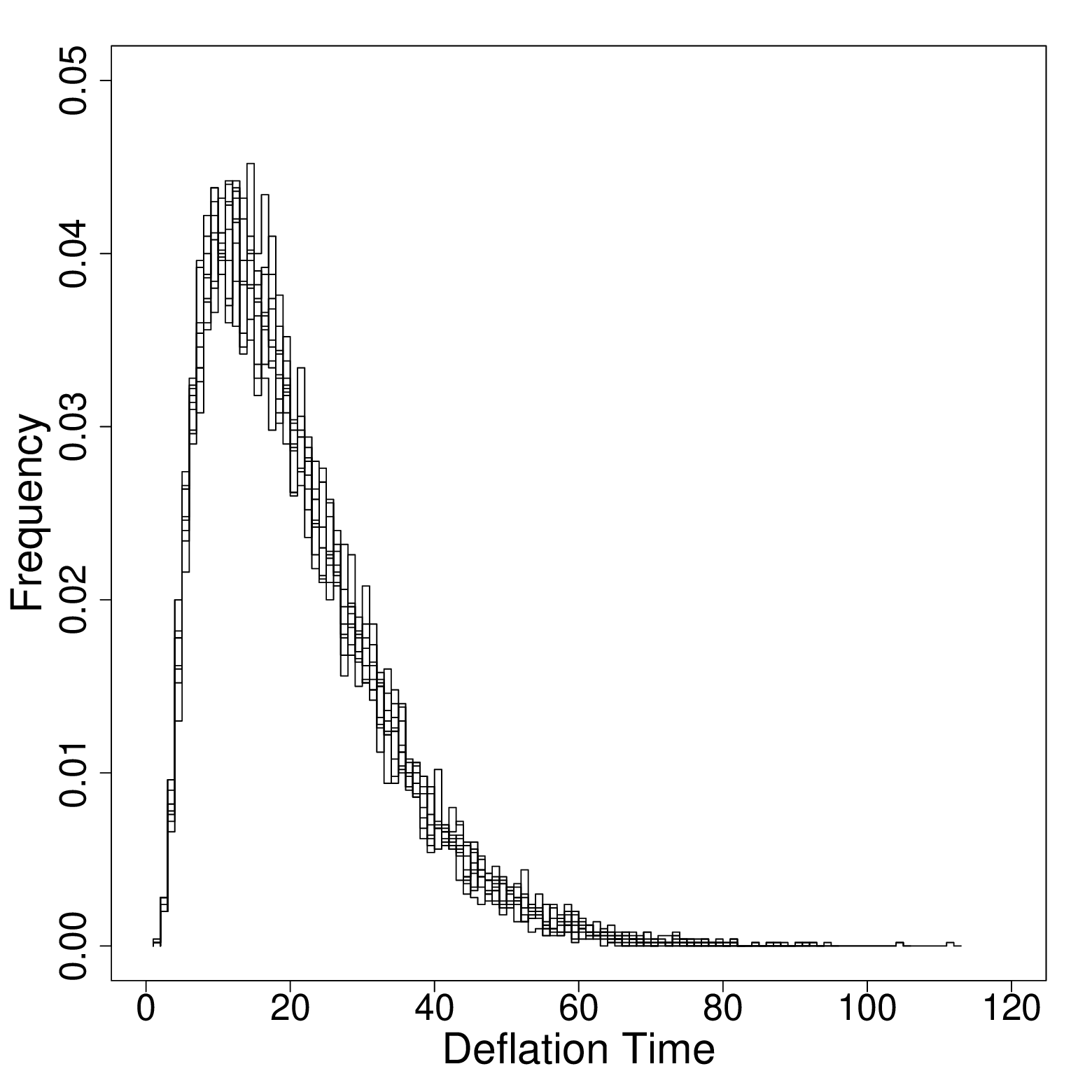}}
\subfloat[][]{
\includegraphics[width=7cm]{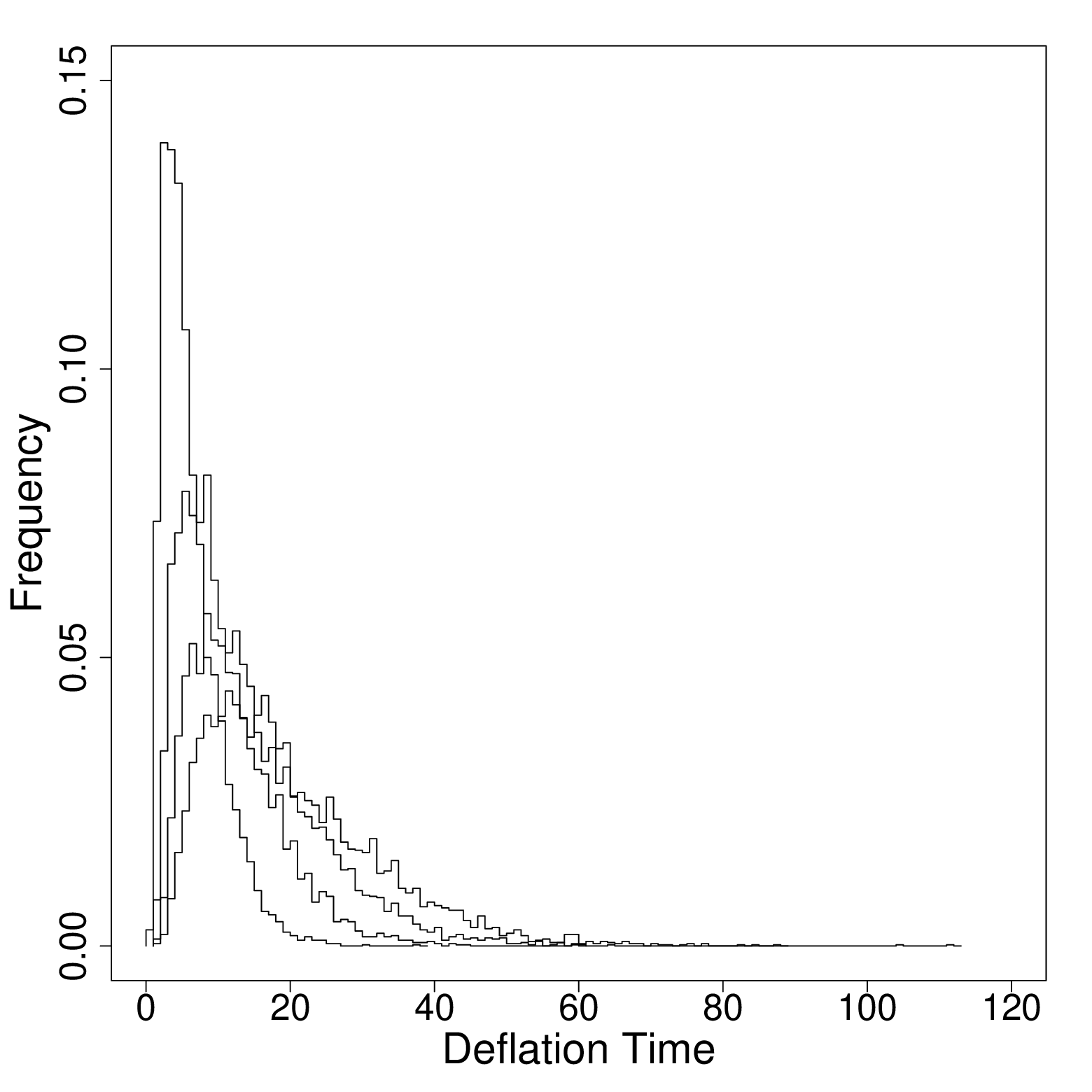}}
\caption{ {\bf The QR algorithm applied to GOE.\/}  
(a) Histogram for the empirical frequency of $\contdeftime_{n,\eps}$ as $n$ ranges from $10,30,\dots,190$ for a fixed deflation tolerance $\eps=10^{-8}$. The curves (there are $10$ of them plotted one on top of another) do not depend significantly on $n$. (b) Histogram for empirical frequency of $\contdeftime_{n,\eps}$ when $\epsilon=10^{-k}$, $k=2,4,6,8$ for fixed matrix size $n=190$. The curves move to the right as $\eps$ decreases.}
\label{fig:goe-qr}
\end{figure}

\begin{figure}[h]
\subfloat[][]{
\includegraphics[width=7cm]{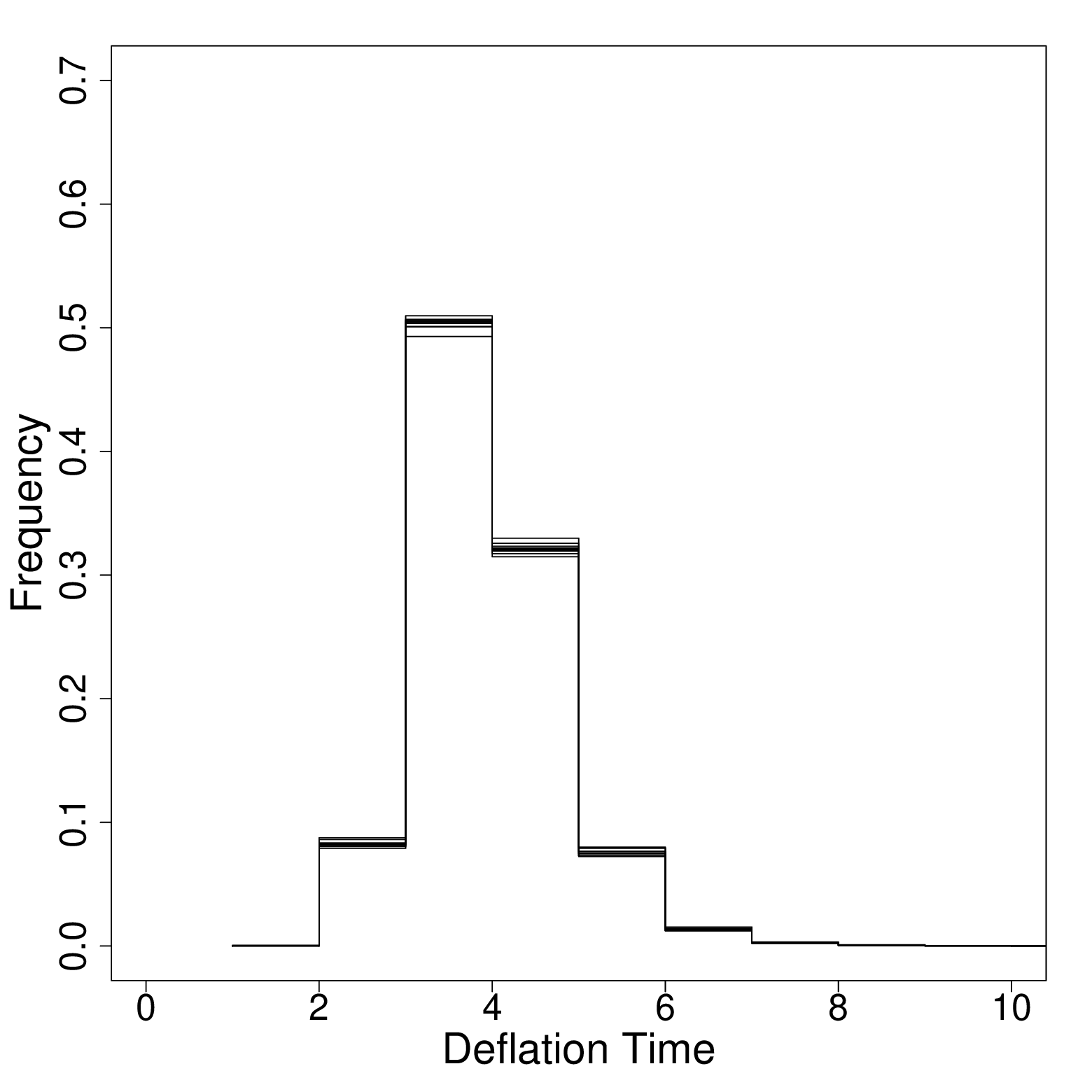}}
\subfloat[][]{
\includegraphics[width=7cm]{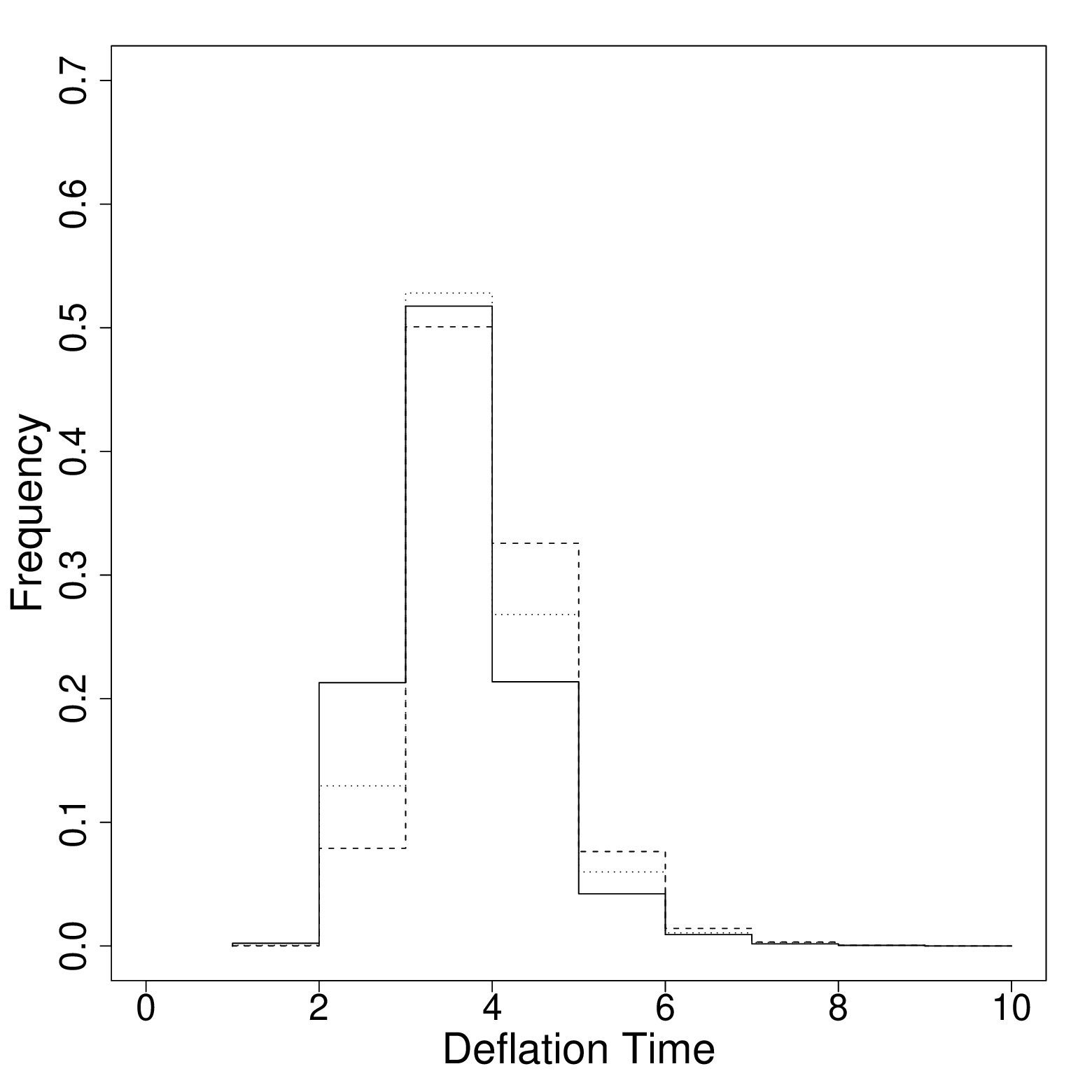}}
\caption{ {\bf The shifted QR algorithm applied to GOE.\/}  
(a) Histogram for the empirical frequency of $\contdeftime_{n,\eps}$ as $n$ ranges from $10,30,\dots,190$ for a fixed deflation tolerance $\eps=10^{-12}$. As for the unshifted QR algorithm, the curves are insensitive to $n$, though the tail becomes more pronounced for larger $n$.  (b) Histogram for empirical frequency of $\contdeftime_{n,\eps}$ when $\epsilon=10^{-k}$, $k=8, 10, 12$ for fixed matrix size $n=190$. The curves move to the right as $\eps$ decreases.}
\label{fig:goe-qrw}
\end{figure}

\begin{figure}[h]
\subfloat[][]{
\includegraphics[width=7cm]{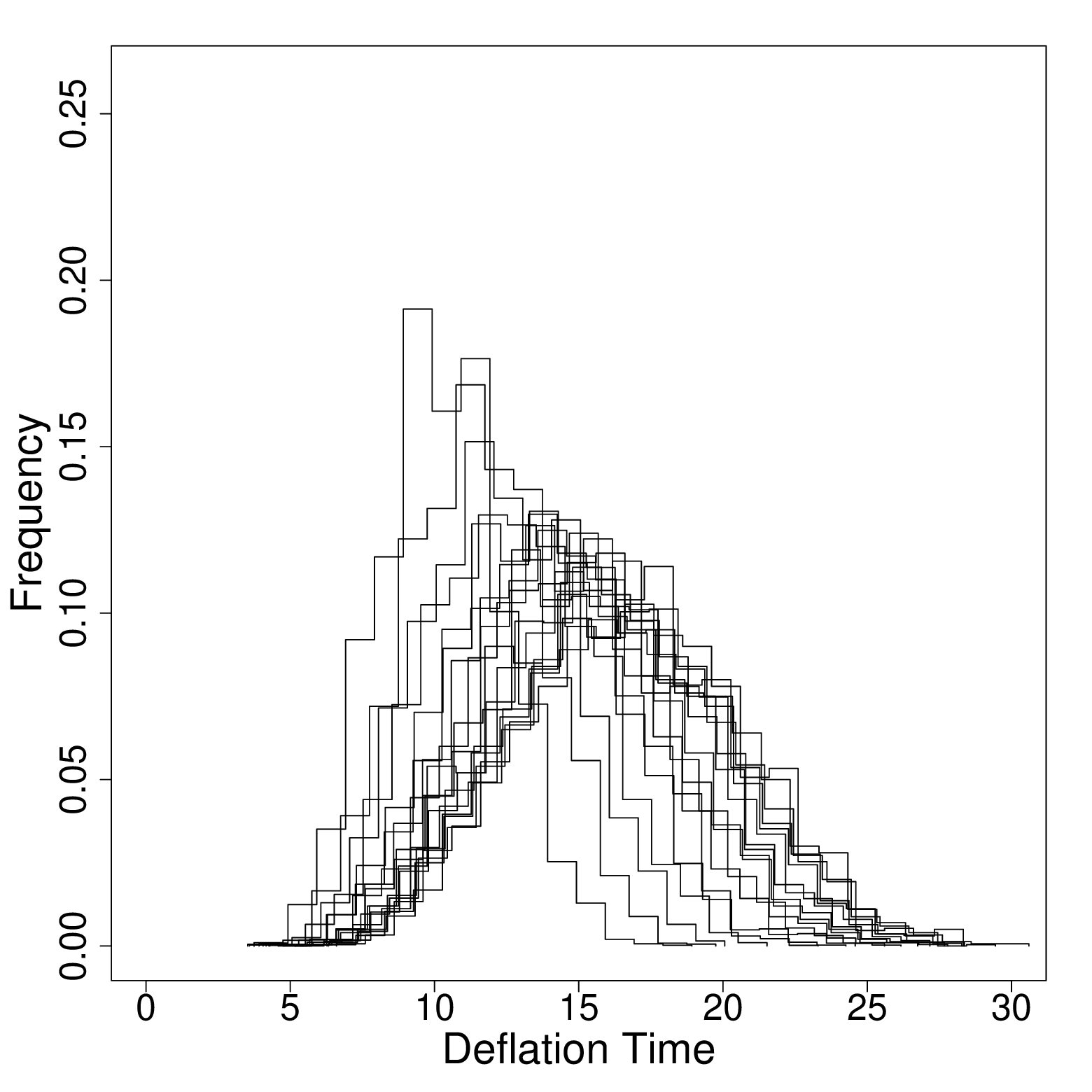}}
\subfloat[][]{
\includegraphics[width=7cm]{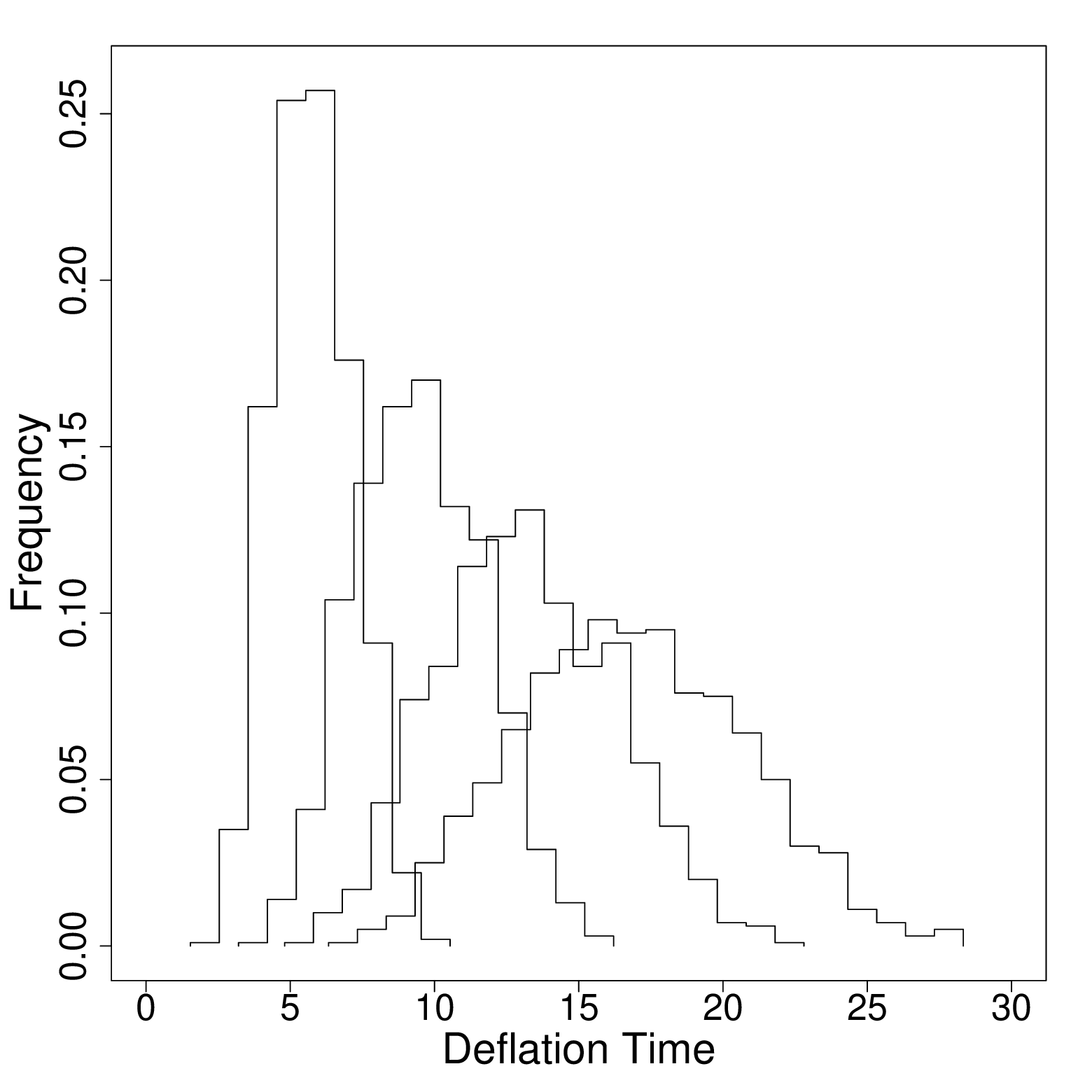}}
\caption{ {\bf The Toda algorithm applied to GOE.\/} 
(a) Histogram for empirical frequency of $\contdeftime_{n,\eps}$ as $n$ ranges from $10,30,\dots,190$ for a fixed deflation tolerance $\eps=10^{-8}$. The curves drift to the right as $n$ increases. (b) Histogram for empirical frequency of $\contdeftime_{n,\eps}$ when $\epsilon=10^{-k}$, $k=2,4,6,8$ for fixed matrix size $n=190$. The curves move to the right as $\eps$ decreases.}
\label{fig:goe-toda}
\end{figure}

We first present deflation time statistics for $\contdeftime_{n,\eps}$ for a fixed ensemble (GOE) for both the QR (shifted and unshifted) and Toda algorithms. The statistics of $\contdeftime_{n,\eps}$ for the unshifted QR algorithm are shown in Figure~\ref{fig:goe-qr}. 
Similar statistics for the QR algorithm with Wilkinson shift and the Toda algorithm are shown in Figure~\ref{fig:goe-qrw} and Figure~\ref{fig:goe-toda} respectively.  These figures reflect the typical dependence of these algorithms on $n$ and $\eps$ for ensembles (1)--(6). Similar statistics for other ensembles  may be found in~\cite[chapter 7]{ChristianDiss}. In all cases, we observe that the histograms for the QR algorithm are relatively insensitive to $n$ and shift to the right as $\eps$ decreases. The effect of the Wilkinson shift is to sharply reduce the number of iterations required (note the different scale of the abcissa in Figure~\ref{fig:goe-qr} and Figure~\ref{fig:goe-qrw}). The values of $\epsilon$ for shifted QR are much smaller than those chosen for QR without shifts. This choice is necessary to generate a viable data set for the shifted QR algorithm with sufficient variation in the deflation time. The histograms for the Toda algorithm shift to the right as $n$ increases and $\eps$ decreases, as discussed below.


\subsection{Normalized deflation time and universality for the Wigner class}
\begin{figure}[h]
\subfloat[][]{
\includegraphics[width=7cm]{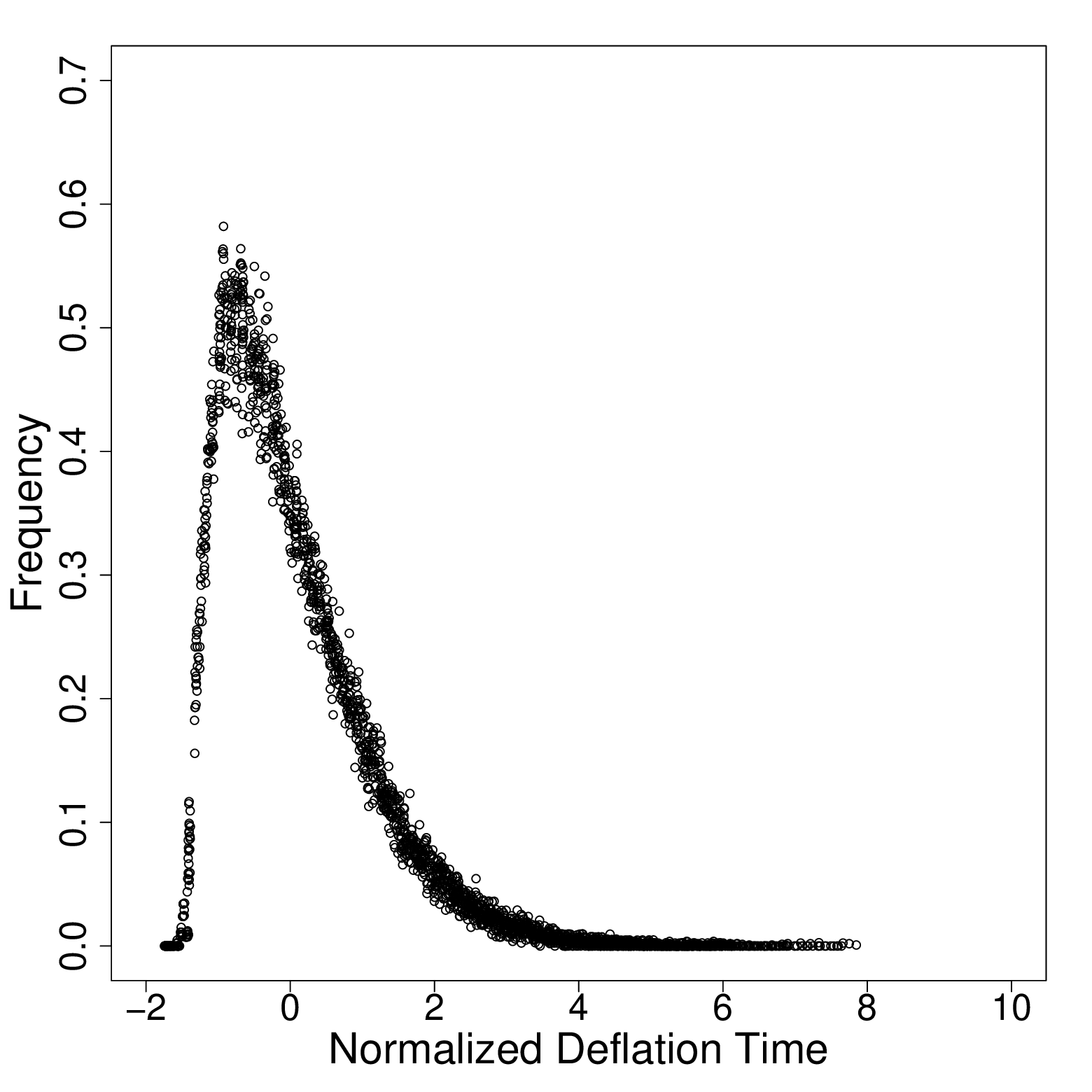}}
\subfloat[][]{
\includegraphics[width=7cm]{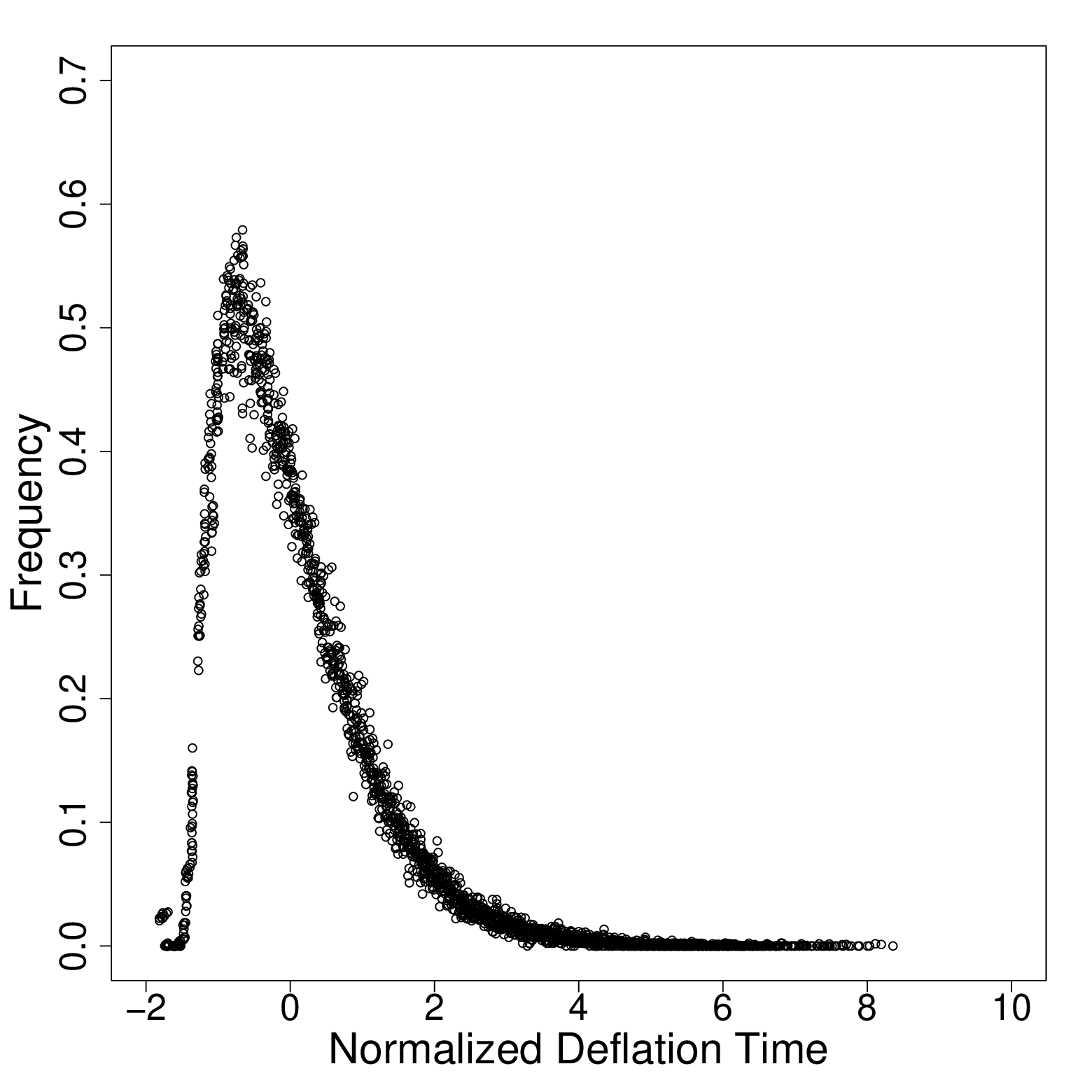}}\\
\subfloat[][]{
\includegraphics[width=7cm]{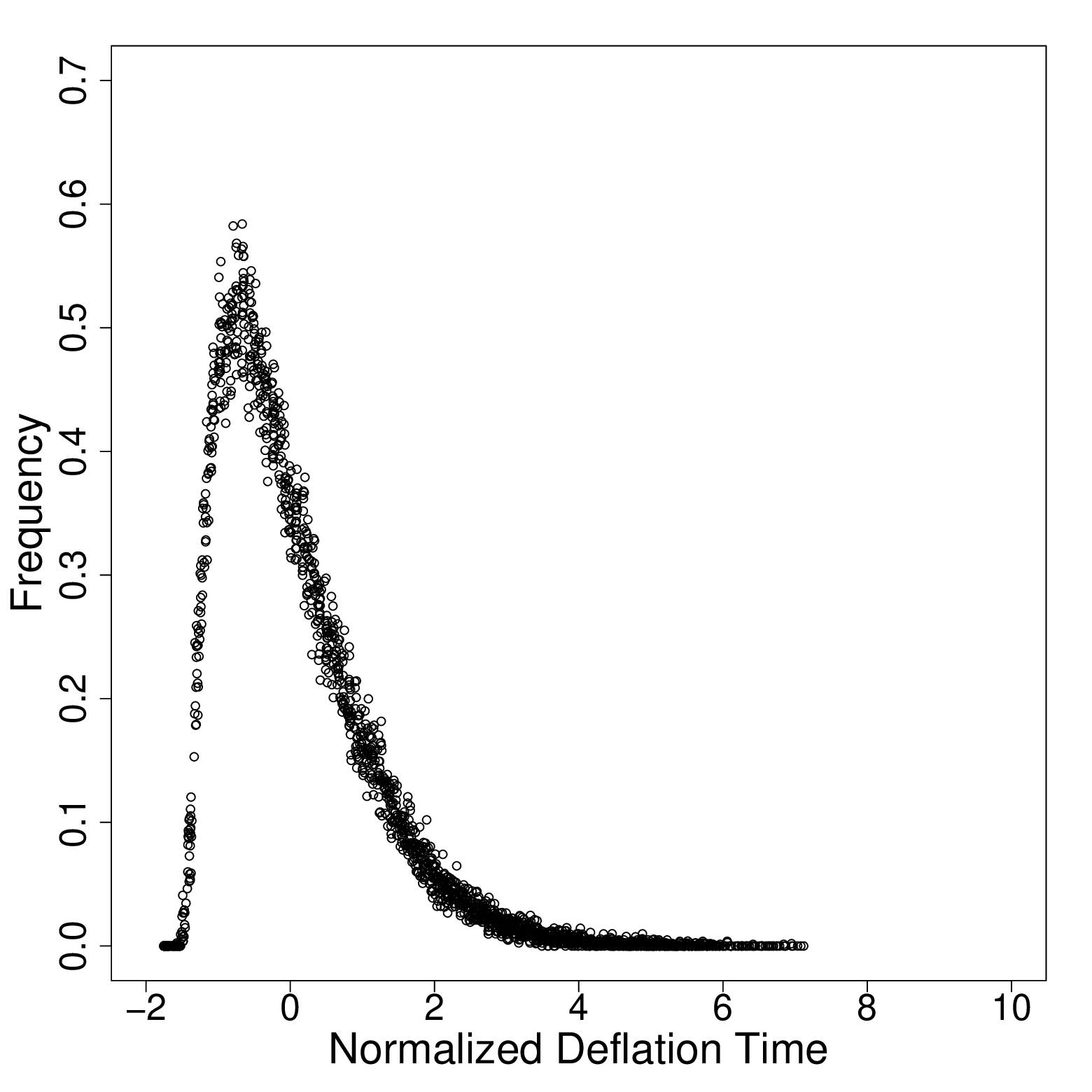}} 
\subfloat[][]{
\includegraphics[width=7cm]{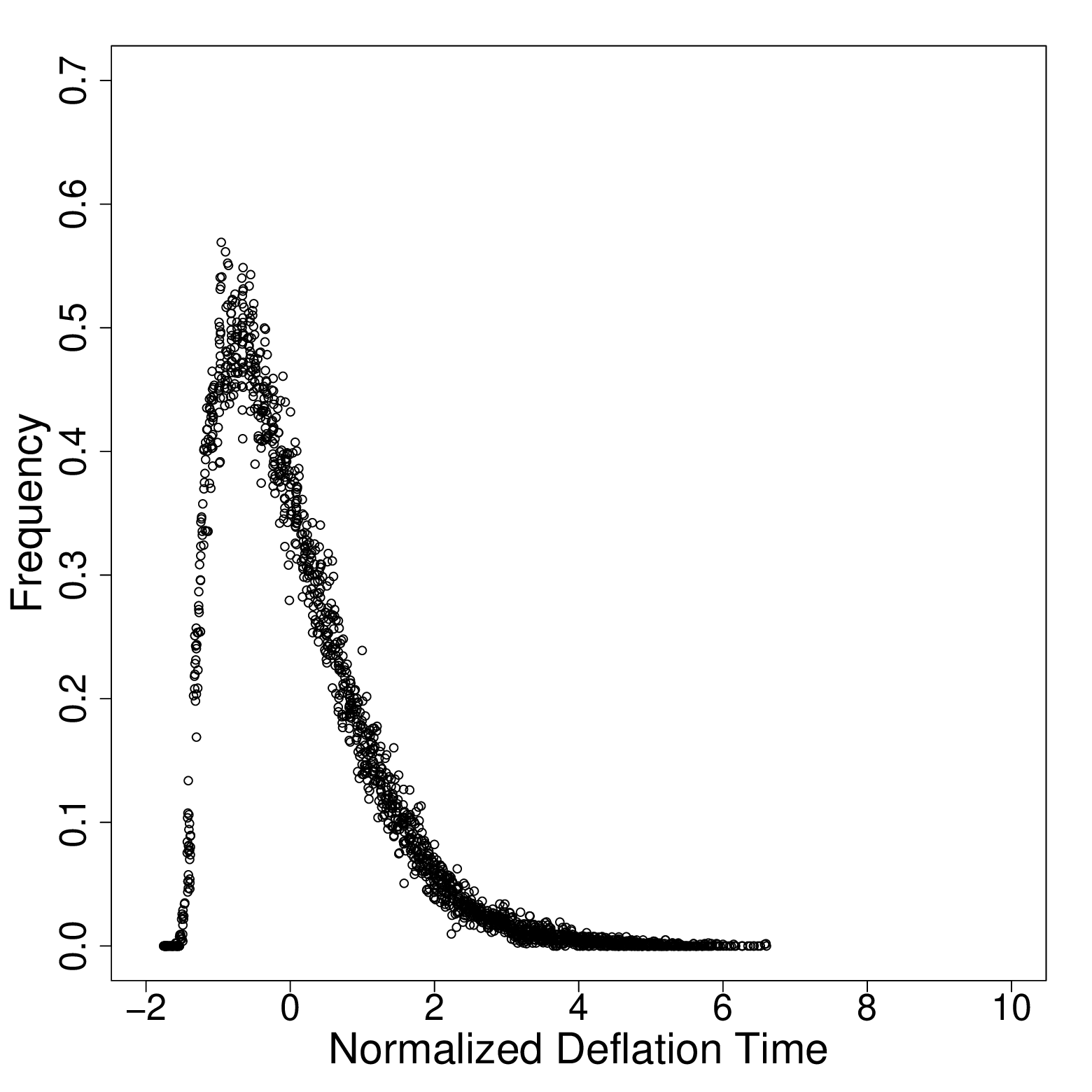}}
\caption{ {\bf Universal deflation time statistics for the QR algorithm applied to the Wigner  class.\/} Empirical deflation time normalized as in \qref{eq:uni-scaling} for $\epsilon=10^{-k}$, $k=2,4,6,8$ and $n$ ranging from $10,30,\dots,190$. The random matrix ensembles are (a) GOE; (b) Hermite-1; (c) Gaussian Wigner; and (d) Bernoulli. Each of the figures
(a), (b), (c), and (d) is obtained by rescaling the data of $10 \times 4$ fixed-$n$ and
fixed-$\epsilon$ histograms and plotting them together. All these data are observed to collapse onto one universal curve. Plotting all $160$
histograms in one figure (as in Figure~\ref{fig:combined-histograms} below) further demonstrates the universality of the deflation algorithm.}
\label{fig:QRnormalized-runtimes-wigner}
\end{figure}

\begin{figure}[h]
\subfloat[][]{
\includegraphics[width=7cm]{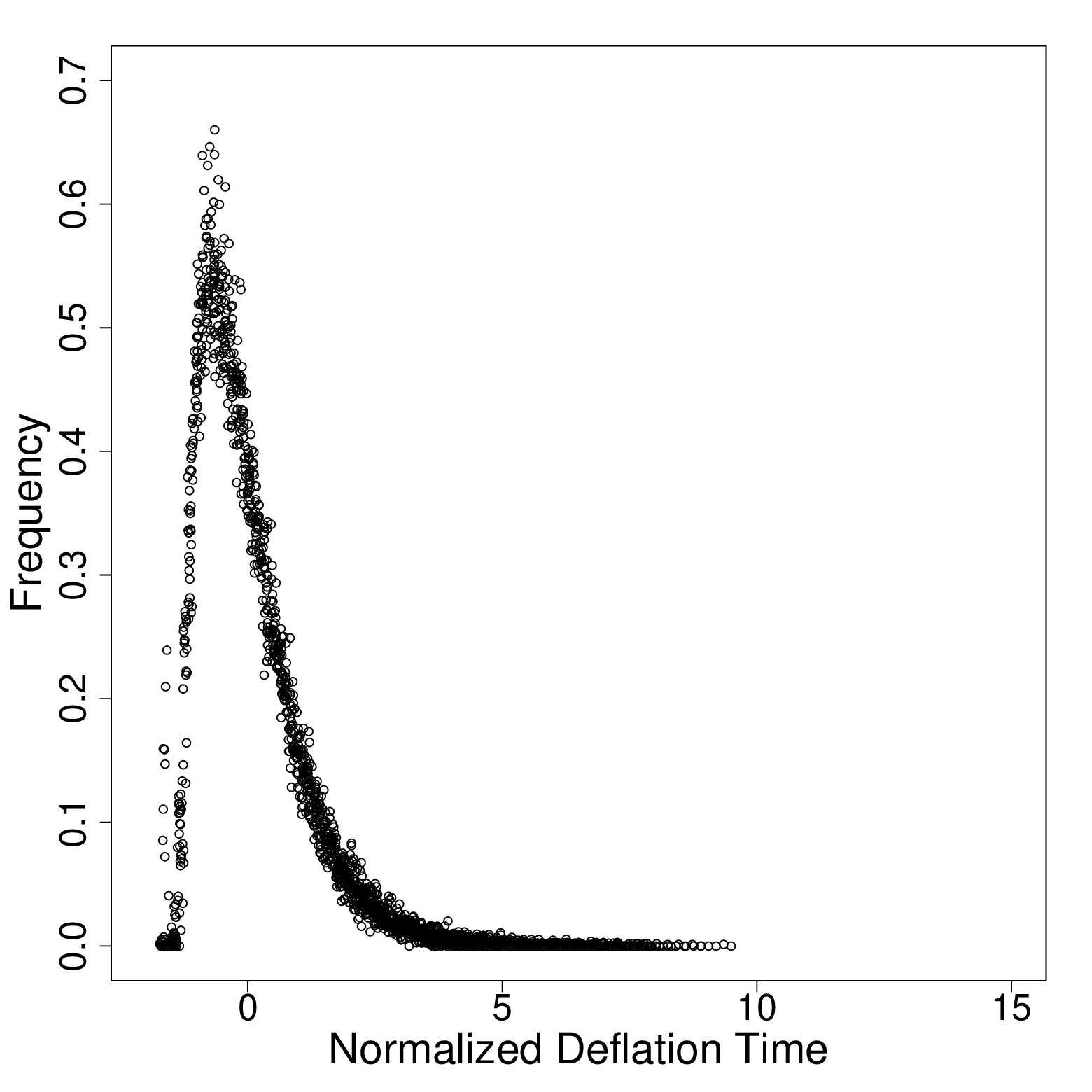}}
\subfloat[][]{
\includegraphics[width=7cm]{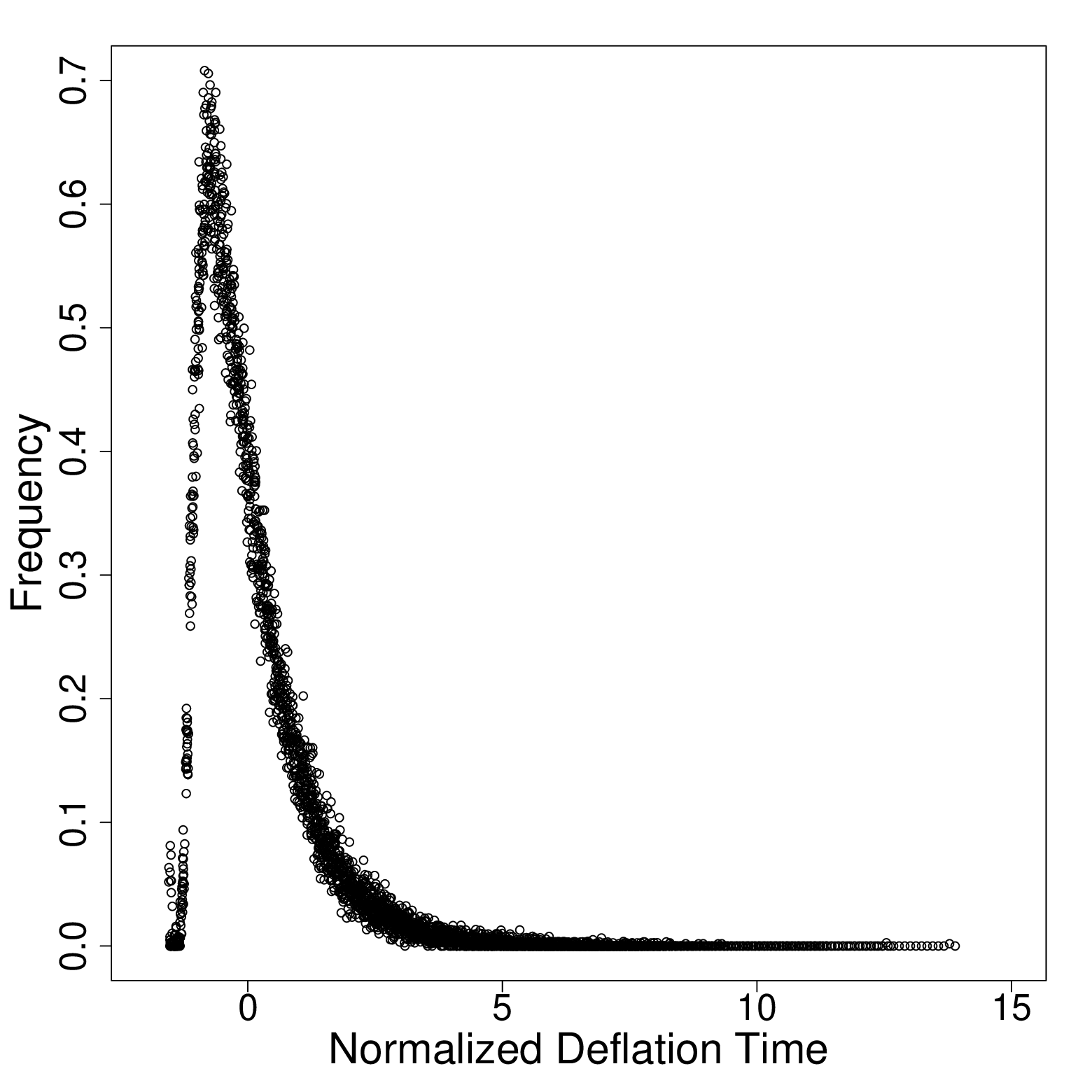}}
\caption{{\bf The QR algorithm applied to non-Wigner ensembles.\/} 
Normalized empirical deflation time distributions for the QR  algorithm with $\epsilon=10^{-k}$, $k=2,4,6,8$ and $n$ ranging from $10,30,\dots,190$.  The random matrix ensembles are (a) UDSJ and (b) JUE. Each figure contains the normalized empirical data of $40$ fixed-$n$ and fixed-$\eps$ histograms. All these data are observed to collapse onto a single curve. However, these curves are not the same for UDSJ and JUE and neither of these coincides with the curve for Wigner data shown in 
Figure~\ref{fig:QRnormalized-runtimes-wigner}.
\label{fig:QR-normalized-runtimes-nonwigner}
}
\end{figure}

We now present results that show the collapse of all data onto universal curves depending only on the algorithm under the rescaling~\qref{eq:uni-scaling}. The statistics of the empirical mean $\mu_{n,\eps}$ and standard deviation $\sigma_{n,\eps}$ are discussed a little later. The empirical distribution of the normalized deflation time $\resctime_{n,\eps}$ for the QR algorithm with intial data from the Wigner ensembles is shown in  Figure~\ref{fig:QRnormalized-runtimes-wigner}. All the data contained in Figure~\ref{fig:goe-qr} collapse onto the single curve seen in Figure~\ref{fig:QRnormalized-runtimes-wigner}(a). Analogous data for the other Wigner class ensembles (2)--(4) collapse onto the {\em same universal curve\/}. The normalized deflation time distributions for UDSJ and JUE are shown in  Figure~\ref{fig:QR-normalized-runtimes-nonwigner}. While we again observe a collapse of the data, it is not onto the curve of Figure~\ref{fig:QRnormalized-runtimes-wigner}(a).  This contrast is amplified in the comparison of the tails of the normalized deflation time (see Figure~\ref{fig:combined-histograms}). QQ plots that directly compare the histograms of these distributions may be found in~\cite{ChristianDiss}.

The most obvious difference between the behavior of the unshifted and shifted QR algorithm is that the spread in the deflation time for the shifted QR algorithm is much narrower. However, this does not seem to affect our general conclusion that there is universality for each Hamiltonian eigenvalue algorithm. The normalized deflation time distribution for shifted QR is shown in  Figure~\ref{fig:S-QRnormalized-runtimes-wigner} and Figure~\ref{fig:S-QR-normalized-runtimes-nonwigner}.  Moreover,  for shifted QR, the deflation times vary far less with the choice of underlying ensemble than the unshifted QR algorithm. In particular, we see a strong similarity for all ensembles in  Figure~\ref{fig:S-combined-histograms}. This behavior is in contrast with that of unshifted QR, shown in Figure~\ref{fig:combined-histograms}.

Finally, we have also observed universality for the Toda algorithm. The  empirical distribution of the normalized deflation time for the Wigner ensembles is shown in Figure~\ref{fig:Todanormalized-runtimes-wigner}. Again, all the data contained in Figure~\ref{fig:goe-toda} collapse onto the single curve seen in Figure~\ref{fig:Todanormalized-runtimes-wigner}(a). Further, analogous data for the other Wigner ensembles (2)--(4) collapse onto the same curve. The data for UDSJ and JUE  collapse under normalization, but not onto the same distribution (see  Figure~\ref{fig:QR-normalized-runtimes-nonwigner} and Figure~\ref{fig:Toda-combined-histograms}).

\begin{remark} We note that for both the QR and Toda algorithms the limiting distribution of the normalized deflation time $\resctime_{n,\eps}$ for UDSJ and 
JUE is distinct from that of ensembles in the Wigner class. This raises the interesting issue in random matrix theory whether UDSJ and JUE are in the same universality class as Wigner ensembles and invariant ensembles. As JUE does not have eigenvalue repulsion built in, this is unlikely to be the case.
\end{remark}


\subsection{Universal tails for deflation times}
We used a hypothesis testing approach to quantify the statement that the rescaled deflation time has an exponential tail for QR 
and a Gaussian tail for Toda. Our approach is modeled on the methodology of \cite{powerlawEmpiricalData}. 
Given deflation time data $D$ we perform maximum likelihood estimation
of parameters for distribution families  conditioned on observing only values above a cutoff value $x_{\min}(D)$ and use
a semiparametric approach to compute p--values for these parameters.
Based on $D$ and our parameter estimate, we compute resampled data sets
and a modified Kolmogorov--Smirnov statistic measuring the distance between
the empirical distribution function and the ones resulting from our 
maximum likelihood estimates. The semiparametric p--value is given
as the proportion of instances that the resampled data sets yield larger
modified KS statistics than the original. If this p--value is large we accept the hypothesis that the original data set has in fact the proposed decay in the right tail.

\begin{figure}[h]
\subfloat[][]{
\includegraphics[width=7cm]{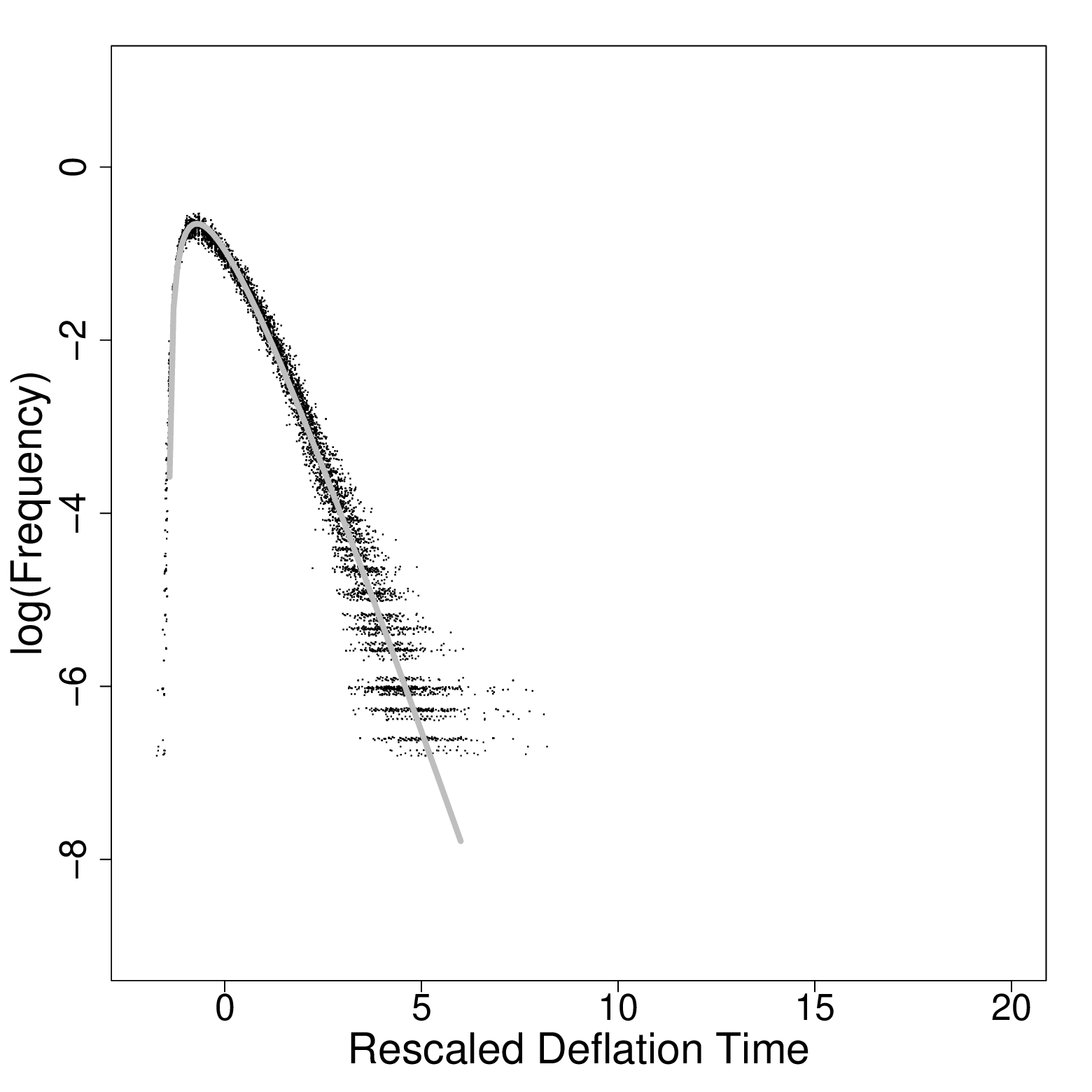}}
\subfloat[][]{
\includegraphics[width=7cm]{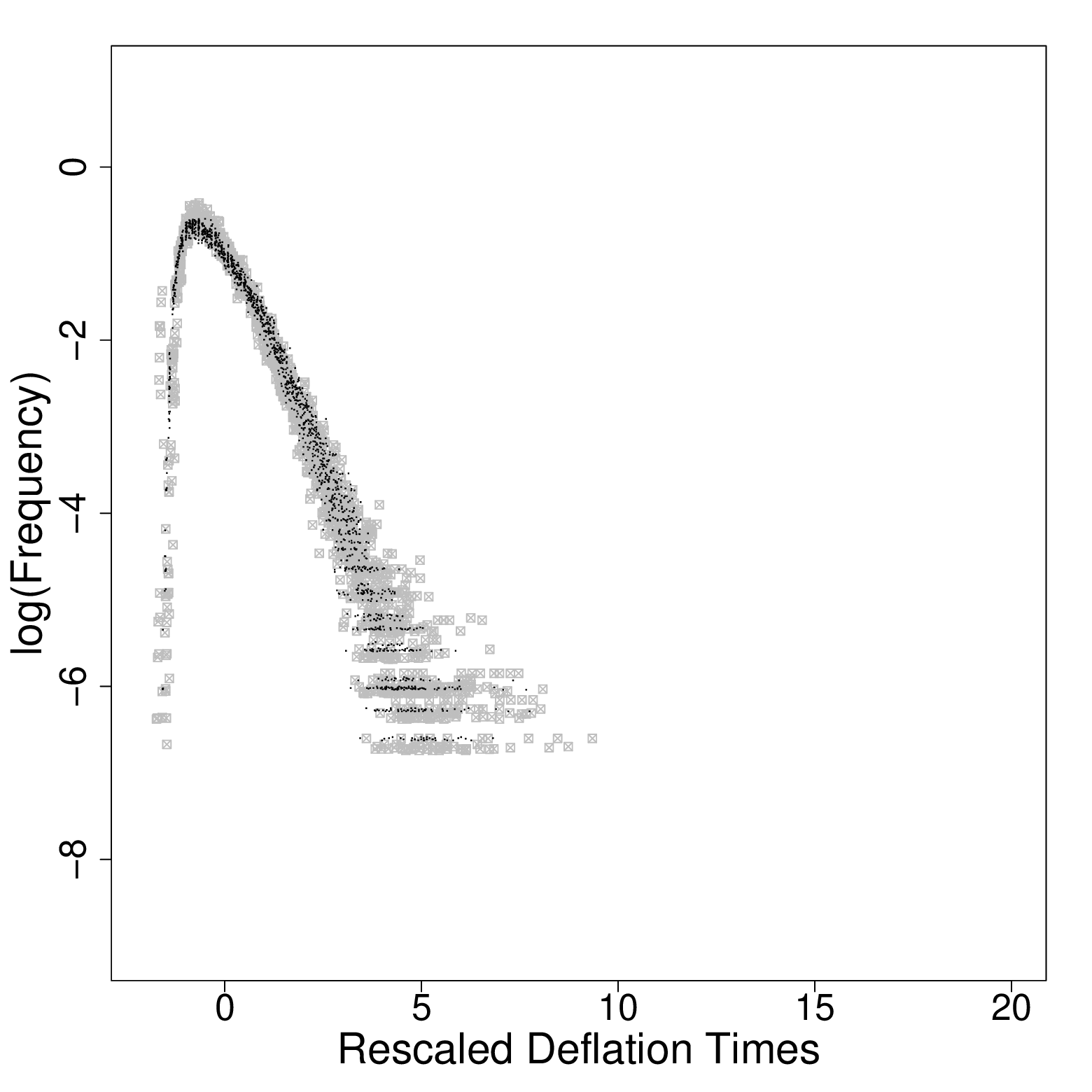}}
\caption{ {\bf Exponential tail for the QR algorithm\/}.
Histograms of the normalized deflation time for the  QR  algorithms on a logarithmic scale. (a) {\bf Wigner data. \/} Empirical normalized deflation time distributions 
from all $160$ histograms of Wigner class initial data (black dots) are compared with a gamma distribution with parameters 
$k=2$ and $\theta=1$ shifted to mean zero (gray line). 
 (b) {\bf non-Wigner data.\/}  Empirical normalized deflation time distributions from $40$ GOE histograms (black
dots) is contrasted with data from $40$ UDSJ histograms (gray squares). 
}
\label{fig:combined-histograms}
\end{figure}

\begin{figure}[h]
\subfloat[][]{
\includegraphics[width=7cm]{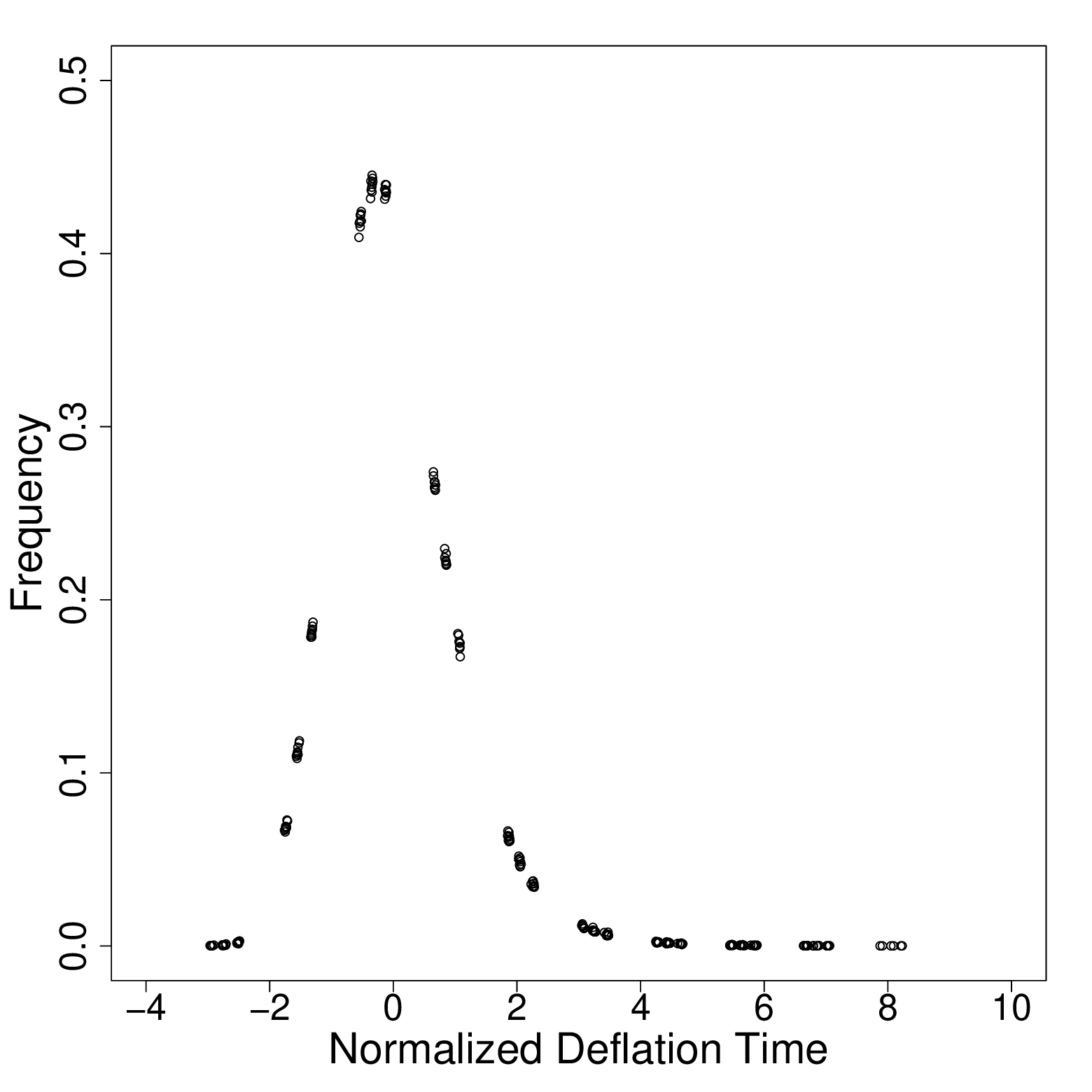}}
\subfloat[][]{
\includegraphics[width=7cm]{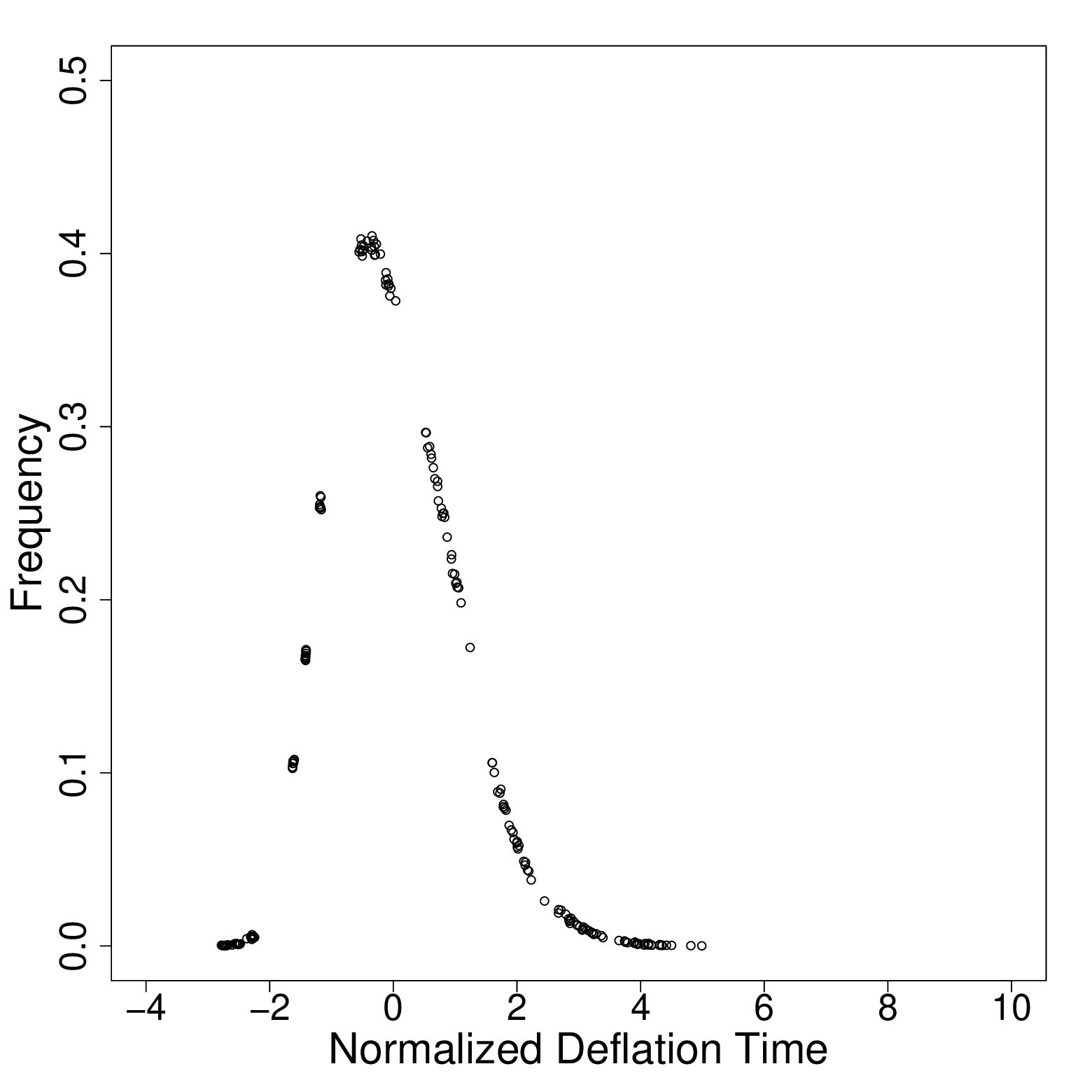}}\\
\subfloat[][]{
\includegraphics[width=7cm]{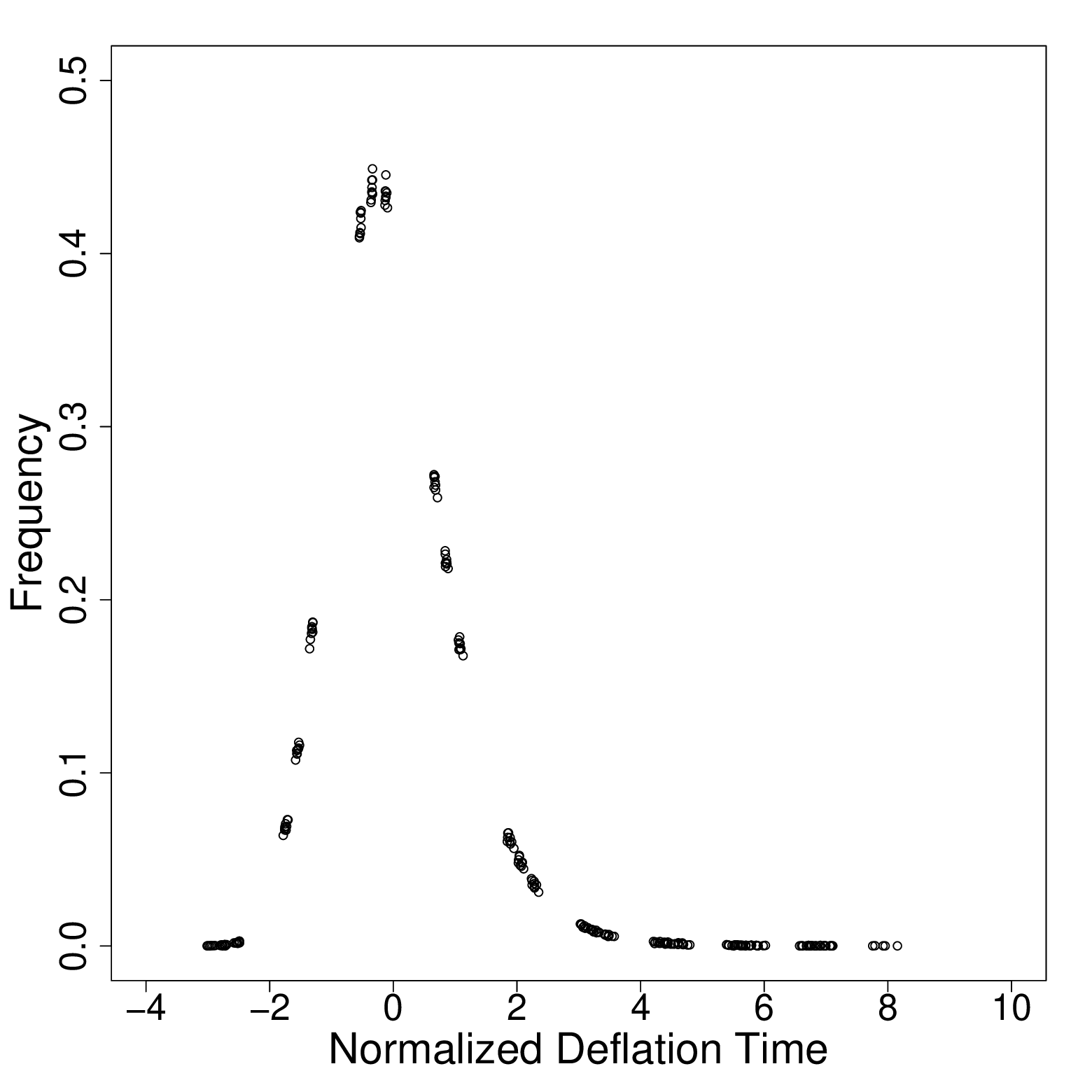}} 
\subfloat[][]{
\includegraphics[width=7cm]{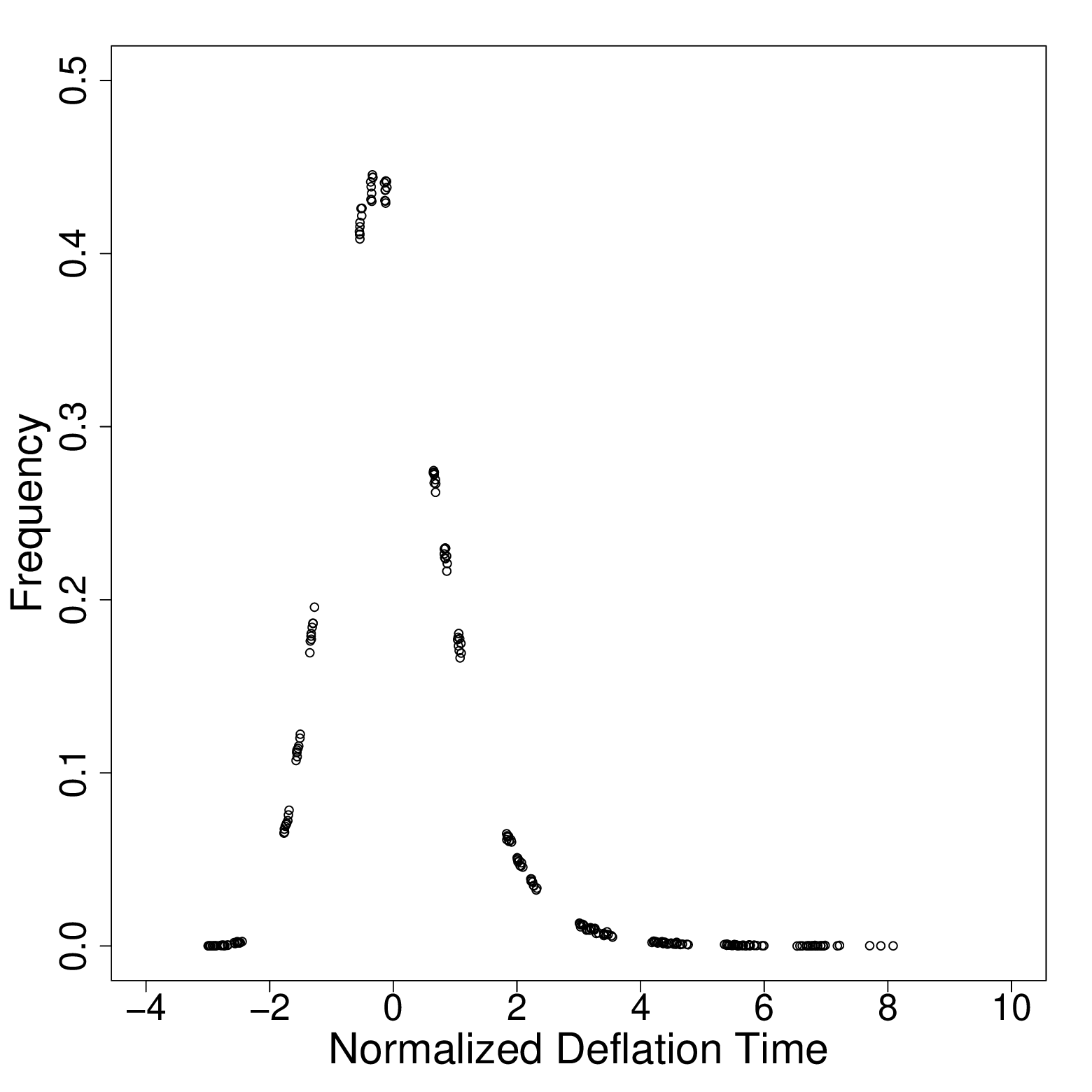}}
\caption{ {\bf Universality for the shifted QR algorithm on the Wigner  class.\/} Empirical deflation time for the QR algorithm with the Wilkinson shift normalized as in \qref{eq:uni-scaling} with $\epsilon=10^{-k}$, $k=8, 10, 12$ and $n$ ranging from $10,30,\dots,190$. Note that $\epsilon$ is significantly smaller than for the unshifted QR algorithm. The ensembles are (a) GOE; (b) Hermite-1; (c) Gaussian Wigner; and (d) Bernoulli. Each of the figures
(a), (b), (c), and (d) is obtained by collapsing the data as in Figure~\ref{fig:QRnormalized-runtimes-wigner}. The peak of the TE1 ensemble is lower, and the tail shorter, than those for the other three ensembles.}
\label{fig:S-QRnormalized-runtimes-wigner}
\end{figure}

\begin{figure}[h]
\subfloat[][]{
\includegraphics[width=7cm]{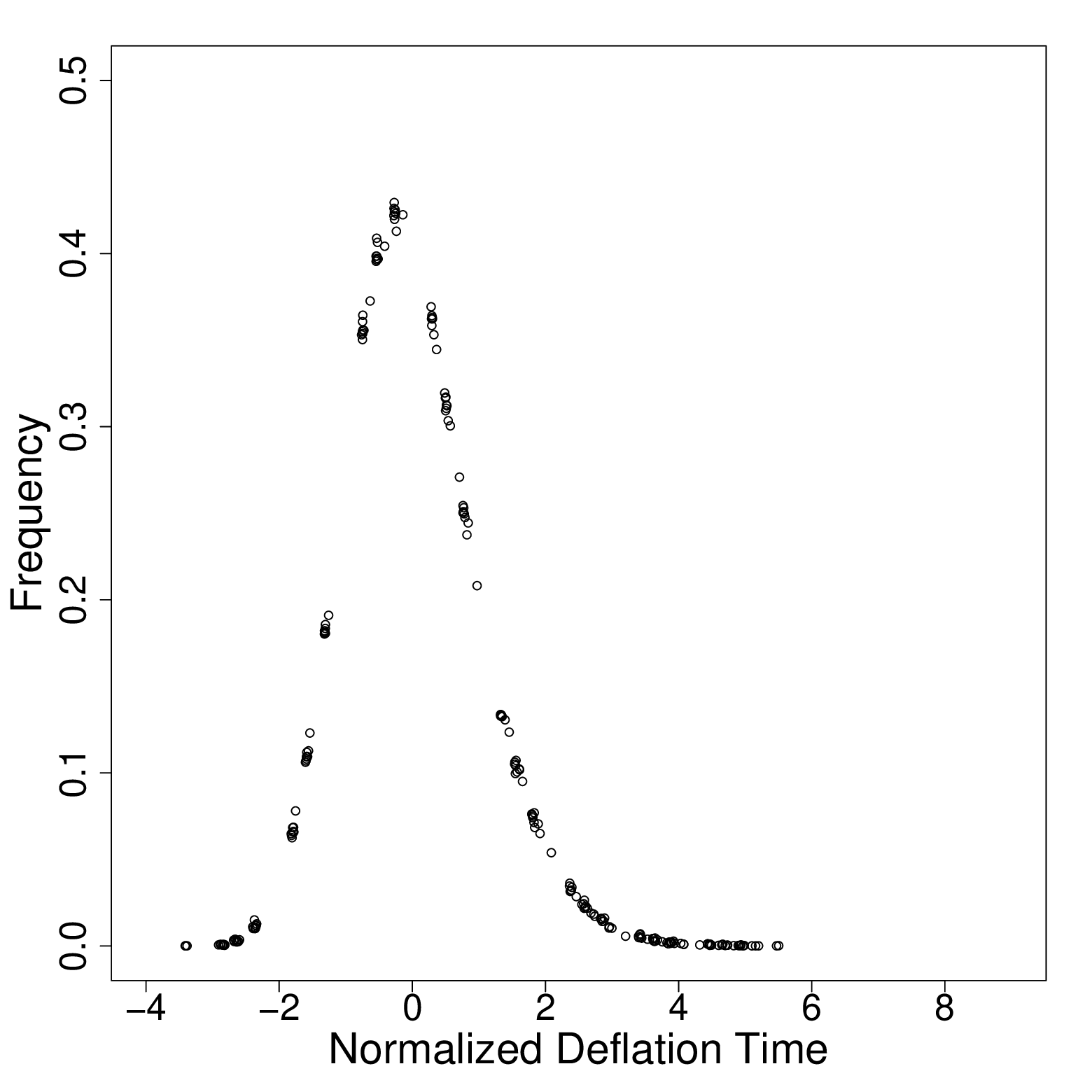}}
\subfloat[][]{
\includegraphics[width=7cm]{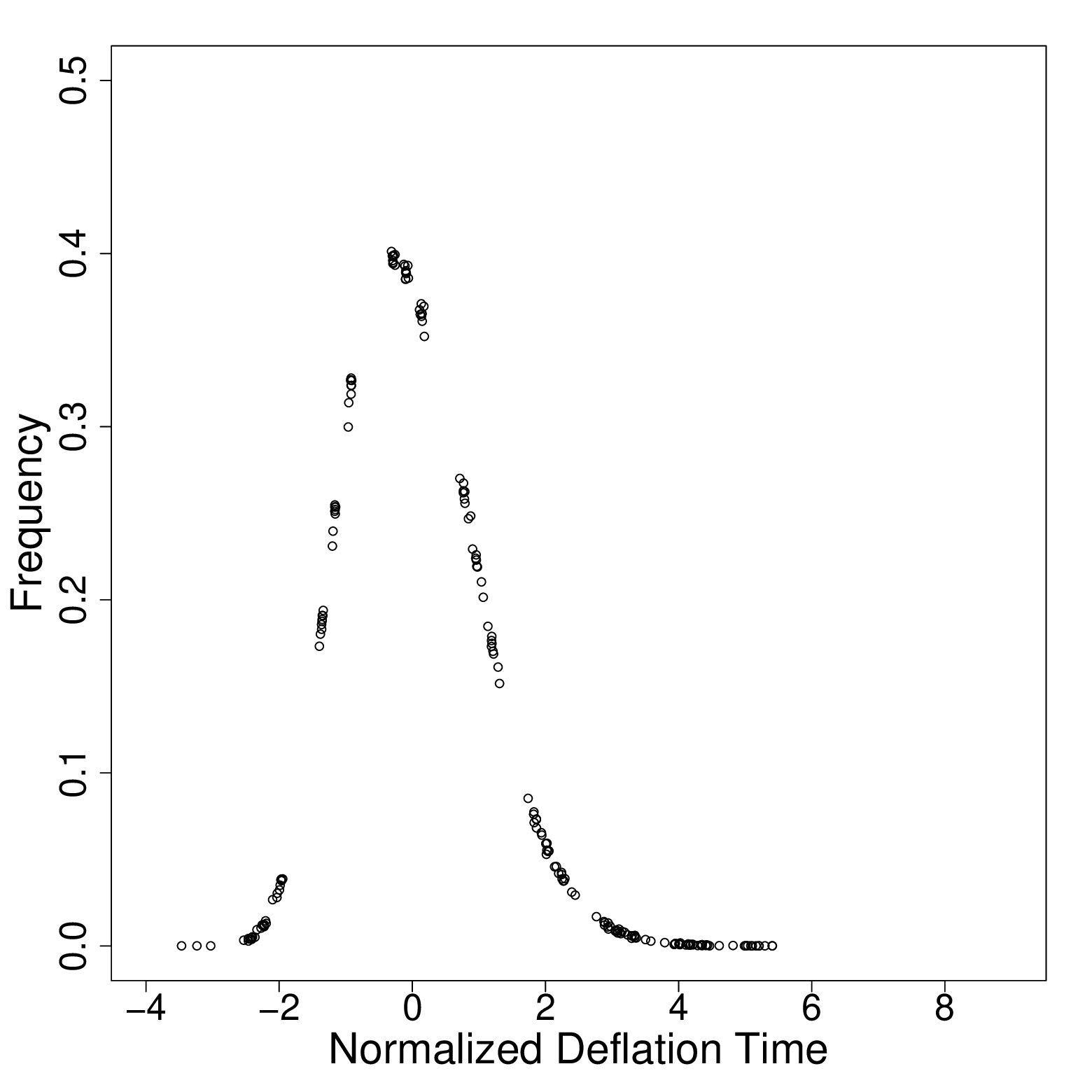}}
\caption{{\bf The shifted QR algorithm applied to non-Wigner ensembles.\/} 
Normalized empirical deflation time distributions for the QR  algorithm with Wilkinson shift for $\epsilon=10^{-k}$, $k=8, 10, 12$ and $n$ ranging from $10,30,\dots,190$.  The random matrix ensembles are (a) UDSJ and (b) JUE. Note that the results for these ensembles seem very similar to those for the Wigner class data shown in Figure~\ref{fig:S-QRnormalized-runtimes-wigner}. UDSJ is similar to the full matrix ensembles, while JUE is similar to TE1, which is also a tridiagonal ensemble.  }
\label{fig:S-QR-normalized-runtimes-nonwigner}
\end{figure}

\begin{figure}[h]
\subfloat[][]{
\includegraphics[width=7cm]{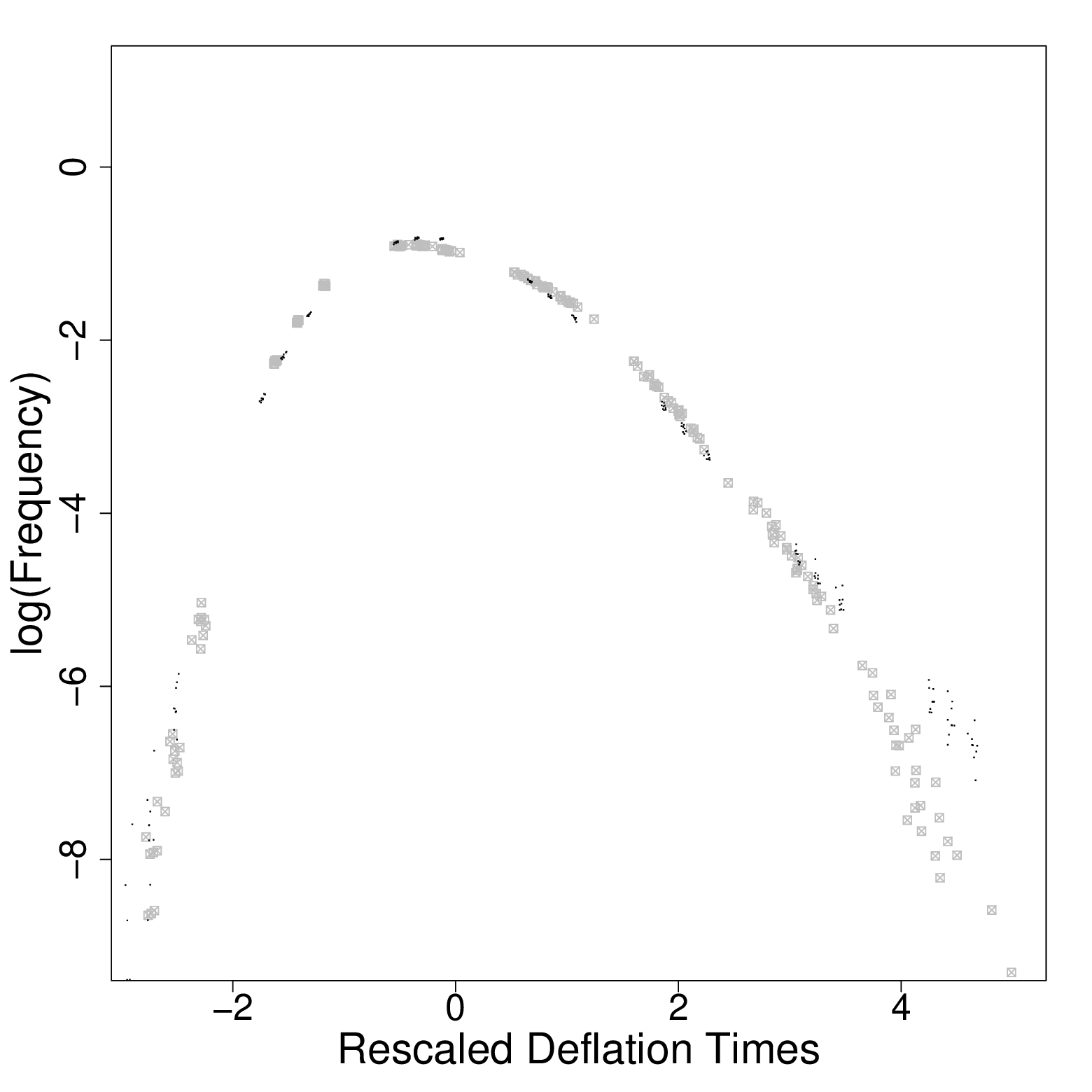}}
\subfloat[][]{
\includegraphics[width=7cm]{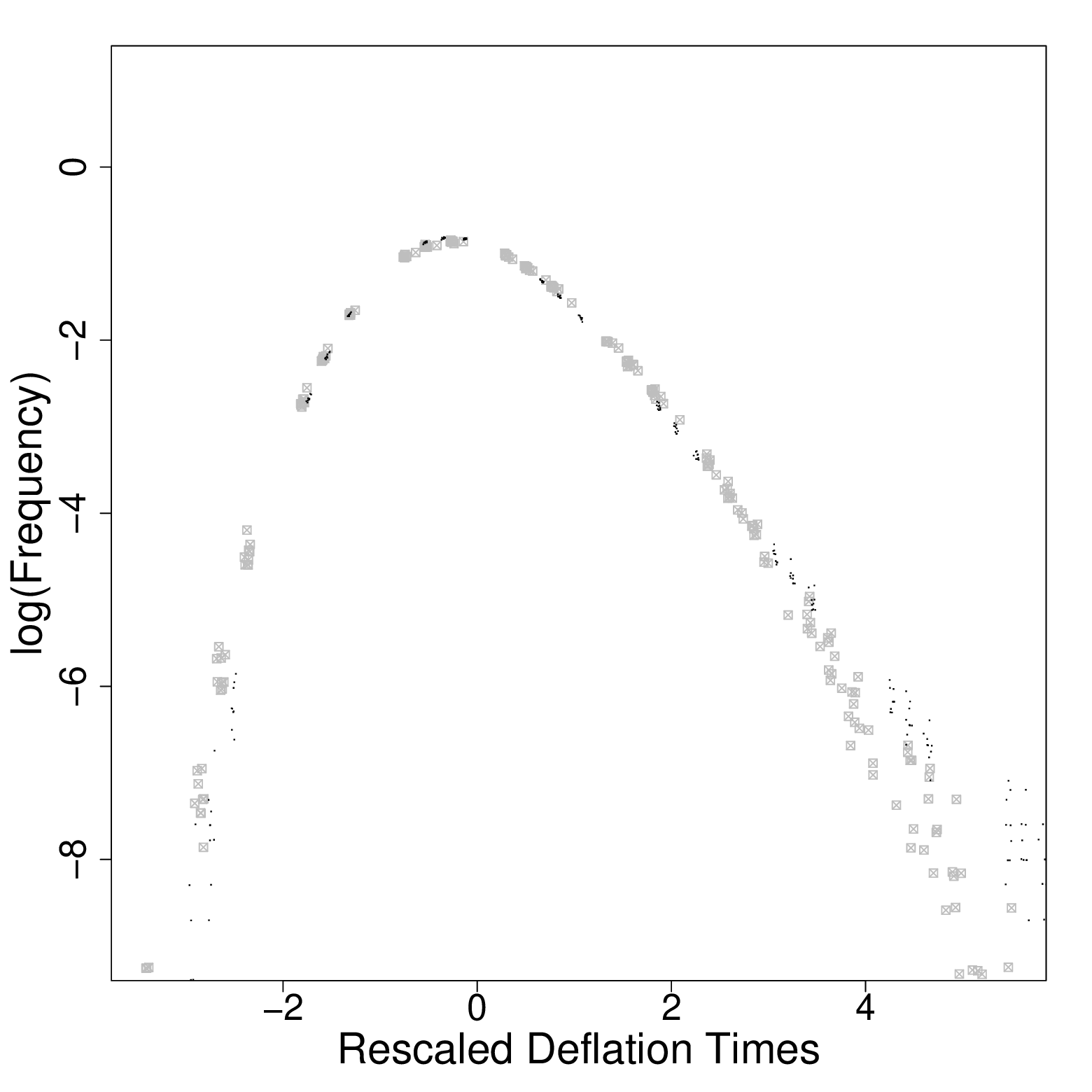}}\\
\subfloat[][]{
\includegraphics[width=7cm]{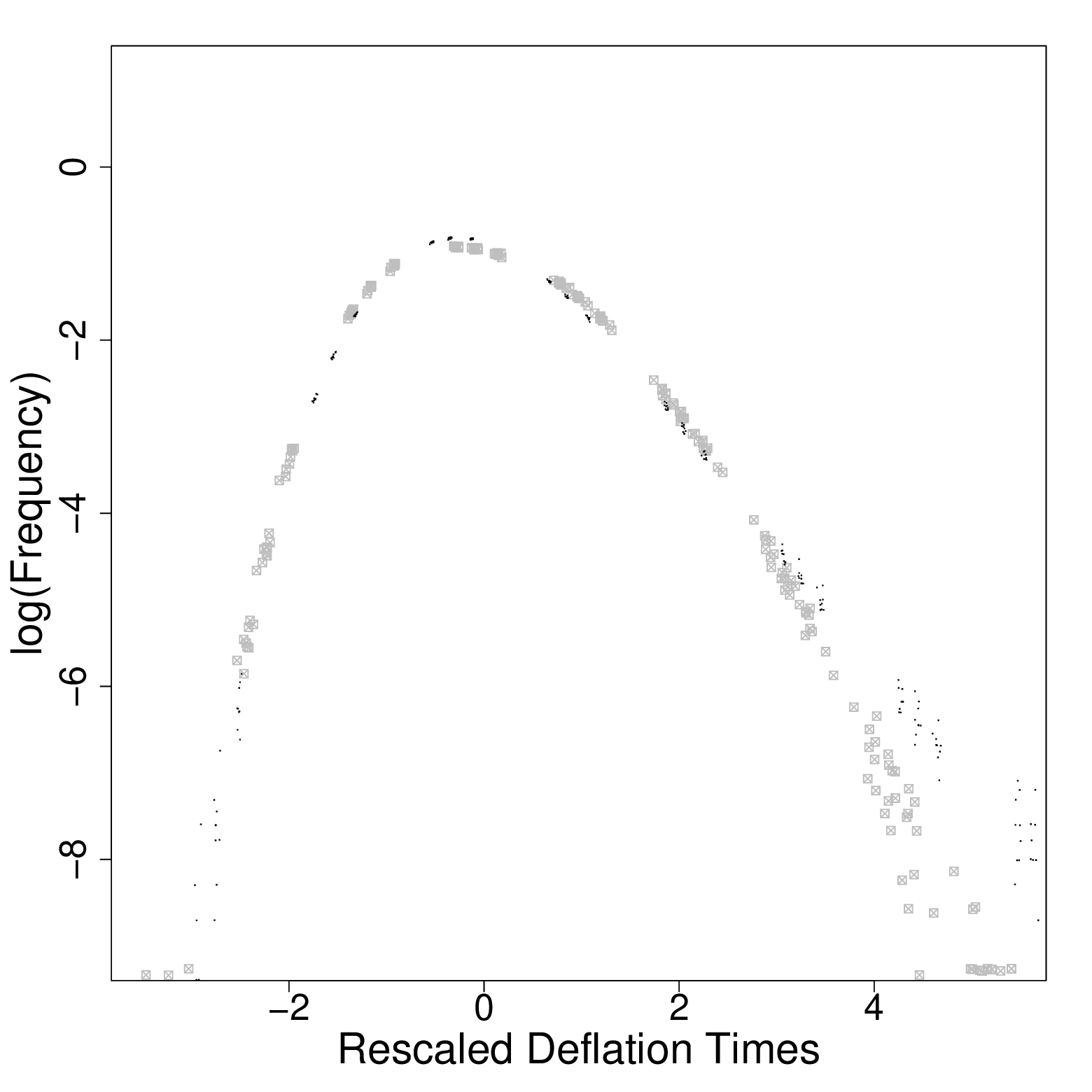}}
\caption{ {\bf Comparison of ensembles for shifted QR \/}.
Histograms of the normalized deflation time for the shifted  QR  algorithms on a logarithmic scale. In each figure GOE (black dots) is contrasted with data from a second ensemble (gray dots) (a)  {\bf GOE and TE1.\/}  Empirical normalized deflation time distributions from $40$ GOE histograms (black
dots) are contrasted with data from $40$ TE1 histograms (gray squares). (b) {\bf GOE and UDSJ.\/} (c) {\bf GOE and JUE.\/} 
\label{fig:S-combined-histograms}
}
\end{figure}

\begin{figure}[h]
\subfloat[][]{
\includegraphics[width=7cm]{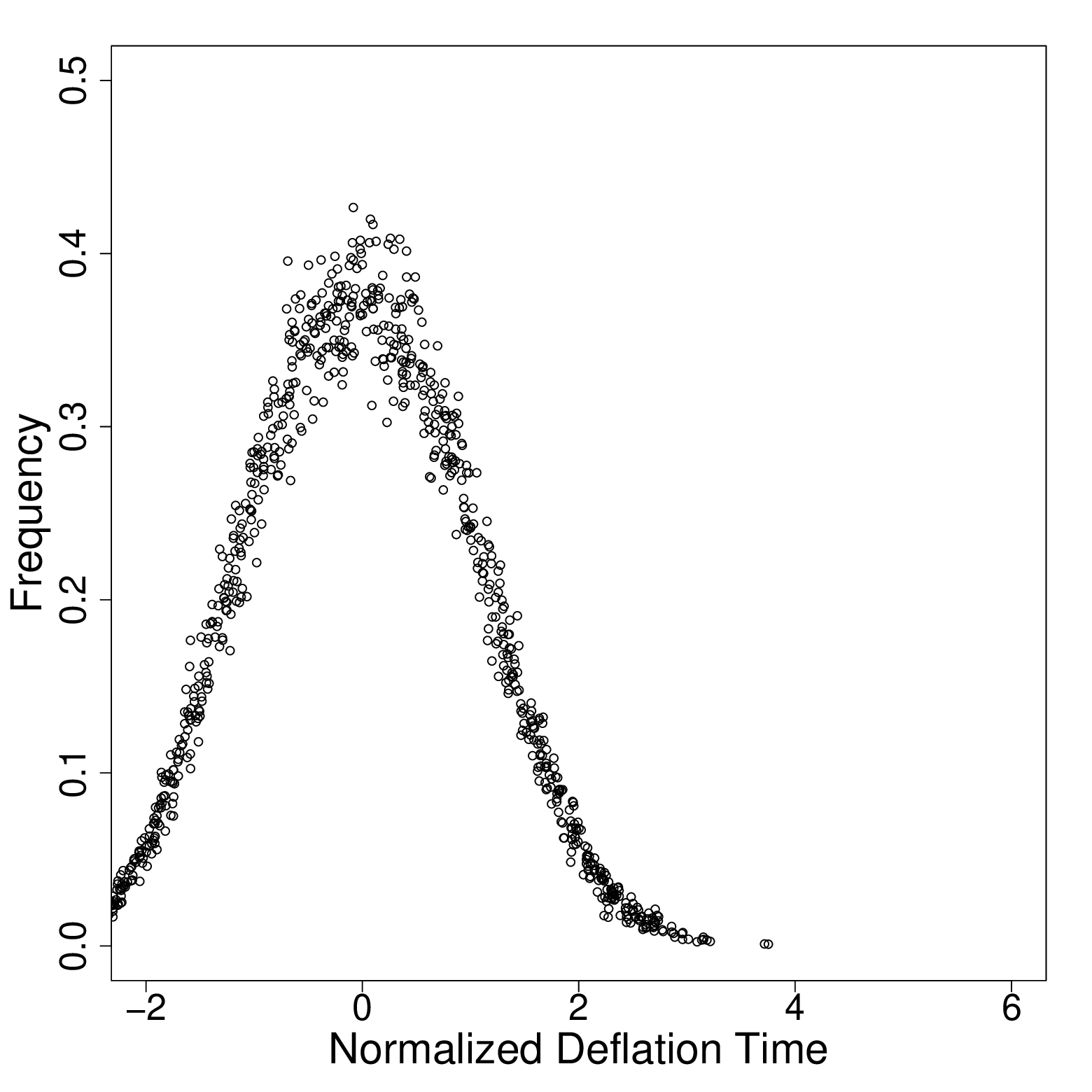}}
\subfloat[][]{ 
\includegraphics[width=7cm]{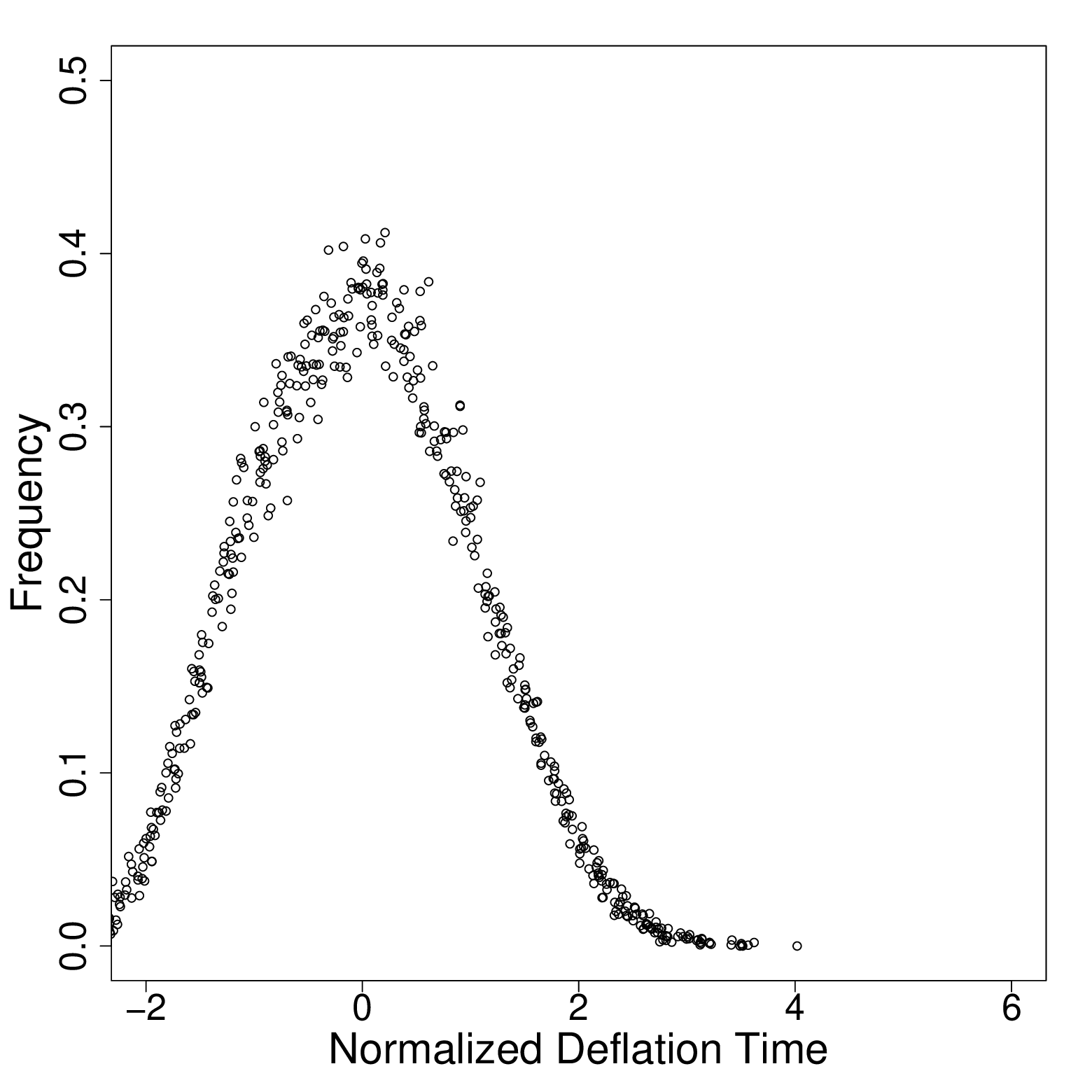}}\\
\subfloat[][]{
\includegraphics[width=7cm]{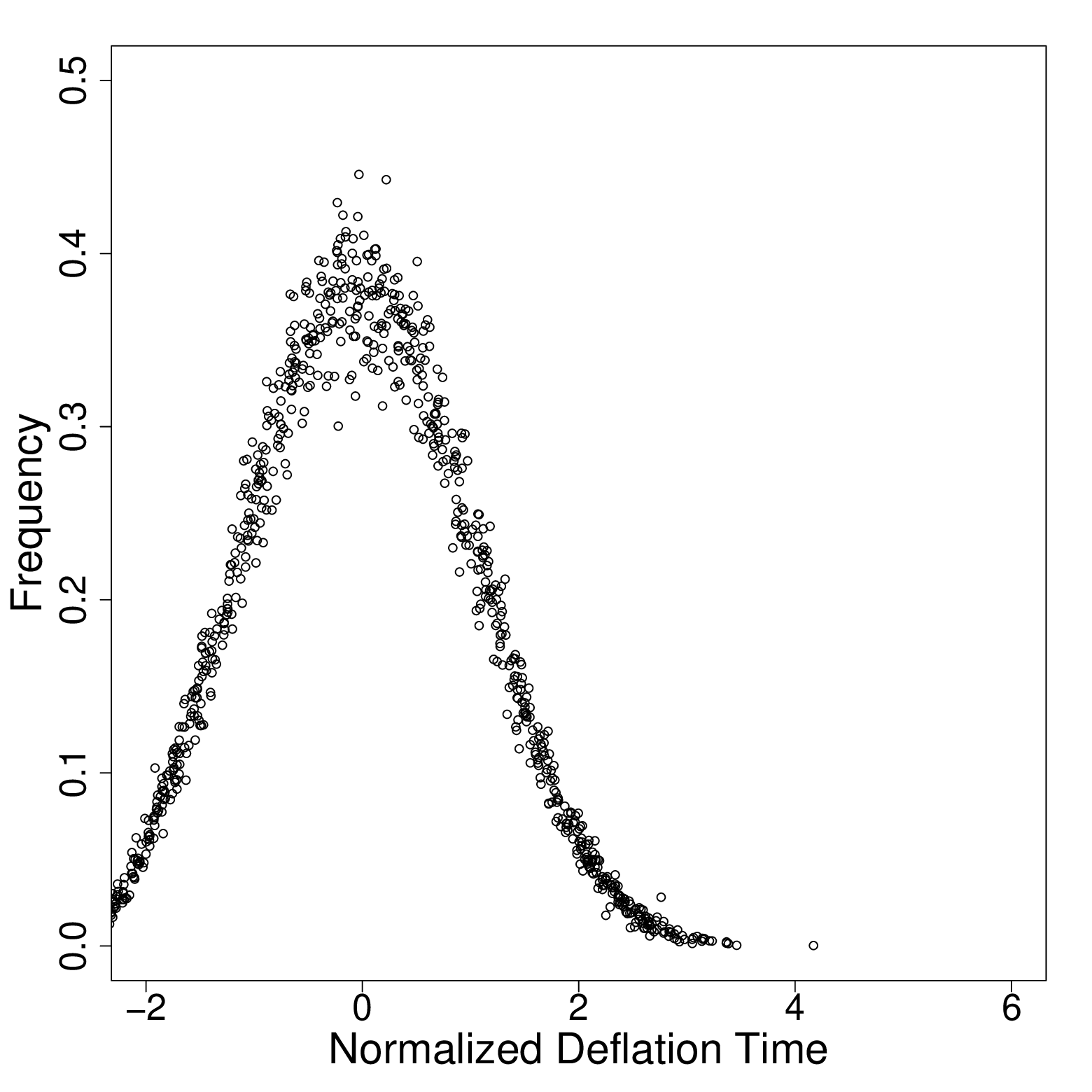}}
\subfloat[][]{
\includegraphics[width=7cm]{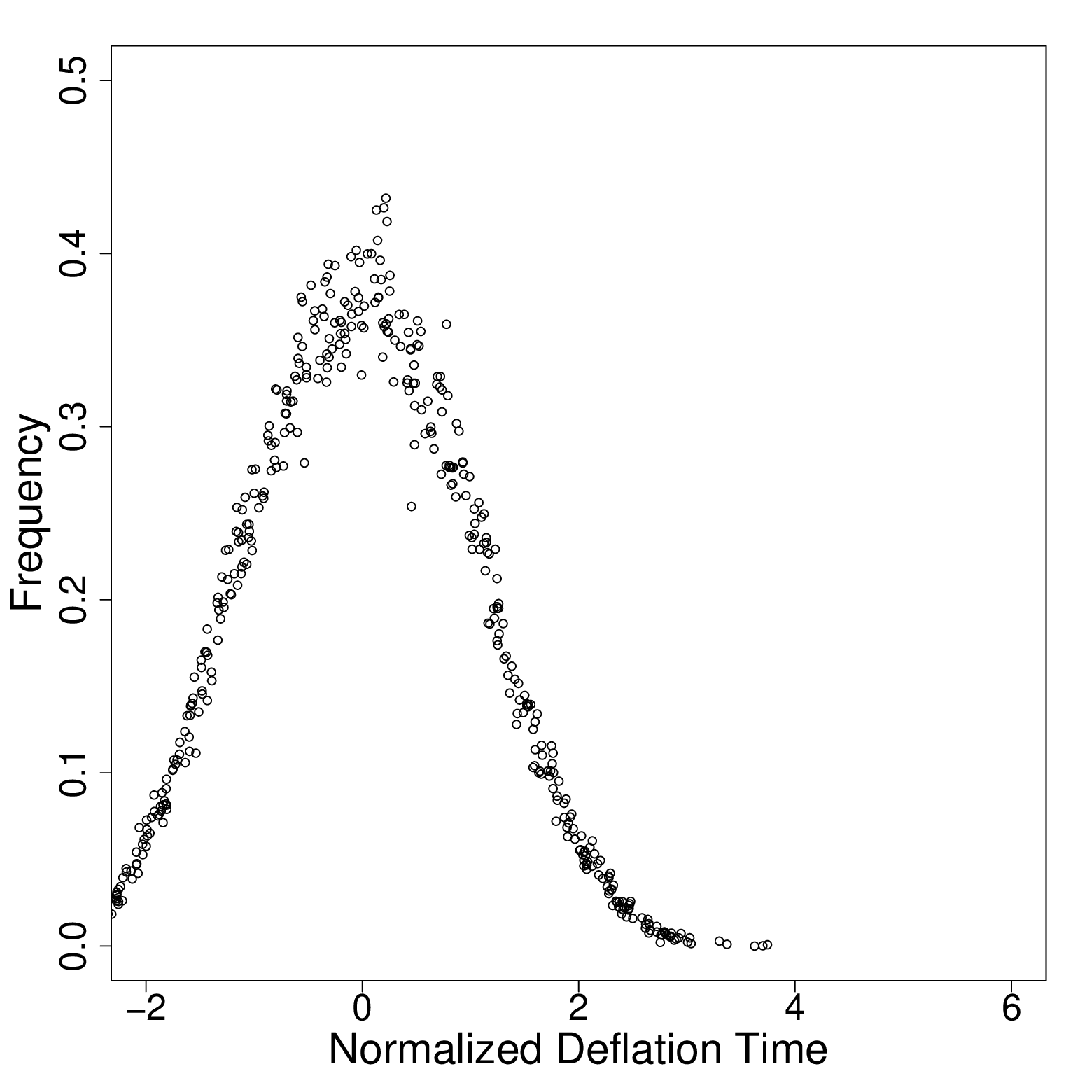}}
\caption{ {\bf Universal deflation time statistics for the Toda algorithm applied to the Wigner  class.\/} Empirical deflation time normalized as in \qref{eq:uni-scaling} for $\epsilon=10^{-k}$, $k=2,4,6,8$ and $n$ ranging from $10,30,\dots,190$. The random matrix ensembles are (a) GOE; (b) Hermite-1; (c) Gaussian Wigner; and (d) Bernoulli. Again each of the figures (a), (b), (c), and (d) is obtained by rescaling the data of $40$ fixed-$n$ and
fixed-$\epsilon$ histograms and plotting them together. All these data are observed to collapse onto one universal curve.  This universality is amplified in  Figure~\ref{fig:Toda-combined-histograms} below.
}
\label{fig:Todanormalized-runtimes-wigner}
\end{figure}

We applied this approach with the Gaussian, Exponential, Weibull and Gamma families.
We found that the exponential tails fit the QR runtime data especially well for 
small values of the deflation tolerance. The fit of the Toda runtime 
data to Gaussian tails is very compelling across most experimental regimes.
Direct pictorial comparisons of the normalized Toda runtimes with the standard normal
as well as normalized QR runtimes with normalized Gamma distributions are shown in
Figure \ref{fig:combined-histograms}. Further details of the statistical tests may be found in~\cite{ChristianDiss}.


\begin{figure}
\subfloat[][]{
\includegraphics[width=7cm]{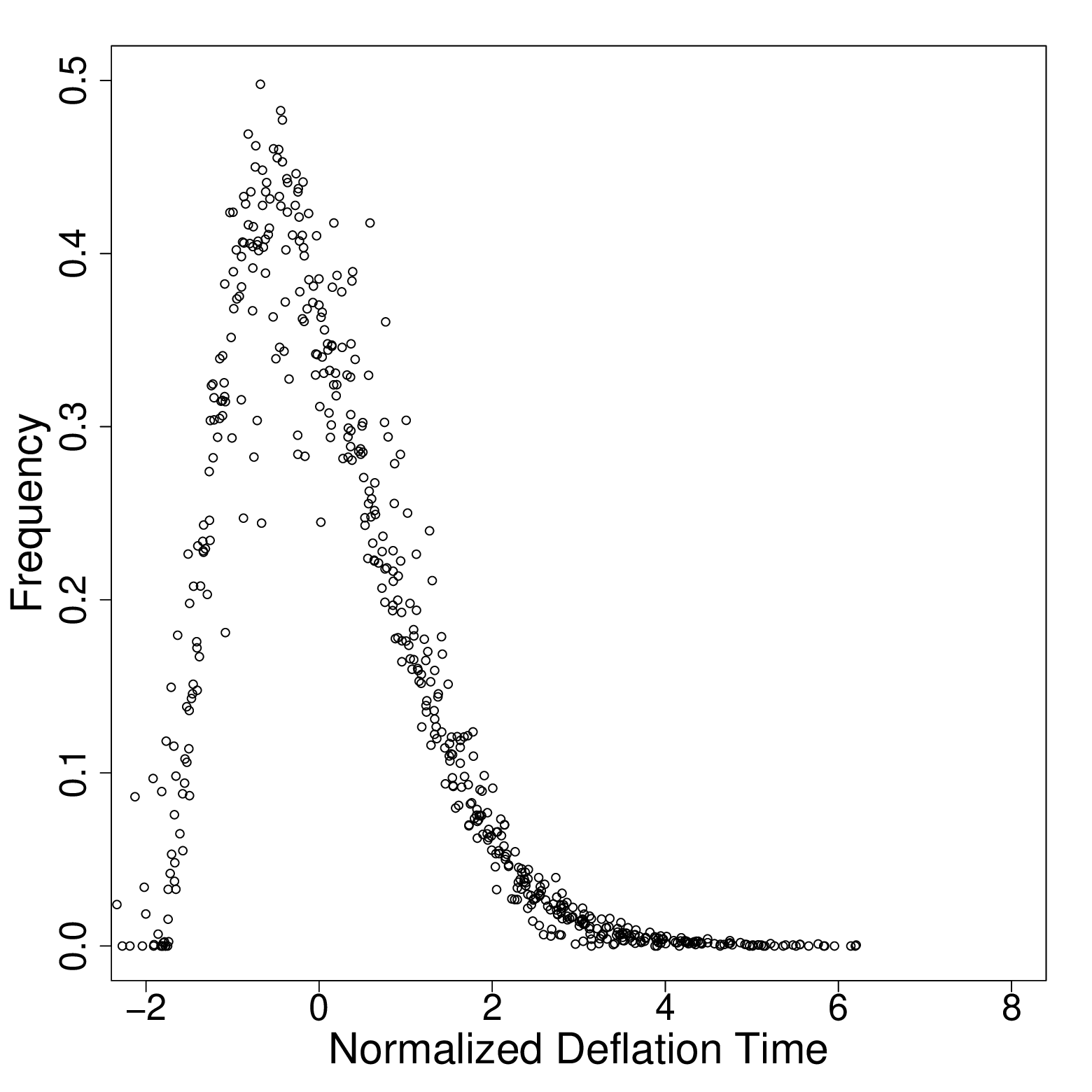}}
\subfloat[][]{
\includegraphics[width=7cm]{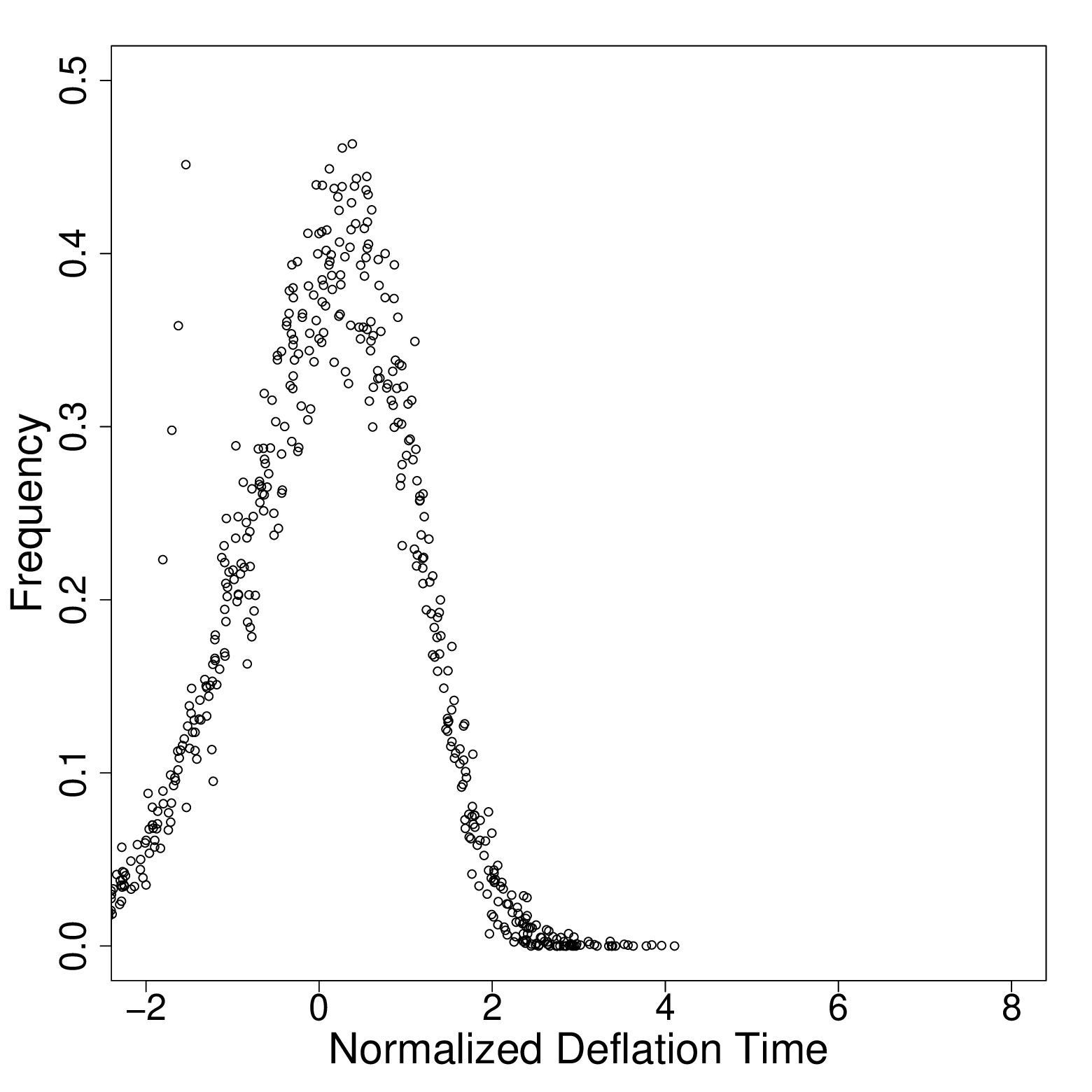}}
\caption{{\bf The Toda algorithm applied to  non-Wigner ensembles.\/} 
Normalized empirical deflation time distributions for the Toda algorithm with $\epsilon=10^{-k}$, $k=2,4,6,8$ and $n$ ranging from $10,30,\dots,190$.  The random matrix ensembles are (a) UDSJ  and  (b) JUE. Each figure contains the normalized empirical data of $40$ fixed-$n$ and fixed-$\eps$ histograms. All these data are observed to collapse onto a single curve. However, these curves are not the same for UDSJ and JUE and neither of these coincides with the curve for Wigner data shown in 
Figure~\ref{fig:Todanormalized-runtimes-wigner}.
	\label{fig:Toda-normalized-runtimes-nonwigner}}
\end{figure}


\begin{figure}[h]
\subfloat[][]{
\includegraphics[width=7cm]{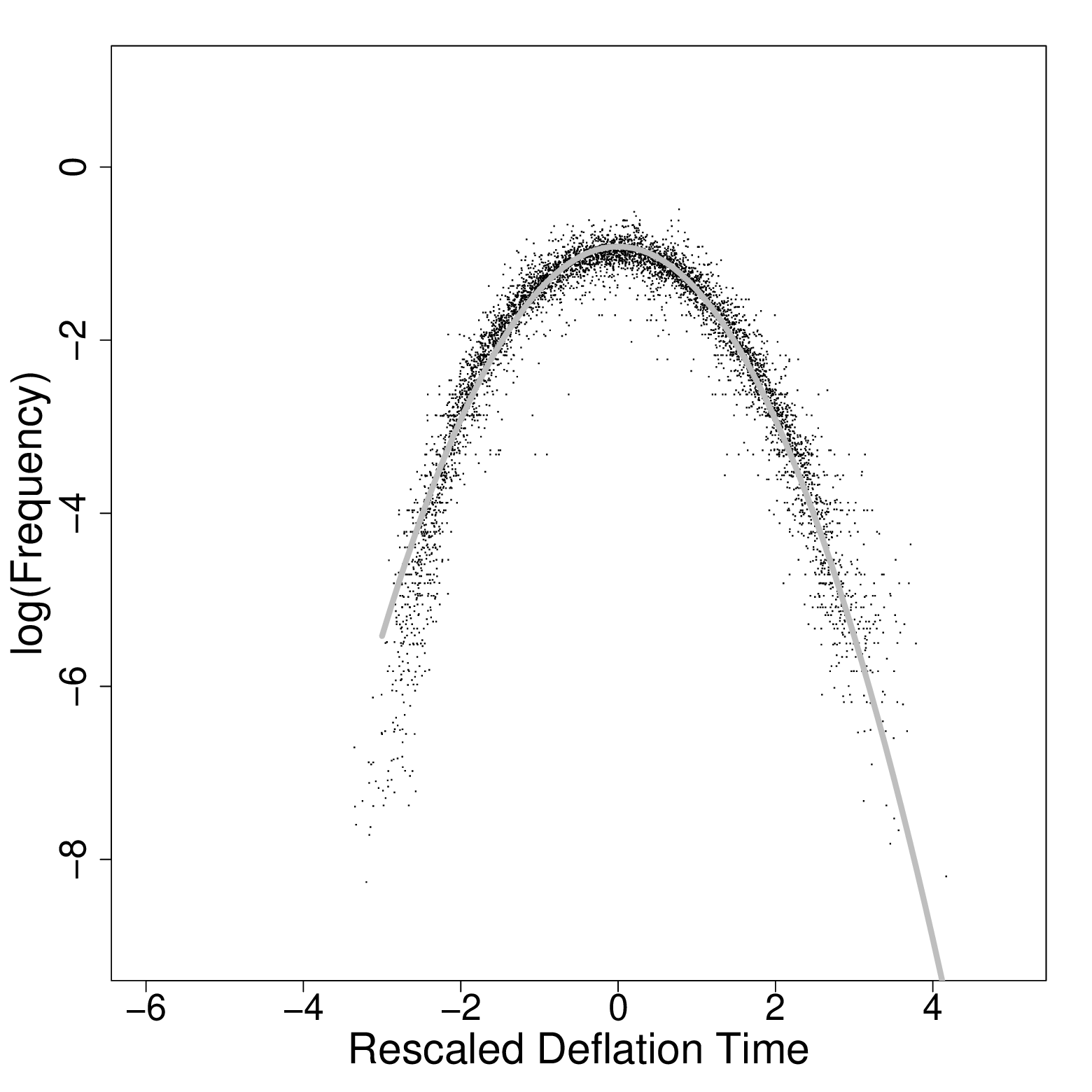}}
\subfloat[][]{
\includegraphics[width=7cm]{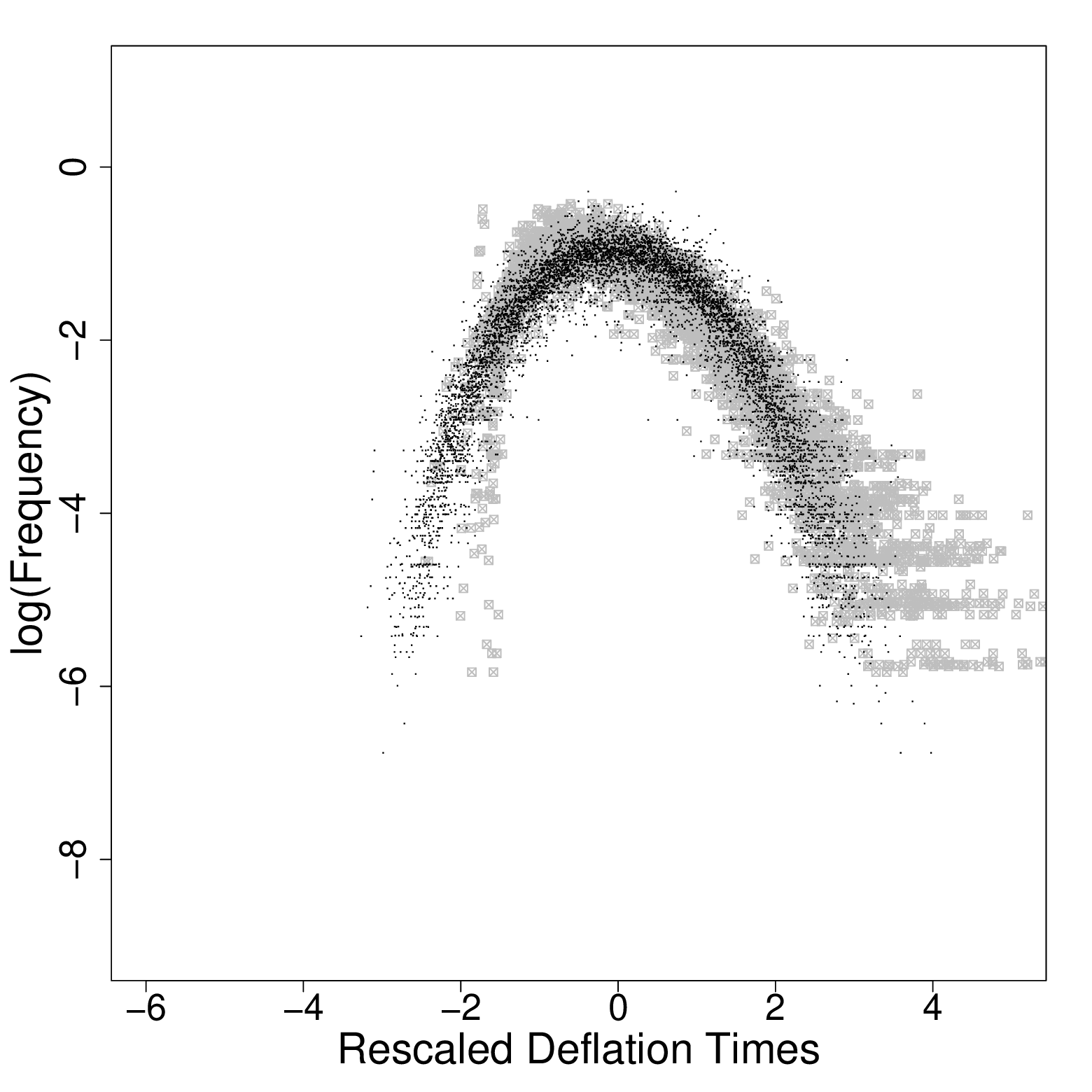}}
\caption{ {\bf Gaussian tail for the Toda algorithm\/}.
Histograms of the normalized deflation time for the  QR  algorithms on a logarithmic scale. (a) {\bf Wigner data. \/} Empirical normalized deflation time distributions 
from all $160$ histograms of Wigner class initial data (black dots) are compared with a standard normal distribution (gray line). 
 (b) {\bf non-Wigner data.\/}  Empirical normalized deflation time distributions from $40$ GOE histograms (black
dots) is contrasted with data from $40$ UDSJ histograms (gray squares). 
}
\label{fig:Toda-combined-histograms}
\end{figure}


\subsection{The dependence of $\mu_{n,\eps}$ and $\sigma_{n,\eps}$ on $n$ and $\eps$}
\label{subsec:means}
\begin{figure}[h]
\subfloat[][]{
\includegraphics[width=7cm]{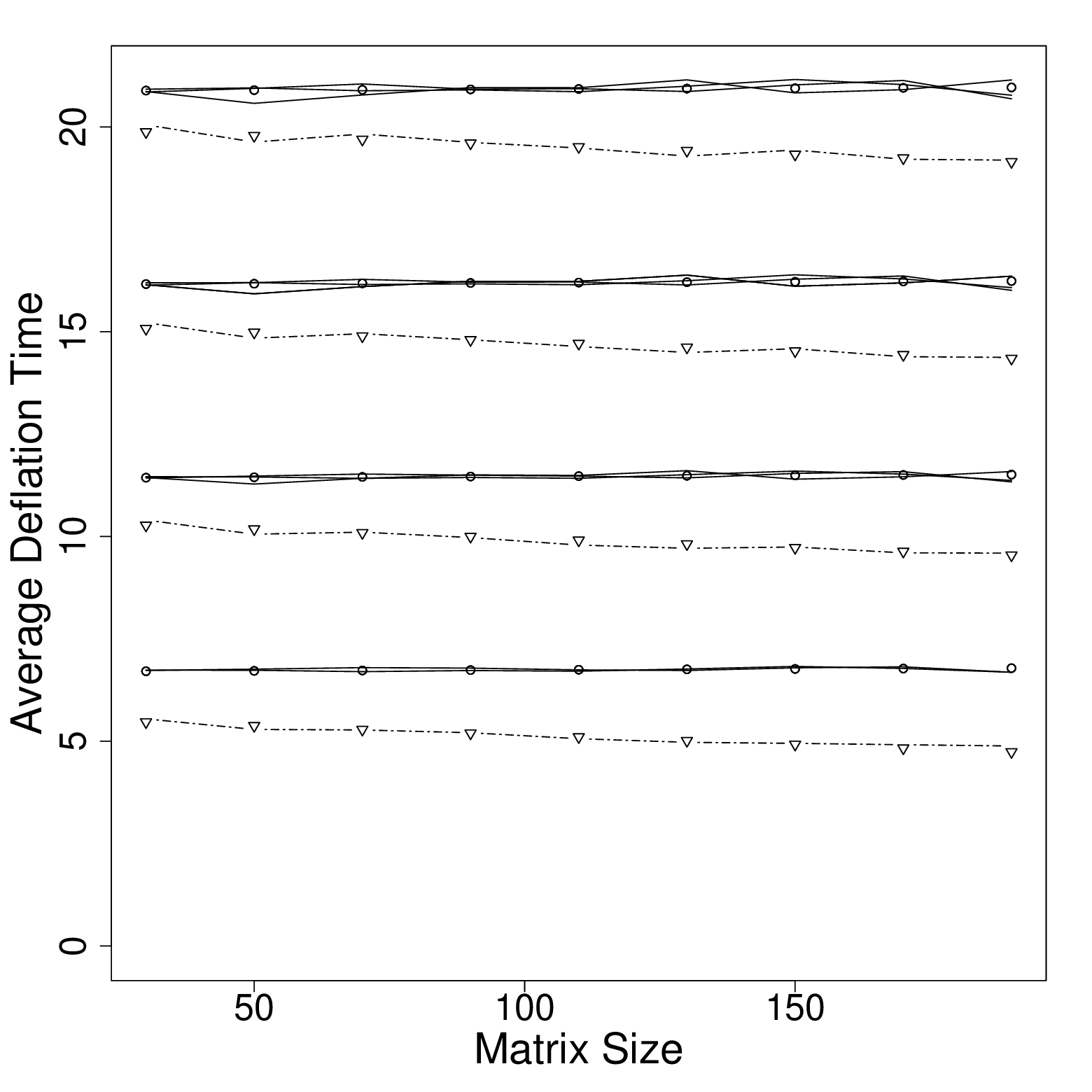}}
\subfloat[][]{
\includegraphics[width=7cm]{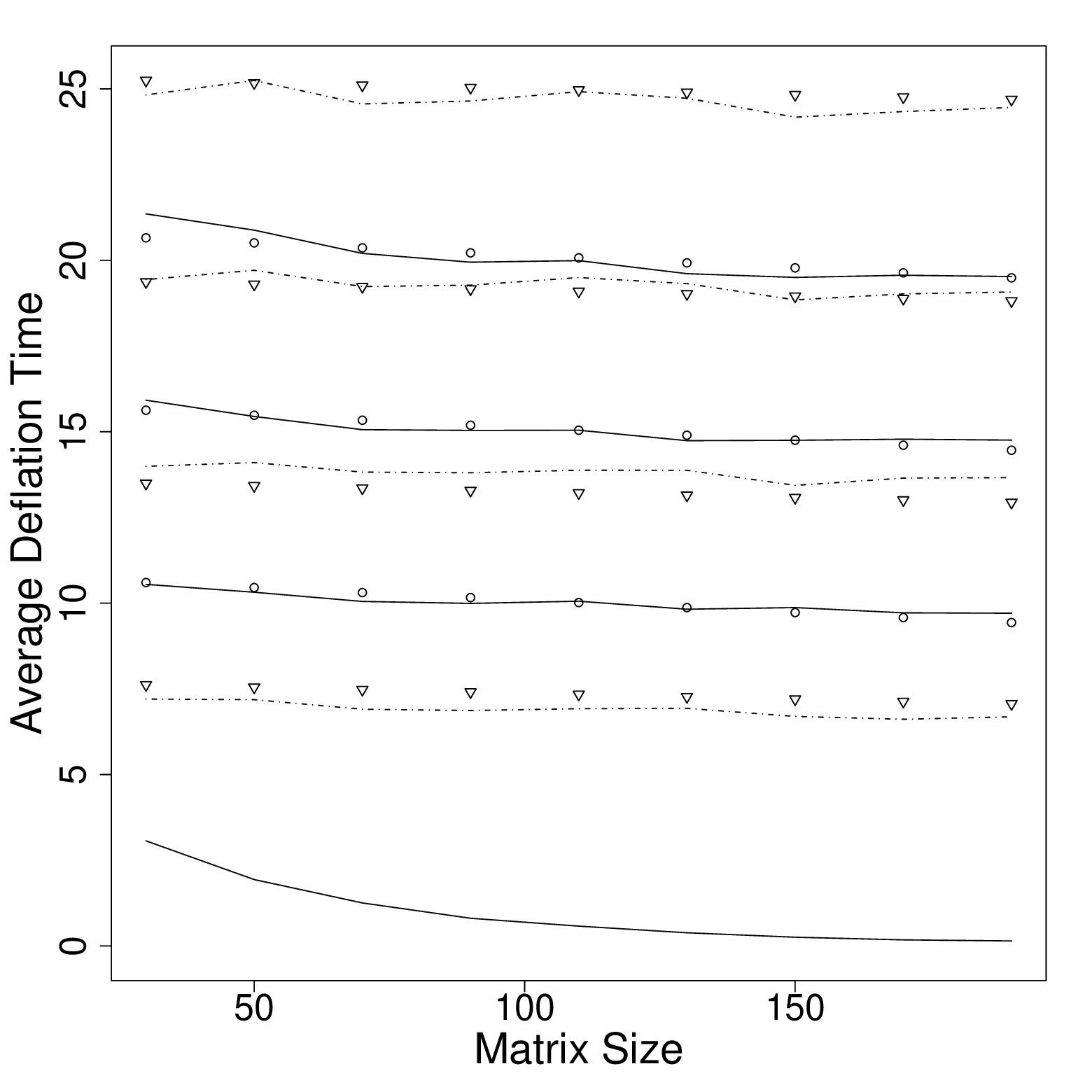}}
\caption{{\bf Mean deflation time $\mu_{n,\eps}$ for the QR algorithm.\/} 
Empirical average of the deflation time for $\epsilon=10^{-k}$, $k=2,4,6,8$ and $n$ in the range $10,30,\dots,190$.
(a) {\bf Wigner class initial data.\/} The full lines are the empirical mean $\mu_{n,\eps}$ for GOE, Gaussian Wigner and Bernoulli ensembles. Note that they seem to align well with one another. The circles are the values obtained from the regression estimate  \qref{eq:qr-mean-regression} with the parameters listed in Table~\ref{tab:means-regression-qr}. The dashed line and triangles represent empirical data and the regression respectively for the Hermite--1 ensemble.
(b) {\bf JUE and UDSJ initial data.\/} The full line and dashed line are the empirical mean $\mu_{n,\eps}$ for UDSJ and JUE data respectively. The circles and triangles are the regression estimates  for UDSJ and JUE respectively. 
As $\eps$ decreases the curves move up monotonically. The regression is not applied to the lowest curve in (b) since $\eps=0.01$ is sufficiently large that several matrices deflate instantaneously.
}
\label{fig:QR-means}
\end{figure}

\begin{figure}[h]
\subfloat[][]{
\includegraphics[width=7cm]{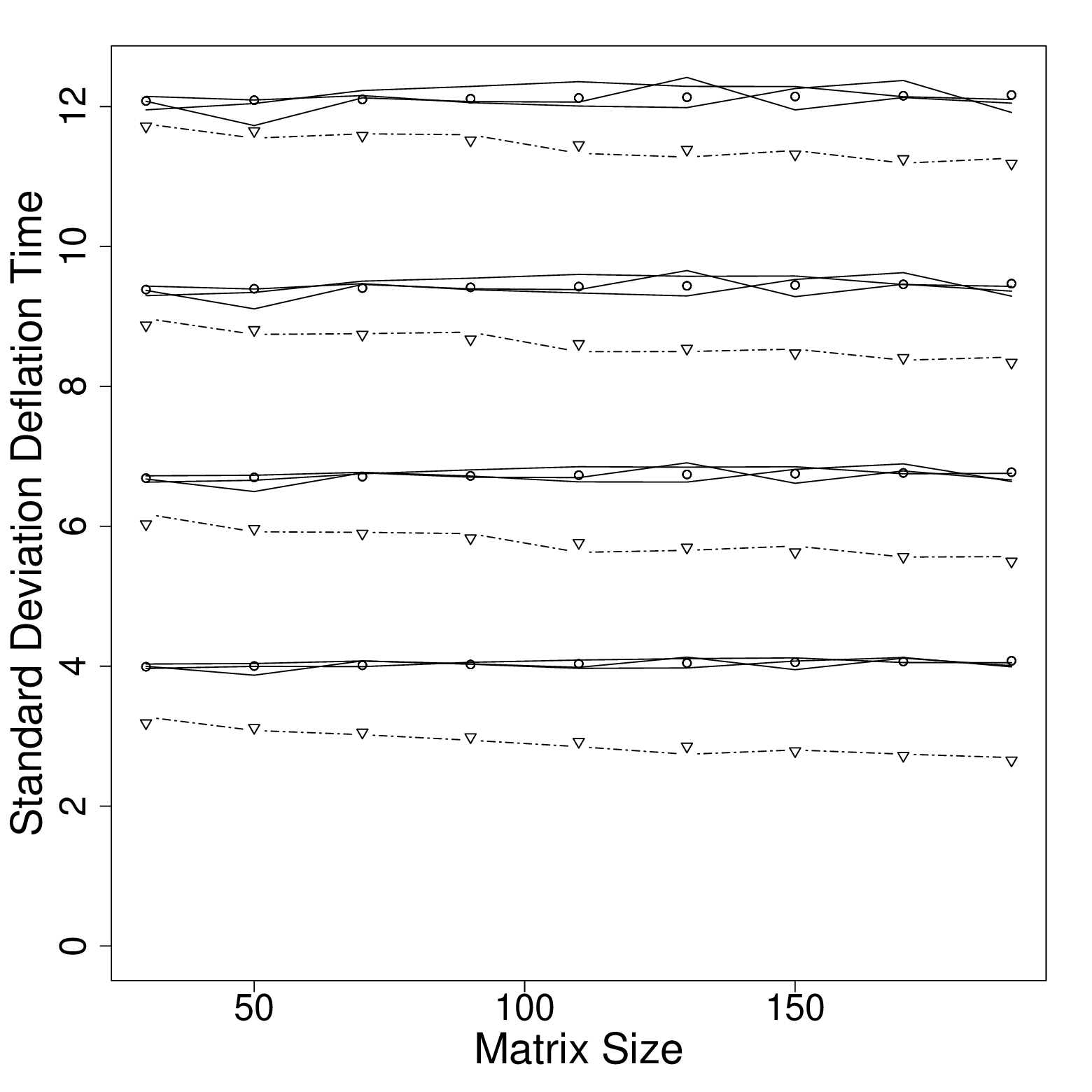}}
\subfloat[][]{
\includegraphics[width=7cm]{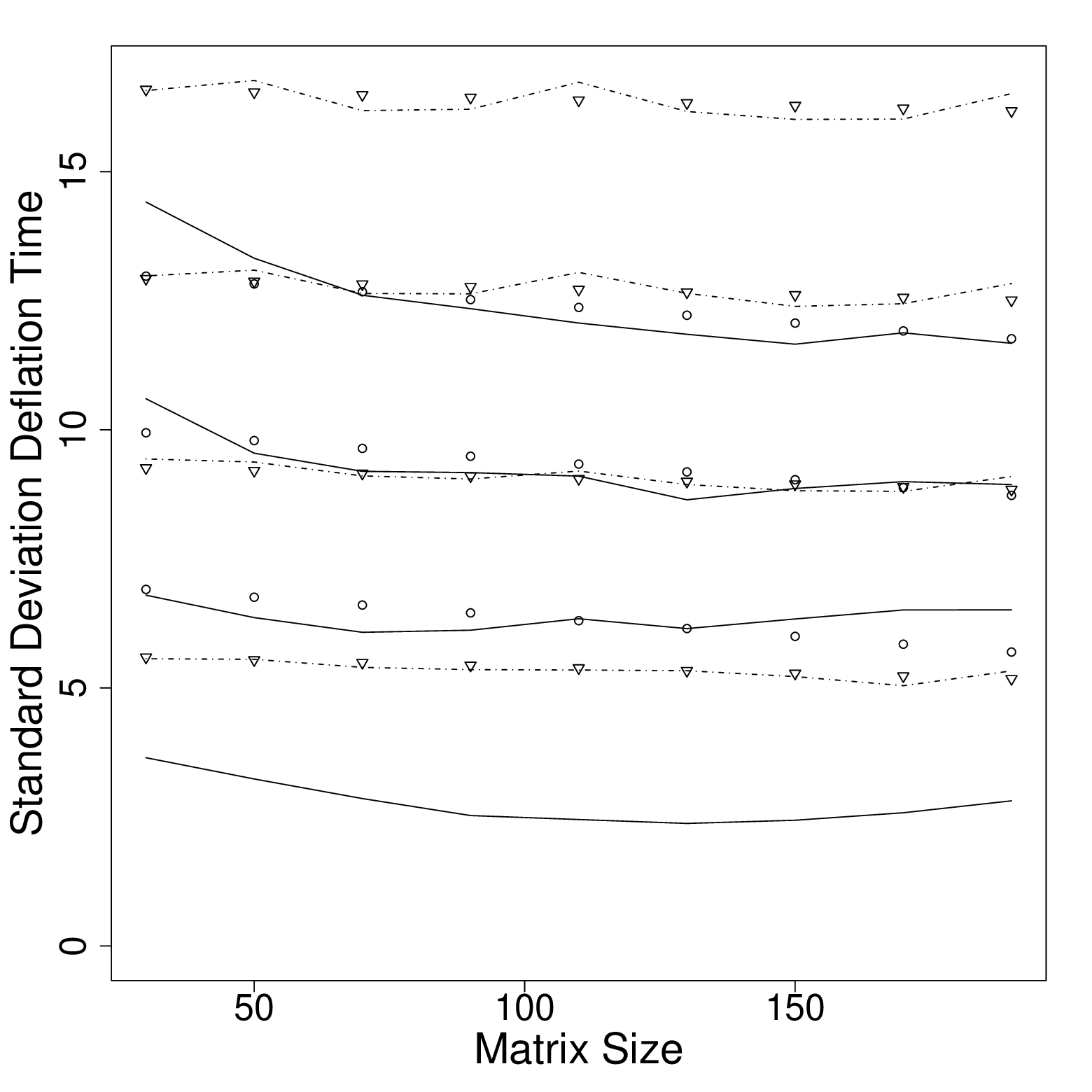}}
\caption{{\bf Standard deviation $\sigma_{n,\eps}$ of the deflation time for the QR algorithm.\/} (a) Ensembles in the Wigner class. (b) JUE and UDSJ. Legend as in Figure~\ref{fig:QR-means} with regression parameters from Table~\ref{tab:sdev-regression-qr}. }
\label{fig:QR-sd}
\end{figure}

We used linear regression to express $\mu_{n,\eps}$ and $\sigma_{n,\eps}$ as  functions of $\log \eps$ and $n$. Only the best fits are reported here. The data for the QR algorithm was matched very well by 
\begin{align} \label{eq:qr-mean-regression}
\mu_{n,\epsilon}&\approx a_0 + a_1 n + a_2 \log\epsilon,\\ \label{eq:qr-sdev-regression}
\sigma_{n,\epsilon}&\approx b_0 + b_1 n + b_2 \log\epsilon.
\end{align}
This regression is compared visually with the numerical data in Figure~\ref{fig:QR-means} and Figure~\ref{fig:QR-sd}. The regression parameters are tabulated in Table~\ref{tab:means-regression-qr} and Table~\ref{tab:sdev-regression-qr}. Since the means and variances do not visually appear to depend on $n$ for ensembles (1)--(3) we have also included the $p$--values for the $t$--test of the hypothesis that the coefficient corresponding to the dimension is zero. Note that $\mu_{n,\eps}$ and $\sigma_{n,\eps}$  are almost identical for the ensembles (1)--(3) in the Wigner class, while for the Hermite--1 initial data both statistics have a slightly larger value.

\begin{table}
\caption{Regression parameters for $\mu_{n,\eps}$ for the unshifted QR algorithm}
\label{tab:means-regression-qr}
\begin{center}
\begin{tabular}{|l||l|l|l|l|}
\hline 
Ensemble & $a_0$ & $a_1$ & $a_2$ & P--value $a_1$\\ \hline
GOE, Gaussian Wigner, Bernoulli   & 1.96824 & 0.0004690  & -1.0263649 &  0.0095\\ \hline
Hermite--1 & .802338 &-.004554 & -1.042907 &   $<2\cdot 10^{-16}$\\ \hline
UDSJ  & 0.7648330 & -0.0072921 & -1.0916354 & $4.16\cdot 10^{-8}$ \\ \hline
JUE  &1.844126 & -0.003467 & -1.276037 & 0.0256\\ \hline
\hline
\end{tabular}
\end{center}
\end{table}
\begin{table}
\caption{Regression parameters for $\sigma_{n,\eps}$ for the unshifted QR algorithm}
\label{tab:sdev-regression-qr}
\begin{center}
\begin{tabular}{|l||l|l|l|l|}
\hline 
Ensemble & $b_0$ & $b_1$ & $b_2$ & P--value $b_1$\\ \hline
GOE, Gaussian Wigner, Bernoulli & 1.2799509   &0.0005311& -0.5854859 & 0.0118  \\ \hline
Hermite--1 &0.442622 &-.003329 & -.617517 &$8\cdot 10^{-15}$  \\ \hline
UDSJ  & 1.066713& -0.007584& -0.658920 & 0.000353 \\ \hline
JUE  & 2.0044243 & -0.0026034& -0.7961700 & 0.000185  \\ \hline
\hline
\end{tabular}
\end{center}
\end{table}

\begin{figure}[h]
\subfloat[][]{
\includegraphics[width=7cm]{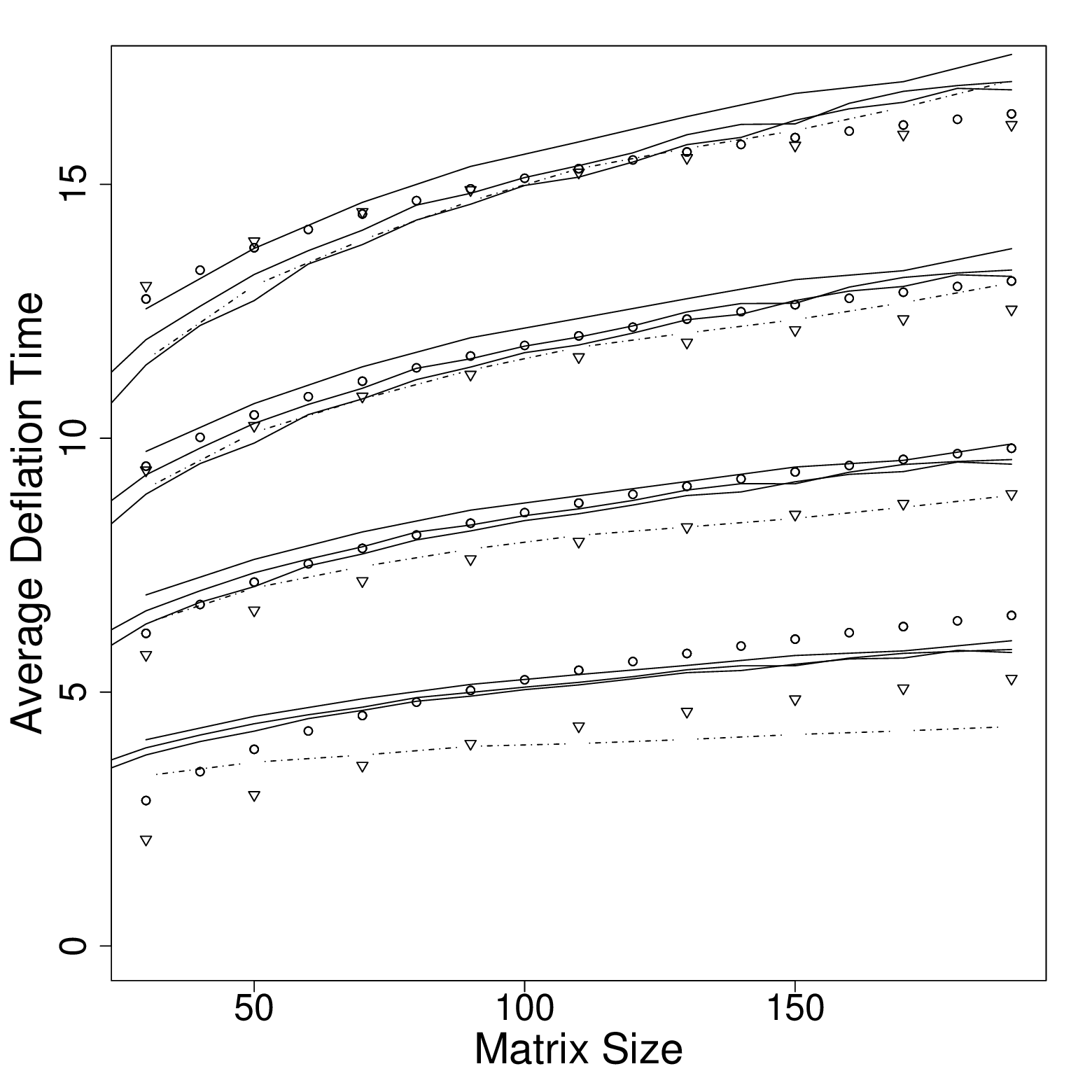}}
\subfloat[][]{
\includegraphics[width=7cm]{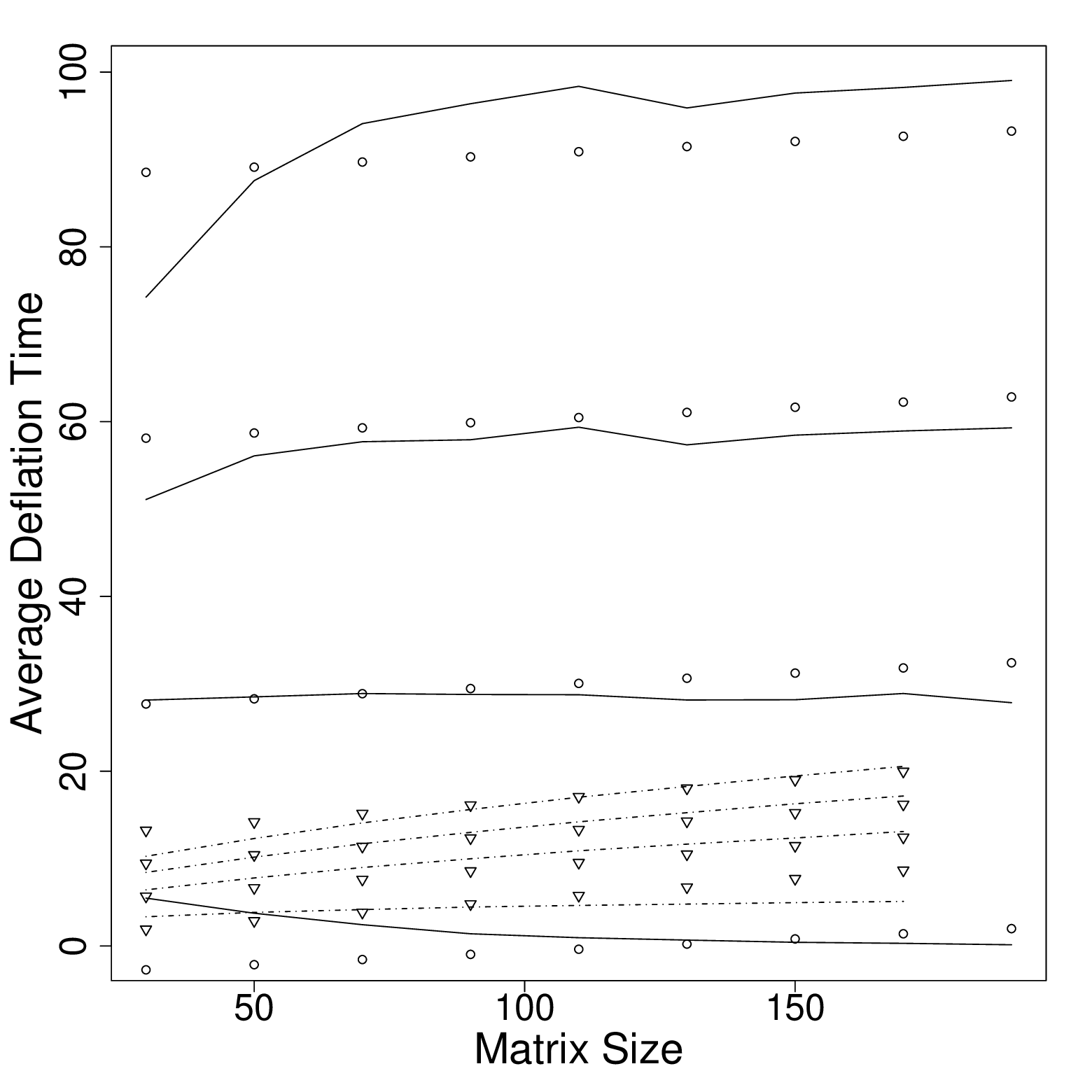}}
\caption{{\bf Mean deflation time $\mu_{n,\eps}$ for the Toda algorithm.\/} Mean deflation time distributions of the Toda algorithm for initial data described in 
Figure~\ref{fig:QR-means}. (a) Wigner class initial data and 
(b) JUE and UDSJ initial data.  Legend as in Figure~\ref{fig:QR-means} and regression parameters as in Table~\ref{tab:mean-regression-toda}.}
\label{fig:Toda-means}
\end{figure}

\begin{figure}[h]
\subfloat[][]{
\includegraphics[width=7cm]{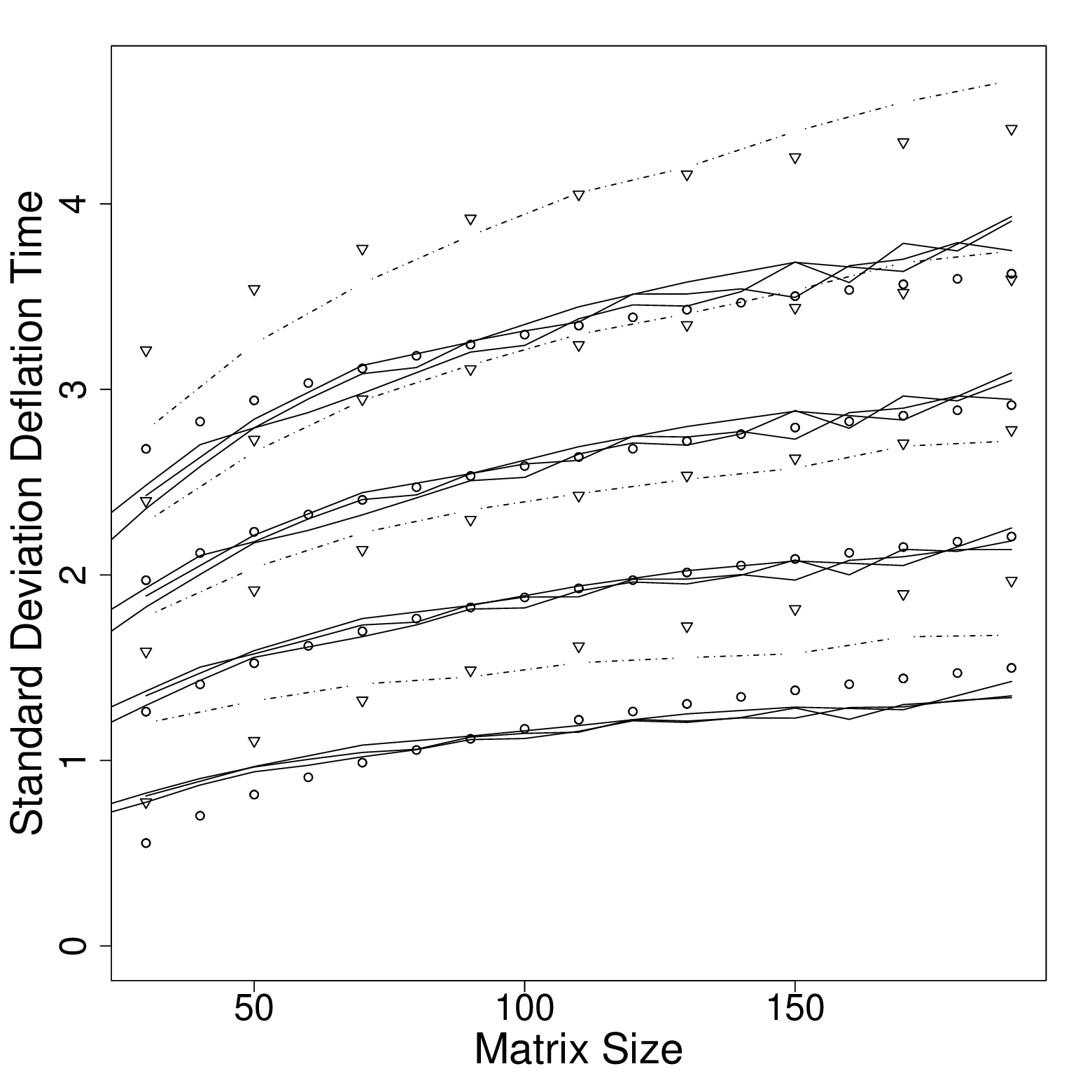}}
\subfloat[][]{
\includegraphics[width=7cm]{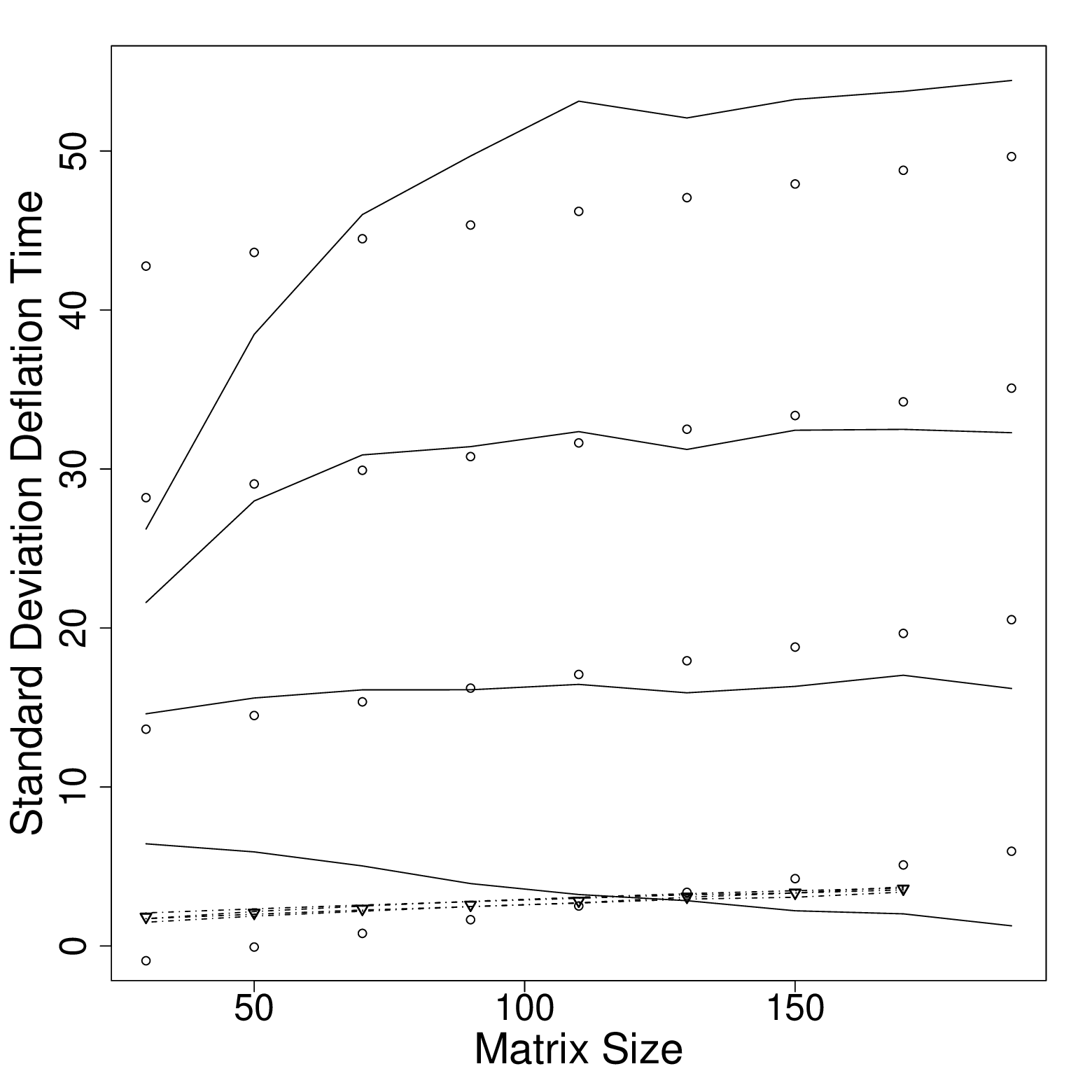}}
\caption{{\bf Standard deviation $\sigma_{n,\eps}$ of the deflation time for the Toda algorithm.\/} (a) Ensembles in the Wigner class. (b) JUE and UDSJ. Legend as in Figure~\ref{fig:QR-means} and regression parameters as in Table~\ref{tab:sdev-regression-toda}.}
\label{fig:Toda-sd}
\end{figure}

The deflation time depends more strongly on $n$ for the Toda algorithm. We explored several regressions but our results for Toda are more ambiguous than for QR. We found that the non-Wigner ensembles (UDSJ and JUE) could be fit with an expression of the form~\qref{eq:qr-mean-regression}--\qref{eq:qr-sdev-regression}. However, the Wigner class ensembles were better suited to the regression
\begin{align} \label{eq:toda-mean-regression-wigner}
\mu_{n,\epsilon}&\approx a_0 + a_1 \log n + a_2 \log\epsilon\\ \label{eq:toda-sdev-regression-wigner}
\sigma_{n,\epsilon}&\approx b_0 + b_1 \log n + b_2 \log\epsilon
\end{align}
The results of this regression are presented in Figure~\ref{fig:Toda-means}, Figure~\ref{fig:Toda-sd}, Table~\ref{tab:mean-regression-toda} and Table~\ref{tab:sdev-regression-toda}.

\begin{table}
\caption{Regression parameters for $\mu_{n,\eps}$ for the Toda algorithm. UDSJ and JUE are fit to \qref{eq:qr-mean-regression}--\qref{eq:qr-sdev-regression} and the Wigner class ensembles are fit to \qref{eq:toda-mean-regression-wigner}--\qref{eq:toda-sdev-regression-wigner}. }
\label{tab:mean-regression-toda}
\begin{center}
\begin{tabular}{|l||l|l|l|}
\hline 
Ensemble & $a_0$ & $a_1$ & $a_2$\\ \hline
GOE, Gaussian Wigner, Bernoulli   & -6.0669 & 1.2888  & -0.7302 \\ \hline
Hermite--1 & -7.0273 &1.6795& -0.7708  \\ \hline
UDSJ & -34.01514 & 0.02984  & -6.60133 \\ \hline
JUE & -2.78614 &0.05318 &-0.74315 \\ \hline
\hline
\end{tabular}
\end{center}
\end{table}
\begin{table}
\caption{Regression parameters for $\sigma_{n,\eps}$ for the Toda algorithm. UDSJ and JUE are fit to \qref{eq:qr-mean-regression}--\qref{eq:qr-sdev-regression} and the Wigner class ensembles are fit to \qref{eq:toda-mean-regression-wigner}--\qref{eq:toda-sdev-regression-wigner}.}
\label{tab:sdev-regression-toda}
\begin{center}
\begin{tabular}{|l||l|l|l|}
\hline 
Ensemble & $b_0$ & $b_1$ & $b_2$\\ \hline
GOE, Gaussian Wigner, Bernoulli & -1.6532  & 0.3347 & -0.1569 \\ \hline
Hermite--1 &-2.1233 &0.6324 &-0.1727  \\ \hline
UDSJ  & -16.46367& 0.04845&  -3.10561 \\ \hline
JUE  &0.97525 &  0.04068 & 0.01451   \\ \hline
\hline
\end{tabular}
\end{center}
\end{table}



\subsection{Deflation index statistics and the effect of the Wilkinson shift}
\begin{figure}[h]
\subfloat[][]{
\includegraphics[width=7cm]{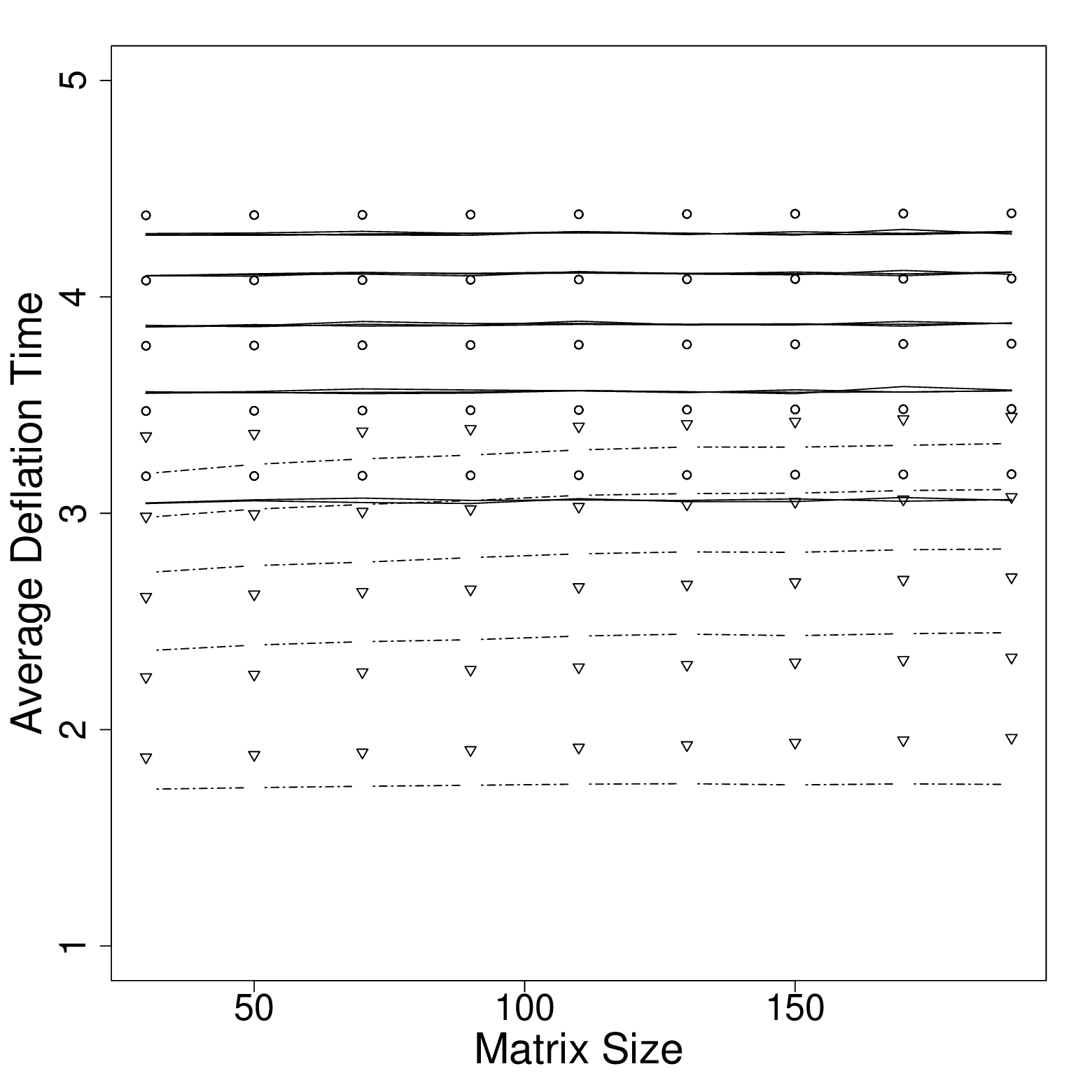}}
\subfloat[][]{
\includegraphics[width=7cm]{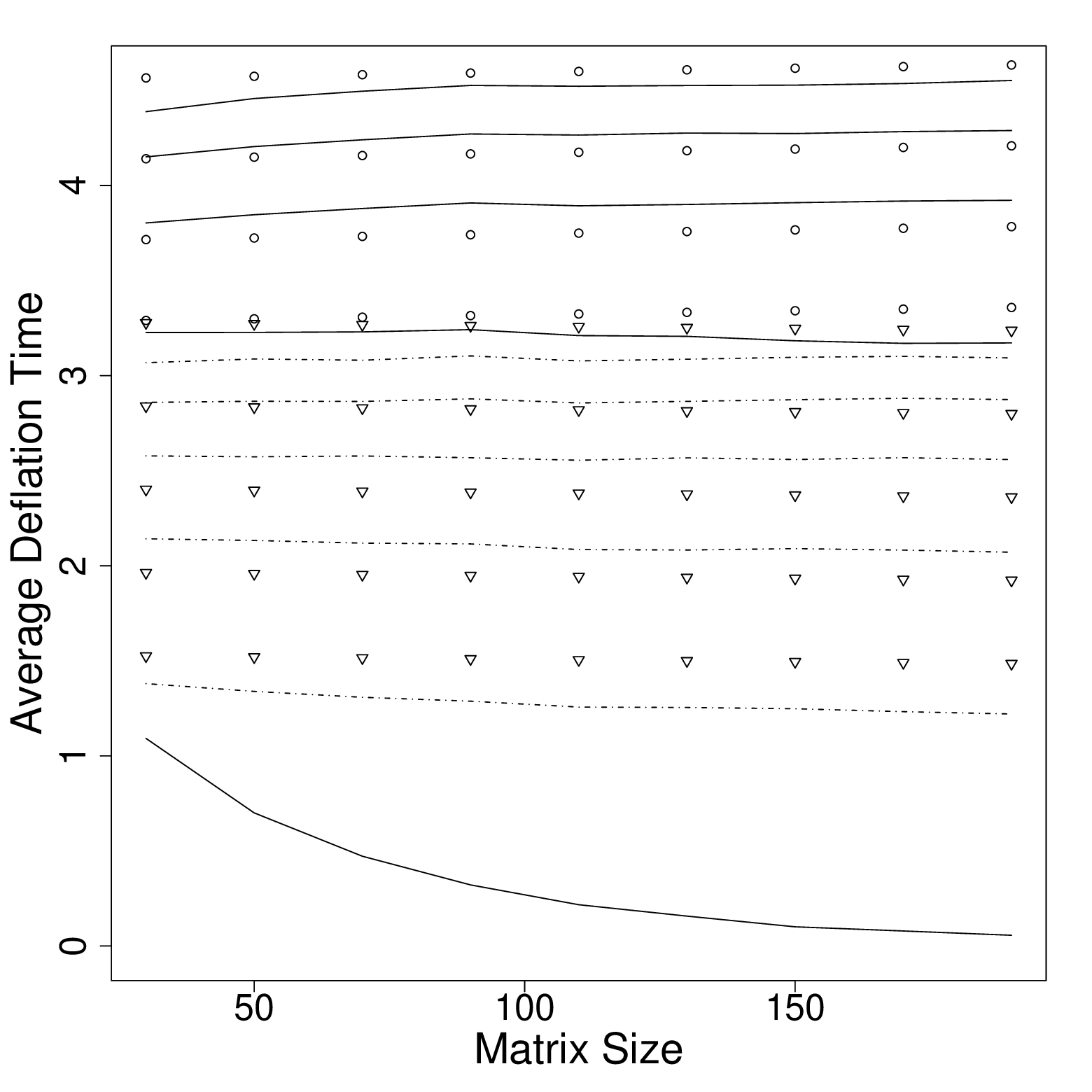}}
\caption{{\bf The effect of the Wilkinson shift.\/} Mean deflation time $\mu_{n,\eps}$ for  the QR algorithm with the Wilkinson shift. (a) Wigner class ensembles. (b) JUE and UDSJ. Empirical data are generated for $\epsilon=10^{-2},...,10^{-8}$ and  $n=20,..., 190$.
The empirical data and a regression of the form \qref{eq:qr-mean-regression} are presented in the same line-styles as Figure \ref{fig:QR-means}. Observe that $\mu_{n,\eps}$ is almost independent of $n$ and that the curves move upwards as $\eps$ decreases as in Figure \ref{fig:QR-means}, but that the scale of the ordinate is different. The regression is not applied to the lowest curve in (b) since $\eps=0.01$ is sufficiently large that several matrices deflate instantaneously. 
}
\label{fig:mean-wilkinson}
\end{figure}

\begin{figure}[h]
\subfloat[][]{
\includegraphics[width=7cm]{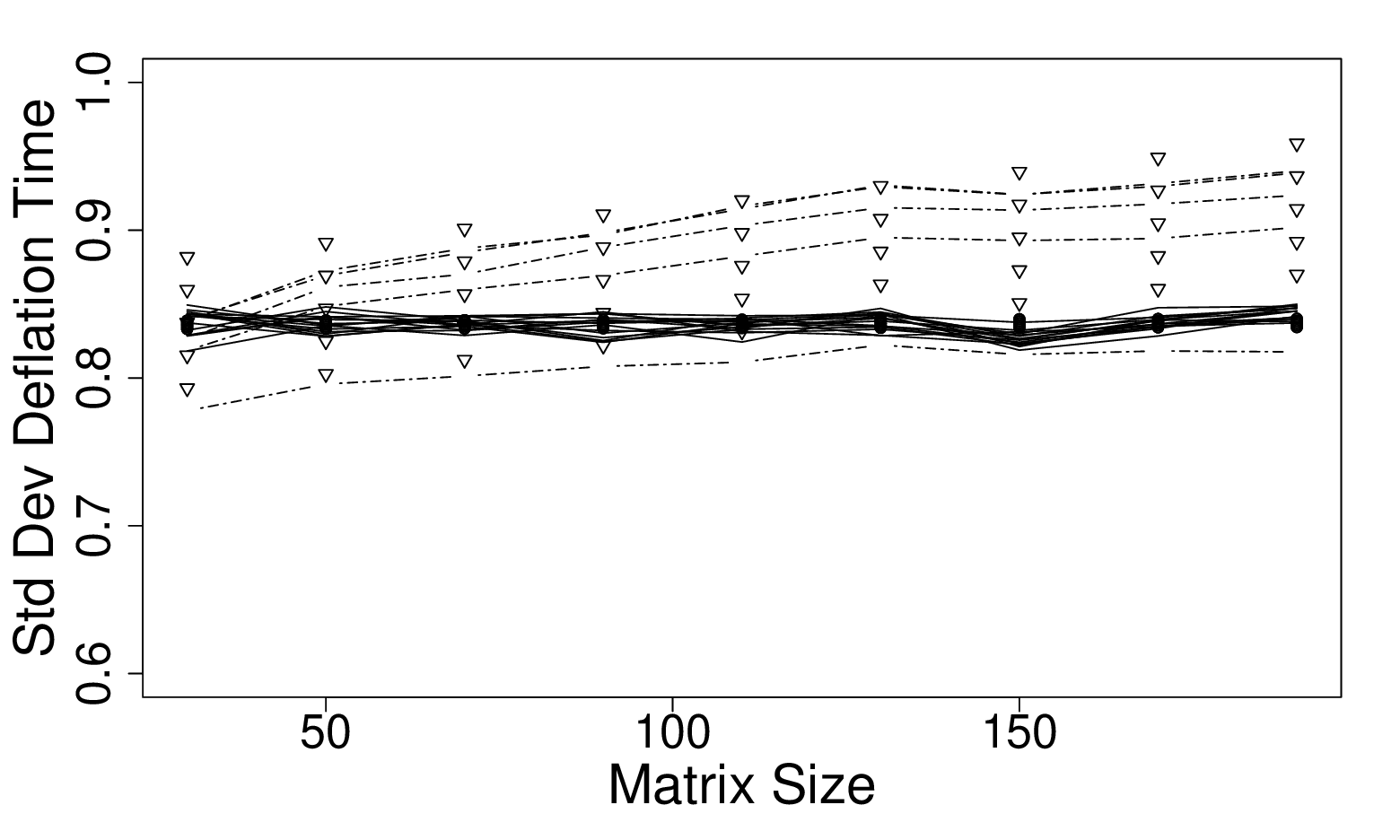}}
\subfloat[][]{
\includegraphics[width=7cm]{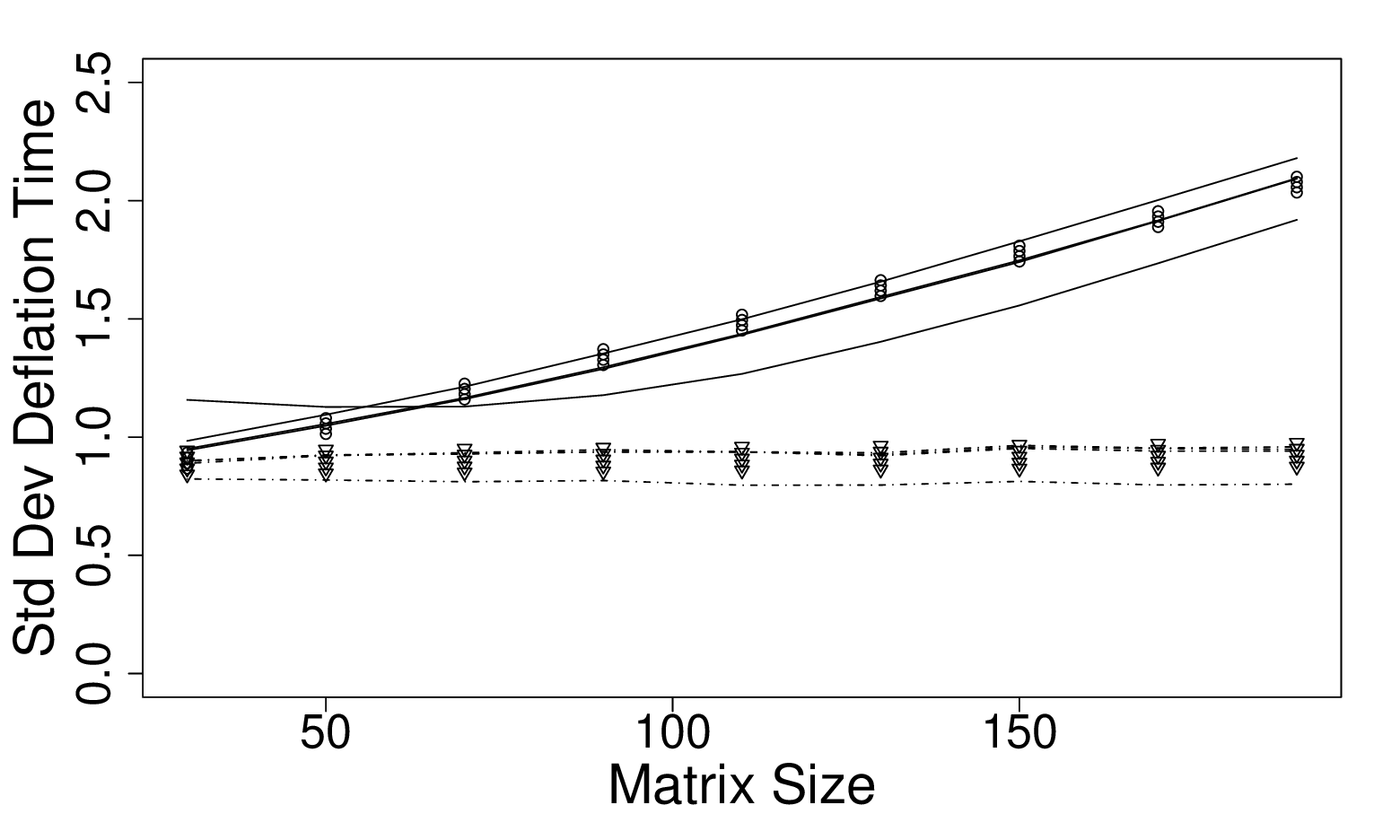}}
\caption{{\bf Standard deviation of deflation time with the Wilkinson shift.\/}
(a) Wigner class ensembles. (b) JUE and UDSJ. 
Line-styles are as in Figure \ref{fig:mean-wilkinson} with a regression of the form \qref{eq:qr-sdev-regression}.
}
\label{fig:sd-wilkinson}
\end{figure}

The remarkable acceleration of QR by shifting is of course well known. Our experiments provide a quantitative statistical picture for the efficacy of the shift. 
Figure~\ref{fig:mean-wilkinson} shows that the deflation time is sharply reduced by the Wilkinson shift. Figure~\ref{fig:sd-wilkinson} shows that the standard deviation of the deflation time is also sharply reduced by the shift.
Deflation takes only a few iterations independent of the size of the matrix.  This is in sharp contrast with the unshifted QR algorithm. 

\begin{figure}[h]
\subfloat[][]{
\includegraphics[width=7cm]{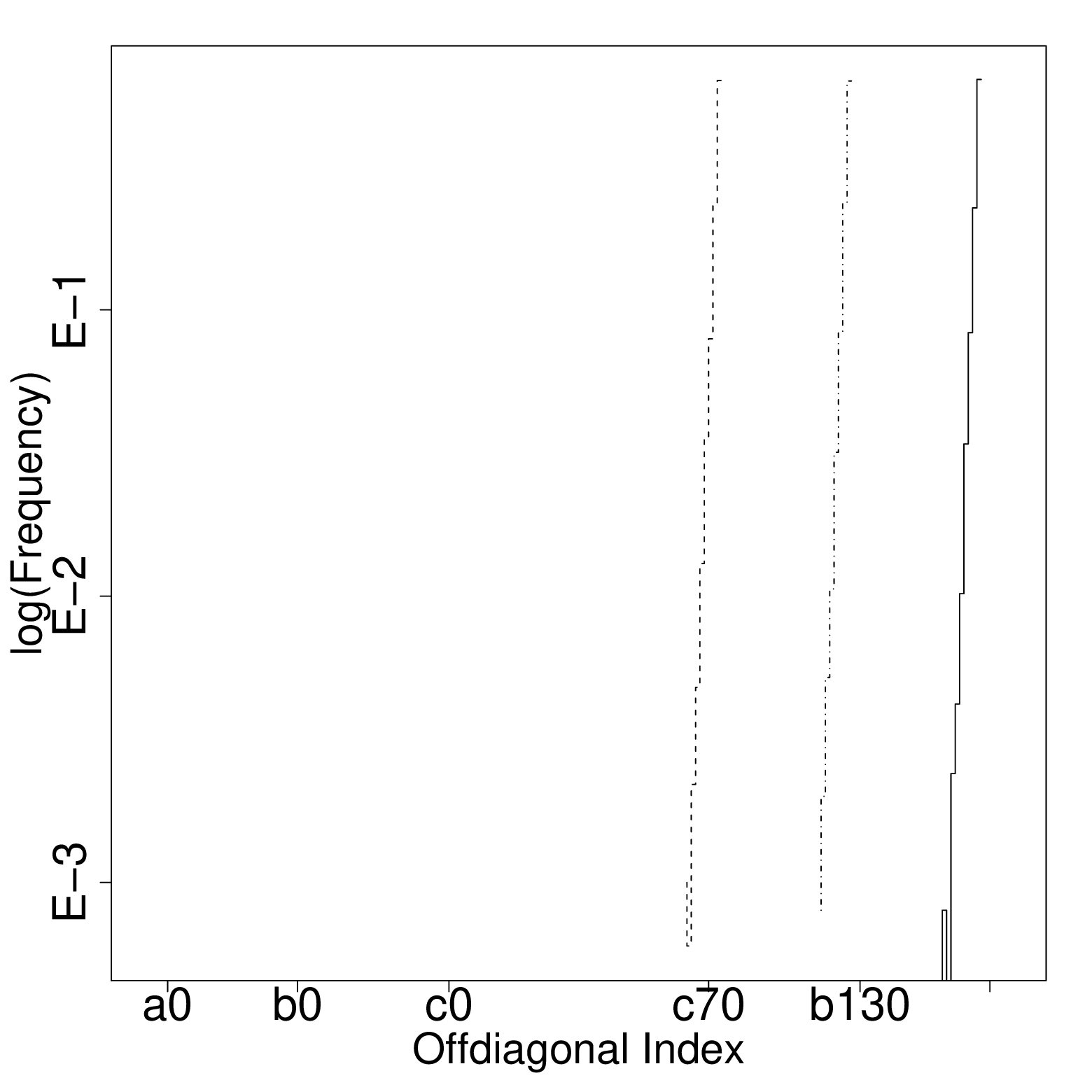}}
\subfloat[][]{
\includegraphics[width=7cm]{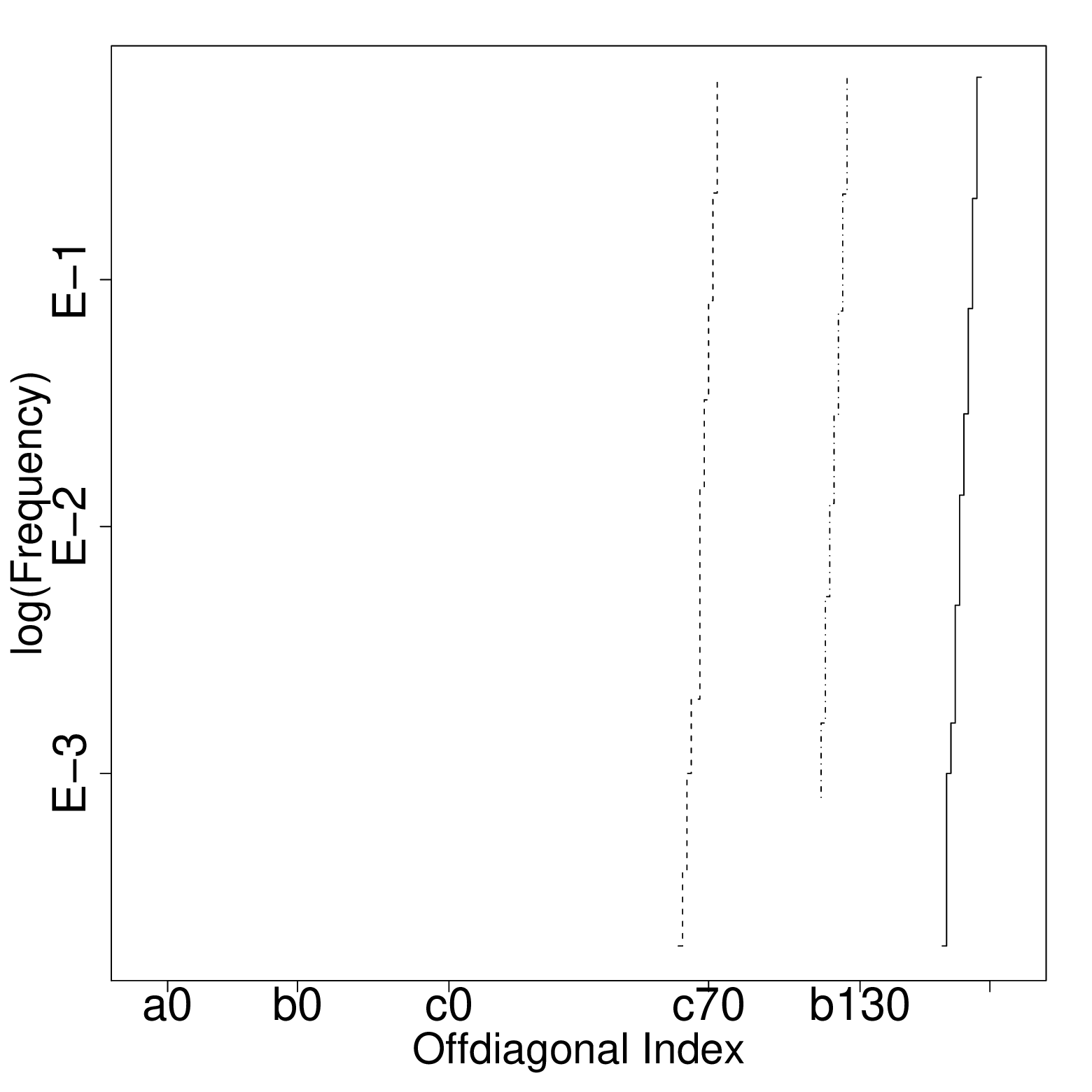}}\\
\subfloat[][]{
\includegraphics[width=7cm]{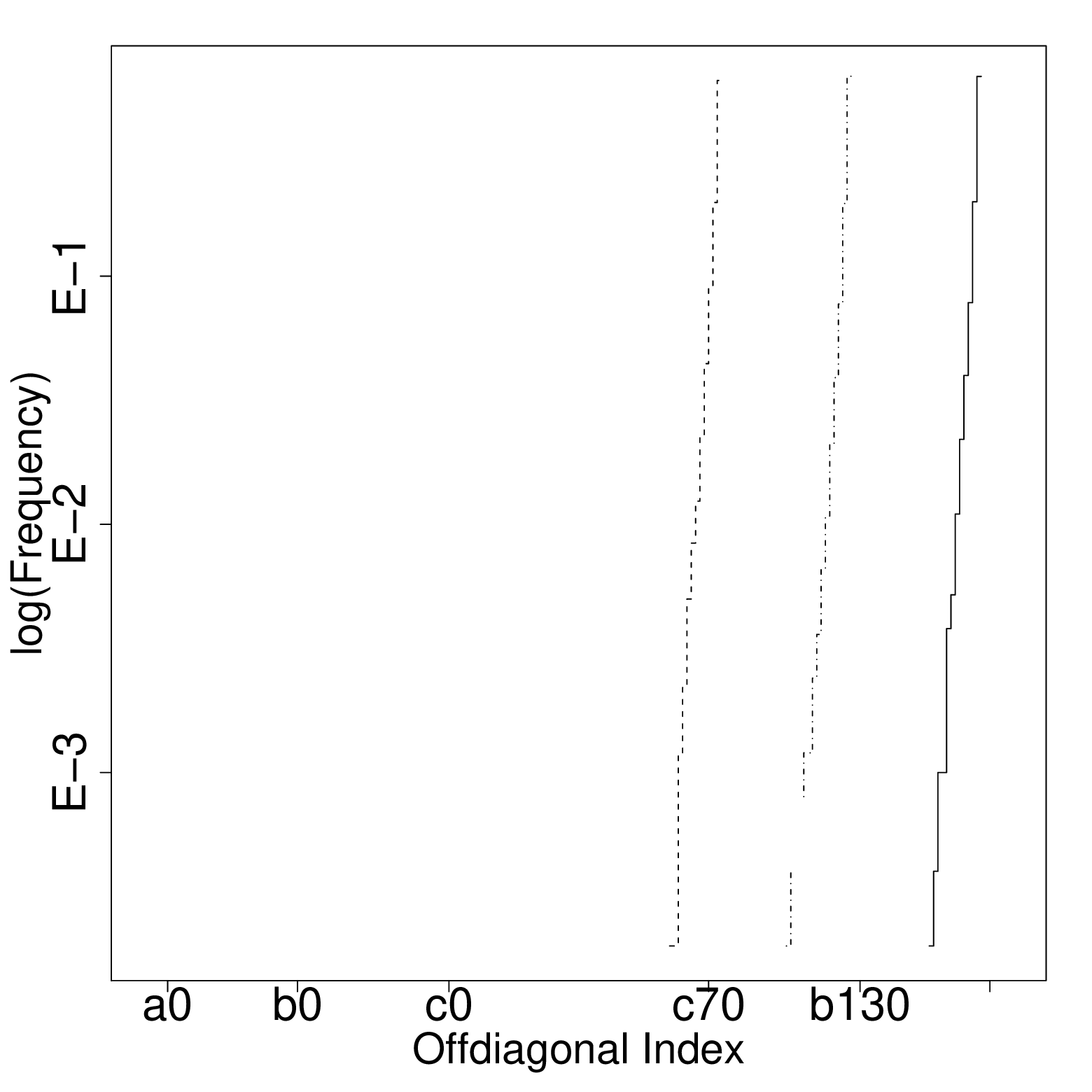}}
\subfloat[][]{
\includegraphics[width=7cm]{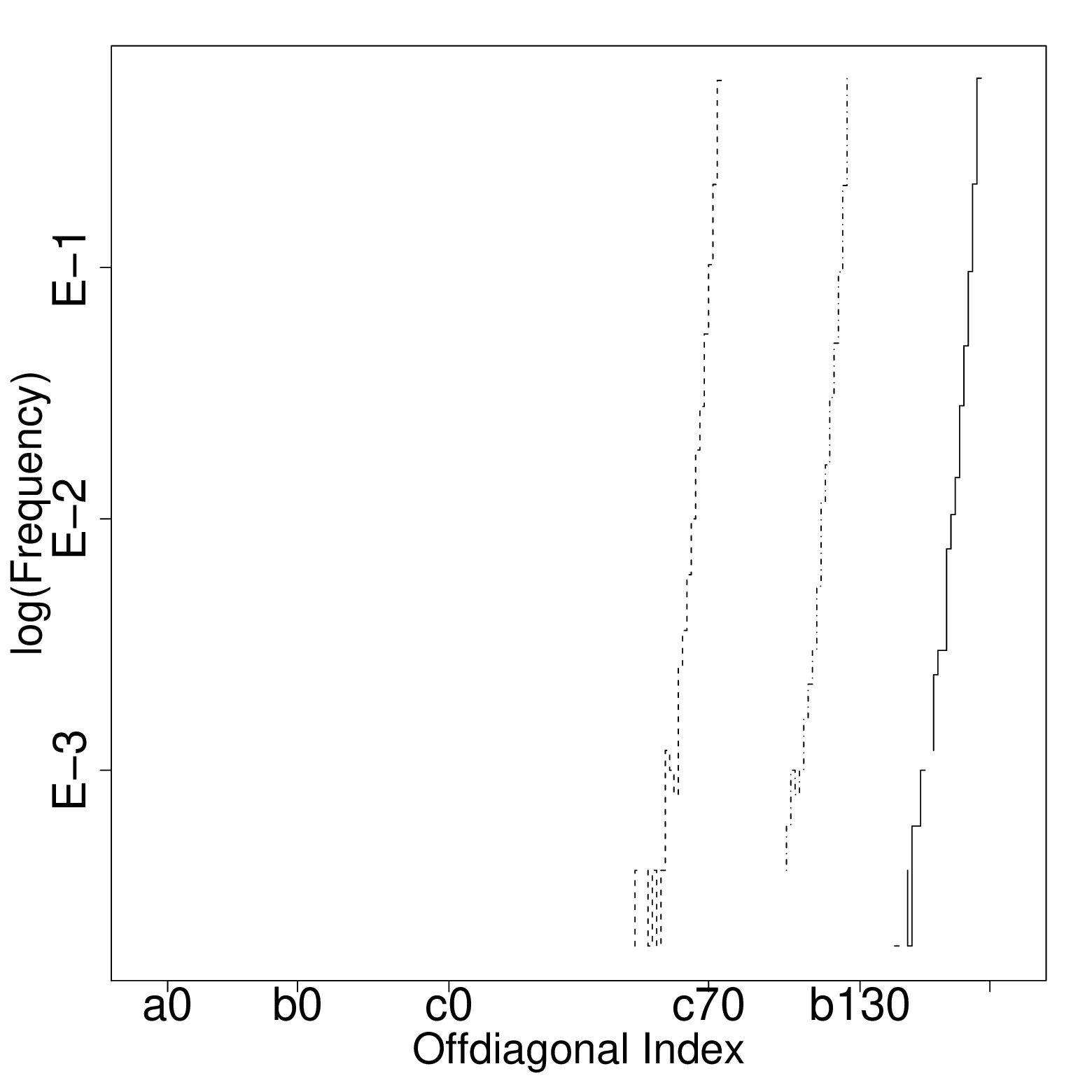}}
\caption{{\bf Empirical distributions of the deflation index $\iota_{n,\epsilon}$ for the unshifted QR algorithm.}
The figures show histograms of the frequency with which deflation occurs at a 
given offdiagonal index. To aid visibility we have centered the distribution so that  the peaks do not overlap. The off-diagonal index takes values between $0$ and $n-2$. Here $a$,$b$ and $c$ refer to ensembles with $n=190$, $130$ and $70$ respectively.
The ensembles shown are (a) Hermite--1; (b) GOE; (c) UDSJ; and (d) JUE.}
\label{fig:QR-deflation-indices}
\end{figure}

\begin{figure}[h]
\subfloat[][]{
\includegraphics[width=7cm]{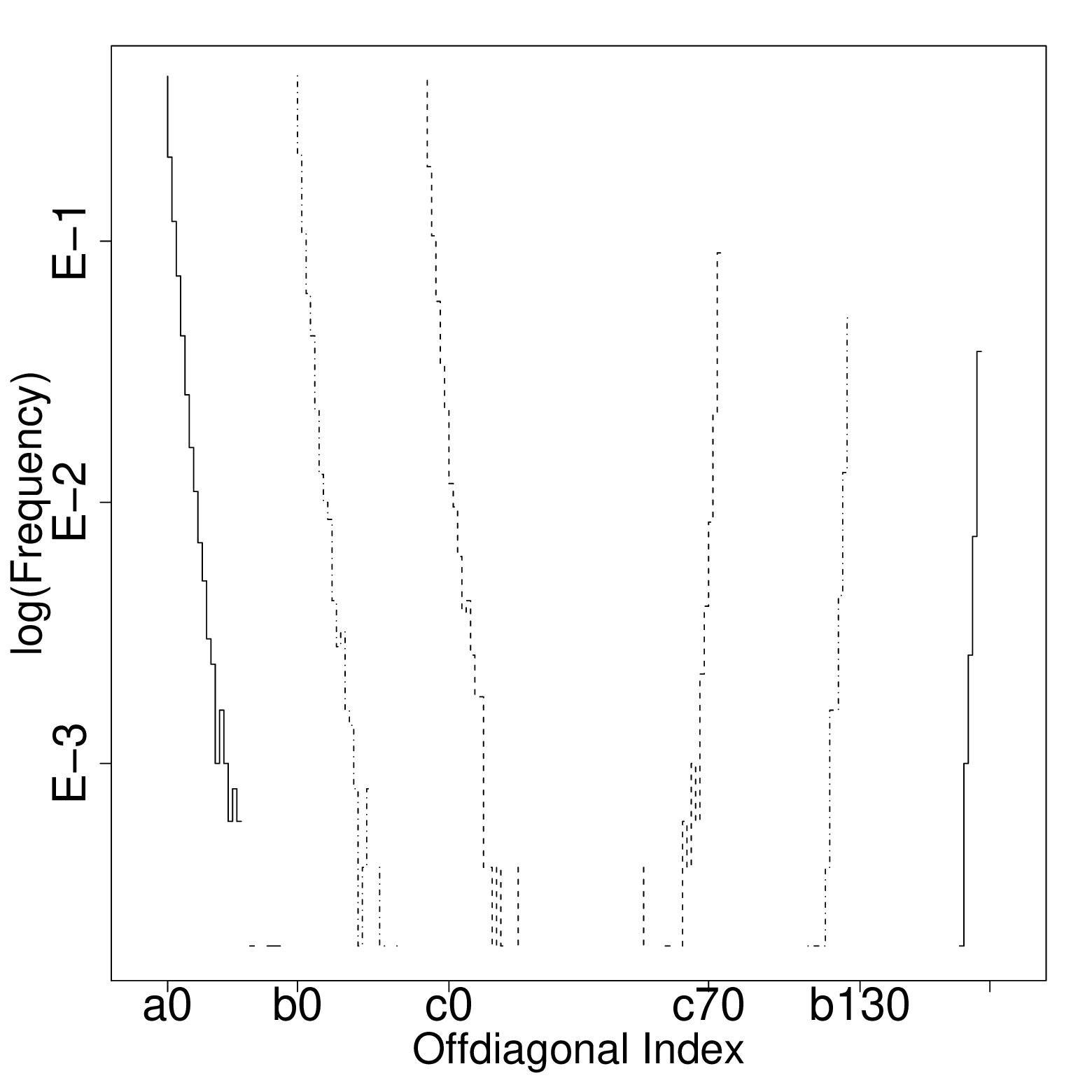}}
\subfloat[][]{
\includegraphics[width=7cm]{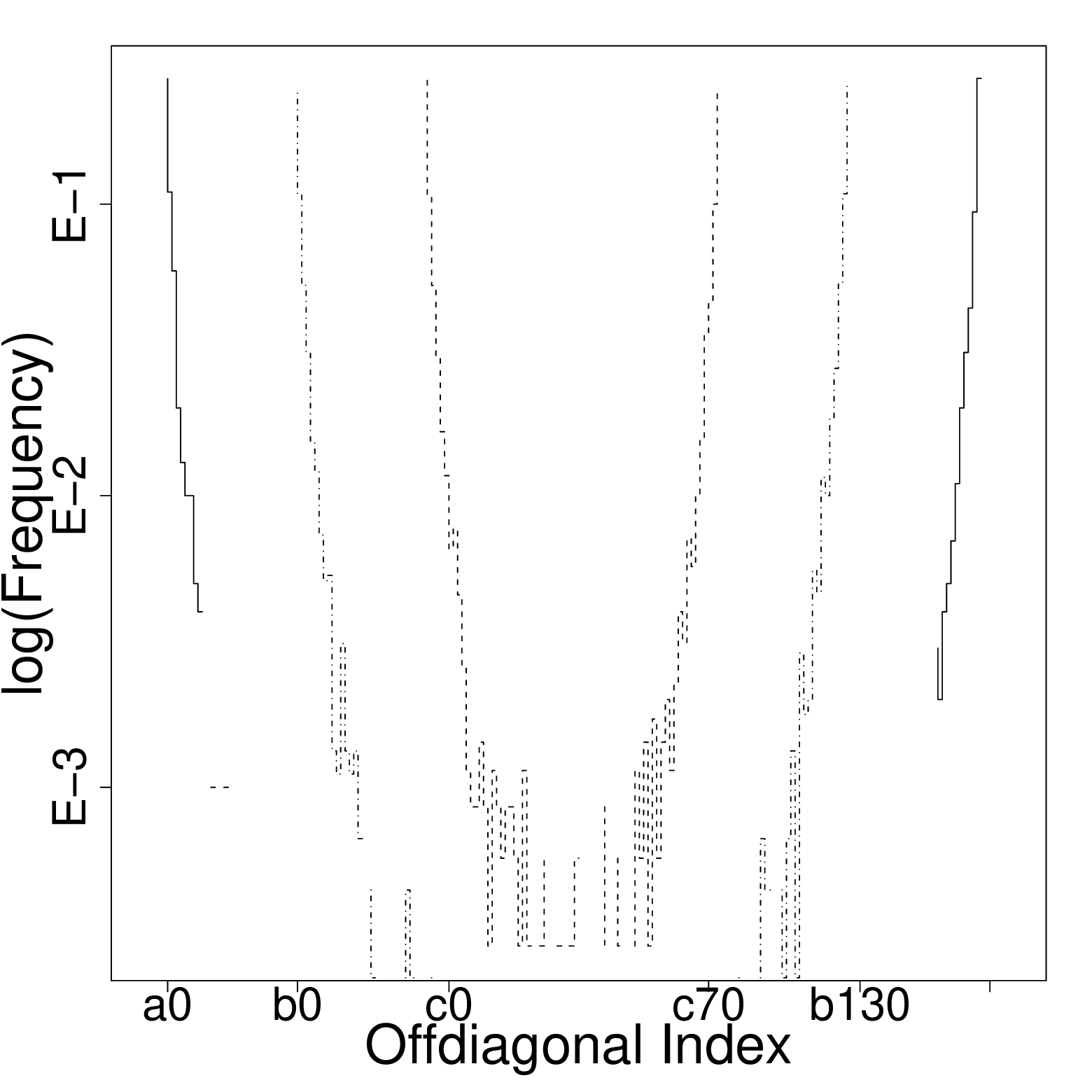}}\\
\subfloat[][]{
\includegraphics[width=7cm]{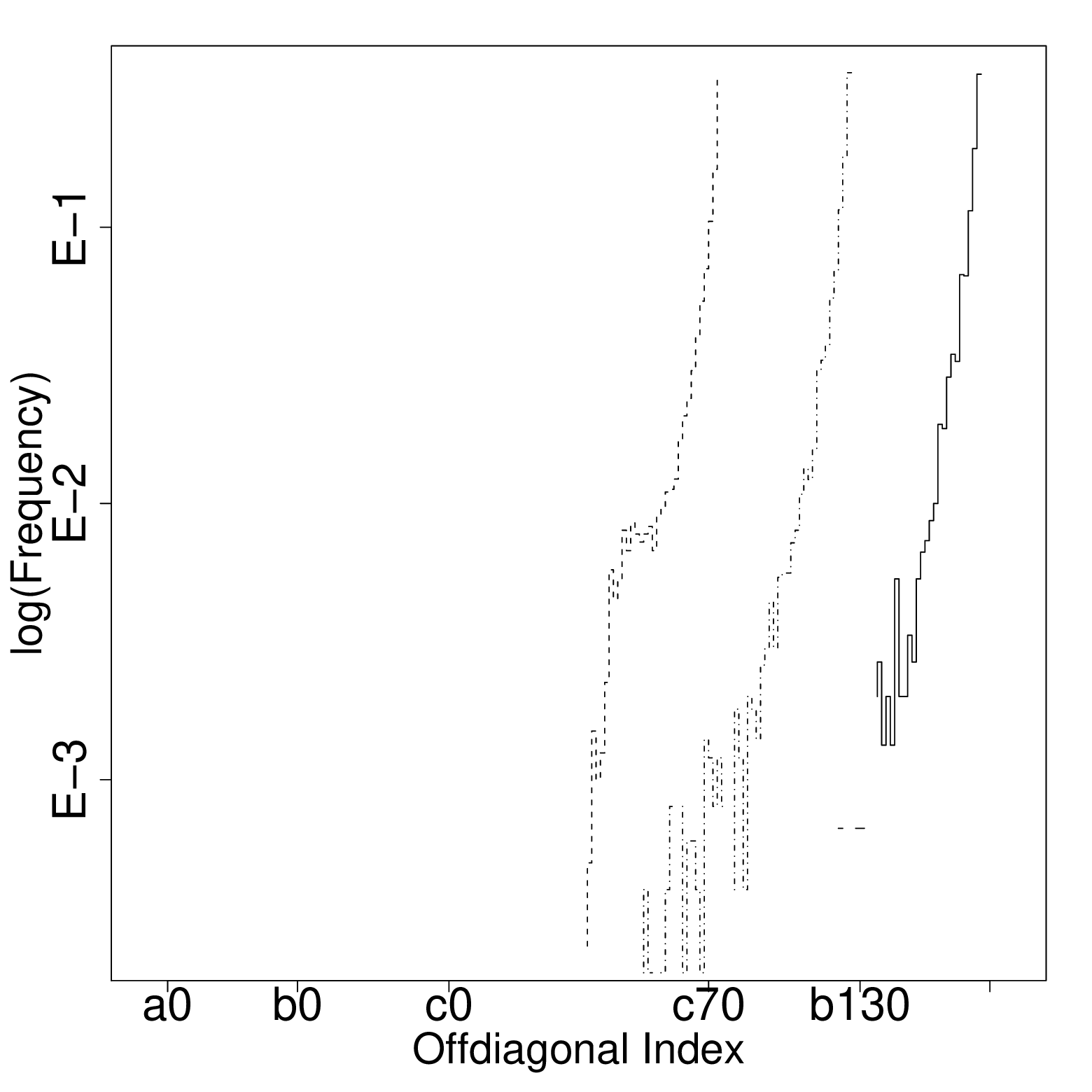}}
\subfloat[][]{
\includegraphics[width=7cm]{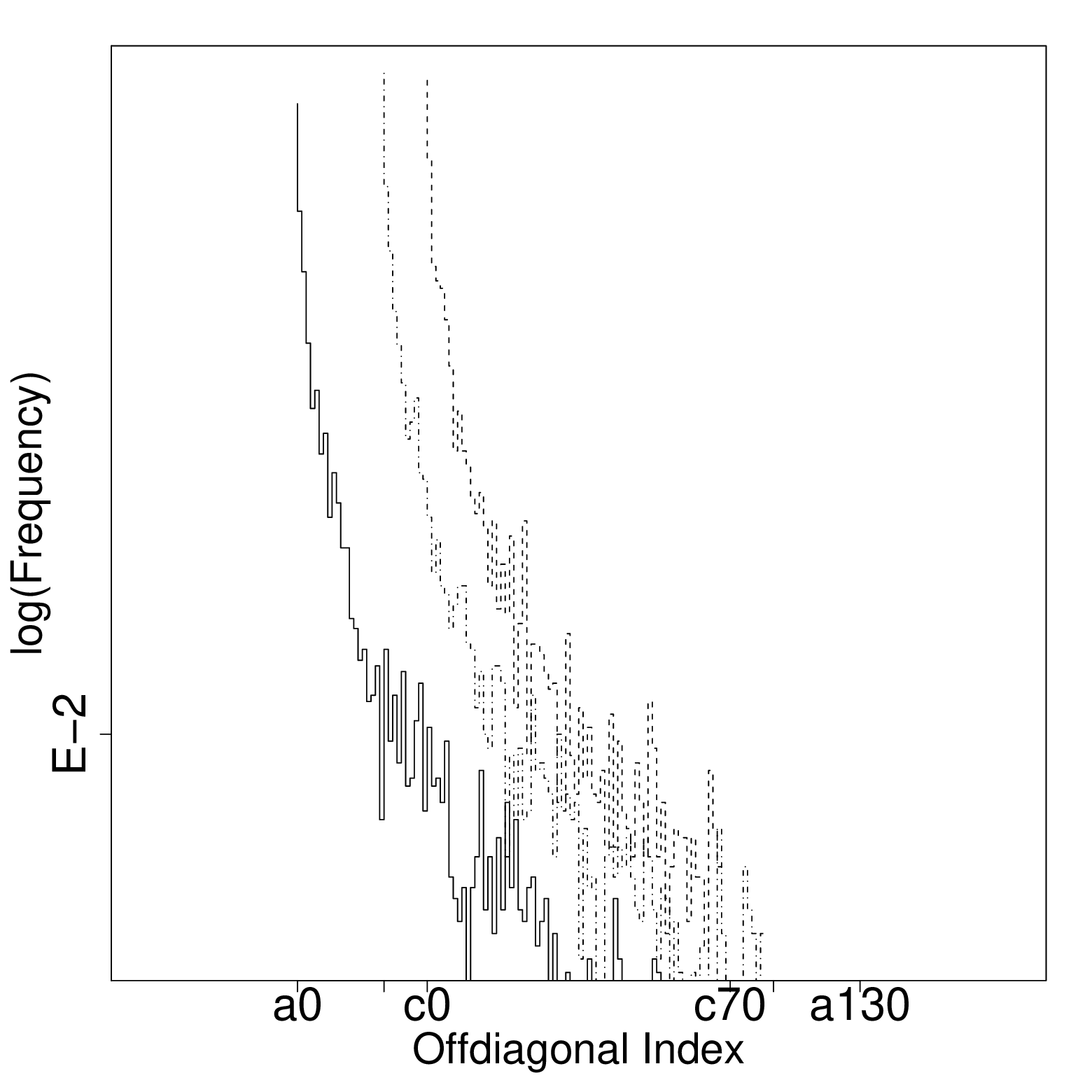}}
\caption{{\bf Empirical distributions of the deflation index $\iota_{n,\epsilon}$ for the Toda algorithm.}
The figures show histograms of the frequency with which deflation occurs at a 
given offdiagonal index. The ensembles are as in Figure~\ref{fig:QR-deflation-indices}.}
\label{fig:Toda-deflation-indices}
\end{figure}

An explanation for the speed-up lies in the statistics of the deflation index shown in Figure~\ref{fig:QR-deflation-indices}. We find that the unshifted QR algorithm deflates at the bottom right corner of the matrix with high probability. Since the Wilkinson shift uses only the $2\times 2$ lower-right block of the matrix, small off-diagonal terms in this block accelerate the unshifted algorithm greatly. 
In contrast with the QR algorithm, the Toda algorithm deflates at both the upper-left and lower-right corner of the matrix (Figure~\ref{fig:Toda-deflation-indices}). Note though that deflation is still predominantly at the corners of the matrix. Similar statistics for other ensembles may be found in~\cite{ChristianDiss}.

\section{Methods and implementation}
\label{sec:methods}
The algorithms were implemented in Python and run on a computing cluster using the module {\tt mpi4py}. For numerical computations we relied
on the {\tt scipy} module except in the case
of the RKPW spectral reconstruction procedure 
\cite{NumStableReconstruction} which was implemented in {\tt C}.
Our simulation strategy was to generate a number $N$ of 
samples for each $(\epsilon,n)$--pair of 
tolerances $\epsilon\in\{10^{-k}:k=2,4,6,8\}$
and matrix dimensions
$n\in\{10,30,\dots,190\}$. One initial matrix sample of
size $n_i\times n_i$ is used to generate deflation time and deflation index 
samples for all pairs $(\tilde{\epsilon},n_i)$, where $\tilde{\epsilon}$
is in our list of tolerances.
To do this we advance the matrix using the algorithm under consideration until
we undercut each of the tolerances in the list and save
the corresponding statistics along the way. Typically for 
each $(\epsilon,n)$--combination we generate between 1000 and 5000 samples.
In the following we present a short summary of the implementation
strategies chosen for the individual algorithms.

\subsection{QR algorithm}
Our simulation code uses the QR decomposition and matrix multiplication methods
provided by {\tt scipy} for the case of full symmetric matrices. For Jacobi
matrices we implemented the efficient (unshifted) QR step presented for example in \cite{GVL}.
We augment these implementations to include the Wilkinson shift by subtracting
(adding) the shift value before (after) the QR step respectively.

\subsection{Toda algorithm}
Both Jacobi and full symmetric matrices are treated similarly for this algorithm.
The implementation uses the QR representation (\ref{eq:dlnt-spectral-explicit})
to generate Toda steps $T_n$ as follows: 
\begin{align}
M_k & = \exp(T_k) = Q_kR_k,\\
M_{k+1} & = R_kQ_k,\\
T_{k+1} & = \log(M_{k+1}).
\end{align}
Our implementation uses {\tt scipy} routines for the matrix
exponentials and matrix logarithms. Note that in general
the matrix exponential of a Jacobi matrix is full symmetric.
{\tt scipy} is also used for the QR decomposition and standard
matrix multiplication routines for the reverse order multiplication.

Note that we do not use an ordinary differential equation solver to solve \qref{eq:dlnt-equation} and diagonalize the matrix as proposed in~\cite{ODESymmetricEig}. This is because our goal here is not to develop a competitive numerical scheme, but to compute reliable statistics of the deflation time for different algorithms. The above numerical scheme based on QR factorization was validated against both an ordinary differential equation solver based method and the use of the explicit solution \qref{eq:dlnt-spectral-explicit} with the \textbf{RKPW} implementation of the inverse spectral map.

\section{Acknowledgments}
The numerical results presented here are part of the first author's Ph.D dissertation at Brown University~\cite{ChristianDiss}.
The help of the support staff at the Center for Computation and Visualization at Brown University is gratefully acknowledged. We 
also thank Jim Demmel, Luen-Chau Li, Irina Nenciu and Nick Trefethen for their interest in this study. 
Finally, we thank the anonymous referee for several valuable comments and suggestions regarding our work.

\bibliographystyle{siam}
\bibliography{qr}

\begin{thebibliography}{10}

\bibitem{Armentano}
{\sc D.~{Armentano}}, {\em {Complexity of path-following methods for the
  eigenvalue problem}}, ArXiv e-prints,  (2011).

\bibitem{Bai-Demmel-Gu}
{\sc Z.~Bai, J.~Demmel, and M.~Gu}, {\em An inverse free parallel spectral
  divide and conquer algorithm for nonsymmetric eigenproblems}, Numerische
  Mathematik, 76 (1997), pp.~279--308.

\bibitem{powerlawEmpiricalData}
{\sc A.~Clauset, C.R. Shalizi, and M.~E.~J. Newman}, {\em Power-law
  distributions in empirical data}, SIAM Review, 51 (2009), pp.~661--703.

\bibitem{OrthPolRH}
{\sc P.~Deift}, {\em {Orthogonal Polynomials and Random Matrices: A Riemann
  Hilbert Approach}}, Courant Lecture Notes in Mathematics, Courant Institute
  of Mathematical Sciences, New York, {First}~ed., 1999.

\bibitem{NotesPercy}
{\sc P.~Deift, C.~D. Levermore, and C.~E. Wayne}, eds., {\em {Dynamical Systems
  and Probabilistic Methods in Partial Differential Equations}}, vol.~31 of
  Lecture Notes in Applied Mathematics, American Mathematical Society,
  Providence, {First}~ed., 1996, ch.~Integrable Hamiltonian Systems,
  pp.~103--138.

\bibitem{TodaFlowGeneric}
{\sc P.~Deift, L.~C. Li, T.~Nanda, and C.~Tomei}, {\em {The Toda Flow on a
  Generic Orbit is Integrable}}, Communications on Pure and Applied
  Mathematics, 39 (1986), pp.~183--232.

\bibitem{Deift-Li-Symplectic}
{\sc P.~Deift, L.-C. Li, and C.~Tomei}, {\em Symplectic aspects of some
  eigenvalue algorithms}, in Important developments in soliton theory, Springer
  Ser. Nonlinear Dynam., Springer, Berlin, 1993, pp.~511--536.

\bibitem{ODESymmetricEig}
{\sc P.~Deift, T.~Nanda, and C.~Tomei}, {\em {Ordinary Differential Equations
  and the Symmetric Eigenvalue Problem}}, SIAM Journal on Numerical Analysis,
  20 (1983), pp.~1--22.

\bibitem{DemmelDifficult}
{\sc J.~W. Demmel}, {\em {The Probability that a Numerical Analysis Problem is
  Difficult}}, {Mathematics of Computation}, 50 (1988), pp.~449--480.

\bibitem{AppliedNumLA}
{\sc James~W. Demmel}, {\em {Applied Numerical Linear Algebra}}, SIAM,
  Philadelphia, {First}~ed., 1997.

\bibitem{RandomDoublyStochasticJacobi}
{\sc P.~{Diaconis} and {P. M.} {Wood}}, {\em {Random Doubly Stochastic Jacobi
  Matrices}}.
\newblock Preprint, 2010.

\bibitem{MatrixModelsBetaEnsembles}
{\sc I.~Dumitriu and A.~Edelman}, {\em {Matrix Models for beta ensembles}},
  Journal of Mathematical Physics, 43 (2002), pp.~5830--5847.

\bibitem{RandMatCondition}
{\sc A.~Edelman}, {\em {Eigenvalues and Condition Numbers of Random Matrices}},
  {Siam J. Matrix Anal. Appl.}, 9 (1988), pp.~543--560.

\bibitem{Edelman-Rao}
{\sc A.~Edelman and N.R. Rao}, {\em Random matrix theory}, Acta Numerica, 14
  (2005), p.~139.

\bibitem{Edelman-Sutton}
{\sc A.~Edelman and B.D. Sutton}, {\em From random matrices to stochastic
  operators}, Journal of Statistical Physics, 127 (2007), pp.~1121--1165.

\bibitem{Yau}
{\sc L{\'a}szl{\'o} Erd{\H{o}}s and Horng-Tzer Yau}, {\em Universality of local
  spectral statistics of random matrices}, Bull. Amer. Math. Soc. (N.S.), 49
  (2012), pp.~377--414.

\bibitem{vonNeumann}
{\sc H.H. Goldstine and J.~Von~Neumann}, {\em Numerical inverting of matrices
  of high order. {II}}, Proceedings of the American Mathematical Society,
  (1951), pp.~188--202.

\bibitem{GVL}
{\sc G.~H. Golub and C.~F. Van~Loan}, {\em Matrix Computations}, The Johns
  Hopkins University Press, Baltimore and London, {Third}~ed., 1996.

\bibitem{NumStableReconstruction}
{\sc W.~B. Gragg and W.~J. Harrod}, {\em {The Numerically Stable Reconstruction
  of Jacobi Matrices from Spectral Data}}, {Numerische Mathematik}, 44 (1984),
  pp.~317--335.

\bibitem{Higham}
{\sc N.J. Higham}, {\em Functions of matrices: theory and computation}, Society
  for Industrial Mathematics, 2008.

\bibitem{Tomei}
{\sc R.S. Leite, N.C. Saldanha, and C.~Tomei}, {\em The asymptotics of
  {W}ilkinson's shift: loss of cubic convergence}, Found. Comput. Math., 10
  (2010), pp.~15--36.

\bibitem{Malyshev}
{\sc A.N. Malyshev}, {\em Parallel algorithm for solving some spectral problems
  of linear algebra}, Linear algebra and its applications, 188 (1993),
  pp.~489--520.

\bibitem{MehtaRMT}
{\sc M.~L. Mehta}, {\em Random Matrices}, Elsevier, San Diego and London,
  third~ed., 2004.

\bibitem{MoserNotesDS}
{\sc J.~Moser and E.~J. Zehnder}, {\em {Notes on Dynamical Systems}}, Courant
  Lecture Notes in Mathematics, American Mathematical Society, Providence,
  {First}~ed., 2005.

\bibitem{DiffEqQRAlgo}
{\sc T.~Nanda}, {\em {Differential Equations and the QR Algorithm}}, SIAM
  Journal on Numerical Analysis, 22 (1985), pp.~310--321.

\bibitem{ChristianDiss}
{\sc C.~W. Pfrang}, {\em Diagonalizing Random Matrices with Integrable
  Systems}, PhD thesis, Brown University, 2011.

\bibitem{Rudelson-Vershynin}
{\sc M.~Rudelson and R.~Vershynin}, {\em The {L}ittlewood-{O}fford problem and
  invertibility of random matrices}, Advances in Mathematics, 218 (2008),
  pp.~600--633.

\bibitem{SmoothedAnalysisCondition}
{\sc A.~Sankar, D.~A. Spielman, and S.-H. Teng}, {\em {Smoothed Complexity of
  the Condition Numbers and Growth Factors of Matrices}}, {Siam J. Matrix Anal.
  Appl.}, 28 (2006), pp.~446--476.

\bibitem{Smale}
{\sc S.~Smale}, {\em On the average number of steps of the simplex method of
  linear programming}, Mathematical Programming, 27 (1983), pp.~241--262.

\bibitem{HamiltonianGroup}
{\sc W.~W. Symes}, {\em Hamiltonian group actions and integrable systems},
  Phys. D, 1 (1980), pp.~339--374.

\bibitem{QRAlgoScattering}
\leavevmode\vrule height 2pt depth -1.6pt width 23pt, {\em The {$QR$} algorithm
  and scattering for the finite nonperiodic {T}oda lattice}, Phys. D, 4
  (1981/82), pp.~275--280.

\bibitem{Tao-Vu-2010}
{\sc T.~Tao and V.~Vu}, {\em Random matrices: The distribution of the smallest
  singular values}, Geometric And Functional Analysis, 20 (2010), pp.~260--297.

\bibitem{Tracy-Widom}
{\sc C.A. Tracy and H.~Widom}, {\em Level-spacing distributions and the airy
  kernel}, Communications in Mathematical Physics, 159 (1994), pp.~151--174.

\bibitem{Trefethen}
{\sc L.N. Trefethen and D.~Bau}, {\em Numerical linear algebra}, Society for
  Industrial Mathematics, 1997.

\bibitem{IsospectralFlows}
{\sc D.~S. Watkins}, {\em {Isospectral Flows}}, {SIAM Review}, 26 (1984),
  pp.~379--391.

\bibitem{Wilkinson}
{\sc J.~H. Wilkinson}, {\em Global convergence of tridiagonal {${\rm QR}$}
  algorithm with origin shifts}, Linear Algebra and Appl., 1 (1968),
  pp.~409--420.

\end{thebibliography}
\end{document}